\documentclass[aap]{imsart}

\RequirePackage{amsthm,amsmath,amsfonts,amssymb}
\RequirePackage[numbers]{natbib}
\RequirePackage[colorlinks,citecolor=blue,urlcolor=blue]{hyperref}
\RequirePackage{graphicx}
\RequirePackage{enumitem}
\RequirePackage{tikz}

\startlocaldefs
\makeatletter\def\journal@name{}\makeatother
\numberwithin{equation}{section}
\theoremstyle{plain}
\newtheorem{theorem}{Theorem}[section]
\newtheorem{proposition}[theorem]{Proposition}
\newtheorem{lemma}[theorem]{Lemma}
\newtheorem{corollary}[theorem]{Corollary}
\theoremstyle{definition}
\newtheorem{definition}[theorem]{Definition}

\newtheorem{remark}[theorem]{Remark}

\newcommand{\R}{\mathbb{R}}

\newcommand{\E}{\mathbb{E}}
\renewcommand{\P}{\mathbb{P}}
\newcommand{\ip}[2]{\langle #1,#2\rangle}
\newcommand{\cZ}{\mathcal{Z}}
\newcommand{\cR}{\mathcal{R}}
\newcommand{\cH}{\mathcal{H}}

\newcommand{\Wt}{\mathsf{W}_2}
\newcommand{\TV}{d_{\mathrm{TV}}}
\newcommand{\eps}{\varepsilon}
\newcommand{\one}{\mathbf{1}}
\DeclareMathOperator{\haus}{haus}

\DeclareMathOperator{\tr}{tr}
\DeclareMathOperator{\Law}{Law}
\DeclareMathOperator{\diam}{diam}
\DeclareMathOperator{\aff}{aff}
\DeclareMathOperator{\Leb}{Leb}
\DeclareMathOperator{\dist}{dist}
\DeclareMathOperator{\osc}{osc}
\endlocaldefs

\begin{document}

\begin{frontmatter}
\title{Metastability for diffusions driven by metric projections onto disjoint convex sets}
\runtitle{Metastability for projection diffusions}

\begin{aug}
\author[A]{\fnms{Narayana}~\snm{Darapaneni}\ead[label=e1]{narayana.darapaneni@northwestern.edu}}
\address[A]{Northwestern University SPS\printead[presep={,\ }]{e1}}
\runauthor{N. Darapaneni}
\end{aug}

\begin{abstract}
Let $C$ be a finite union of disjoint compact convex sets in $\R^d$. We study the diffusion whose drift pulls the state towards its nearest point in $C$, with noise intensity $T$, as $T\to0$ and as time grows. Let $\delta_i$ be the distance from $C_i$ to the nearest other component. The exit time from the basin of $C_i$ is $\exp((\delta_i^2/8+o(1))/T)$ in probability. Up to times $\exp(c/T)$ with $c<\min_i\delta_i^2/8$, the basin weights are frozen and inside each basin the state behaves as for a single convex set. After a time $\exp(c/T)$, with an explicit $c$, the law of the state is close to the Gibbs state in total variation. For two and three components any $c$ larger than the largest $\delta_i^2/8$ works. For two components this is sharp. At zero temperature the Gibbs state charges only the components of largest dimension, in proportion to their volumes. In the associated control problem, the two limits of the value function commute for every bounded Lipschitz terminal cost exactly when $C$ is convex. The proofs use results of a companion paper.
\end{abstract}



\end{frontmatter}

\tableofcontents

\section{Introduction}\label{sec:intro}

\subsection{The model}
Let $C$ be a nonempty closed subset of $\R^d$. We write
\[
d_C(x)=\min_{c\in C}|x-c| ,\qquad V_C(x)=\tfrac12\,d_C(x)^2 ,
\]
and we denote by $P_Cx$ a nearest point of $x$ in $C$. We study the diffusion
\begin{equation}\label{eq:SDEintro}
dX_t=-(X_t-P_CX_t)\,dt+\sqrt{2T}\,dW_t ,
\end{equation}
where $W$ is a standard Brownian motion in $\R^d$ and $T>0$ is the intensity of the noise. We call $T$ the temperature and $X_t$ the state. The drift moves each state along the straight segment to its nearest point in $C$. At every point $x$ where the nearest point is unique we have $\nabla V_C(x)=x-P_Cx$. So \eqref{eq:SDEintro} is a gradient diffusion with potential $V_C$. The law $\rho_{T,t}$ of $X_t$ solves the forward Kolmogorov equation
\begin{equation}\label{eq:FP}
\partial_t\rho=\nabla\cdot(\rho\,\nabla V_C)+T\,\Delta\rho .
\end{equation}
When $C$ is compact, the invariant law is the Gibbs state
\[
\gamma_T^C(dx)=\cZ_T(C)^{-1}e^{-d_C(x)^2/2T}\,dx ,
\]
where $\cZ_T(C)$ is the normalising constant. Throughout the paper we also use one more process. We denote by $\zeta$ the process of Ornstein and Uhlenbeck (OU process for short) which solves $d\zeta=-\zeta\,dt+\sqrt{2T}\,dW$, $\zeta_0=0$. It is driven by the same Brownian motion as the state.

The diffusion \eqref{eq:SDEintro} appears in several places. After the change of time $s=tT$ it is the penalised equation used to construct reflected diffusions in convex domains \cite{Menaldi1983,LionsSznitman1984}. The function $d_C^2/2T$ is the Moreau envelope, with parameter $T$, of the function which is $0$ on $C$ and $+\infty$ outside $C$. For this reason \eqref{eq:SDEintro}, in the time $s$, is the continuous time limit of the proximal Langevin algorithm of \cite{Brosse2017}. One can also read \eqref{eq:SDEintro} as a model for the internal state of a system perturbed by noise. The system is pulled towards the set $C$ of states that are feasible in a given context. Models with a latent state which evolves in time are common in machine learning, see for example \cite{HaSchmidhuber2018,Song2021}. In this reading a disconnected set $C$ describes several incompatible groups of feasible states.

In this paper $C$ is a finite union of pairwise disjoint sets,
\[
C=C_1\cup\dots\cup C_m,\qquad m\ge2,
\]
where each $C_i$ is nonempty, compact and convex. The potential $V_C$ vanishes on every component, so it has $m$ wells and these wells are flat. Its gradient jumps across the set of points which are equally close to two components. The basin $B_i$ of $C_i$ is the set of points which are nearer to $C_i$ than to every other component, that is, the Voronoi cell of $C_i$. The ridge is the set of points which lie in no basin. On the ridge the nearest component is not unique. Let $\delta_i$ be the distance from $C_i$ to the nearest other component. The lowest value of $V_C$ on the boundary of $B_i$ is $\delta_i^2/8$. Metastability means that, at low temperature, the state stays in one basin for a very long time before it passes to another basin.

There are two parameters, the time $t$ and the temperature $T$. We ask how the law of $X_t$ behaves when $t\to\infty$ and $T\to0$ together. The companion paper \cite{DarapaneniConvex} gives the answer for a single compact convex set $C$. Let the initial law $\rho_0$ have a finite second moment. When $T\to0$ and $t\to\infty$ with $tT\to s\in[0,\infty]$, the law $\rho_{T,t}$ converges to
\begin{equation}\label{eq:introdiagram}
\mu_s=\begin{cases}(P_C)_\#\rho_0,& s=0,\\[2pt] \Law(R_s),&0<s<\infty,\\[2pt] u_C,& s=\infty,\end{cases}
\end{equation}
where $u_C$ is the uniform probability on $C$ and $R$ is the reflected Brownian motion in $C$ started from the law $(P_C)_\#\rho_0$, the image of $\rho_0$ under $P_C$ (Theorem \ref{thm:diagram} below). In a time of order one the state comes close to $C$. It needs a time of order $1/T$ to explore $C$. For a disconnected set there is one more time scale, the time which the state needs to pass from one component to another. This time is exponentially long in $1/T$.

\subsection{Results}
We prove the following.
\begin{enumerate}[label=\textup{(\roman*)}]
\item \emph{Exit times} (Theorem \ref{thm:exit}). Let the state start in $B_i$. We write $\tau_i$ for the exit time from $B_i$ and $\theta_i$ for the first time at which the state reaches a neighbourhood of another component. Then
\[
T\log\tau_i\longrightarrow\frac{\delta_i^2}8,\qquad T\log\theta_i\longrightarrow\frac{\delta_i^2}8\qquad\text{in probability as }T\to0 .
\]
\item \emph{Freezing} (Theorem \ref{thm:meta}). Suppose that $\rho_0$ has a finite second moment and does not charge the ridge. Let $T\to0$ and $t\to\infty$ with $tT\to s\in[0,\infty]$ and $t\le\exp(c/T)$, where $c<\min_i\delta_i^2/8$. Then the law of $X_t$ converges weakly to a limit in which each basin keeps its initial mass $a_i=\rho_0(B_i)$. Inside each basin this limit is given by the picture \eqref{eq:introdiagram} for the component $C_i$, started from $\rho_0$ conditioned on the basin. In particular the limit is $\sum_ia_iu_{C_i}$ when $tT\to\infty$.
\item \emph{Equilibrium} (Proposition \ref{prop:gibbsunion}). As $T\to0$, $\gamma_T^C\to\sum_iw_iu_{C_i}$ weakly. Only the components of largest dimension have a positive weight $w_i$. These weights are proportional to their volumes.
\item \emph{Convergence to equilibrium} (Theorem \ref{thm:equil}). Let $\Delta_*$ be the smallest number such that $\{d_C\le\Delta_*/2\}$ is connected (this agrees with Definition \ref{def:direct}, see Remark \ref{rem:connected}). If $t\ge\exp(c/T)$ with $c>\lfloor m/2\rfloor\Delta_*^2/8$, then $\rho_{T,t}-\gamma_T^C\to0$ in total variation, for every initial law.
\end{enumerate}
In general the limits in (ii) and (iii) are different. They agree exactly when $a_i=w_i$ for all $i$ (Corollary \ref{cor:three}). For two and for three components we have $\lfloor m/2\rfloor=1$ and $\Delta_*=\max_i\delta_i$. The threshold $\Delta_*^2/8$ in (iv) then cannot be lowered as soon as a component $C_i$ with $\delta_i=\Delta_*$ has $w_i<1$ (Remark \ref{rem:equil}). For two components this always happens. Consider now two components and put $\delta=\dist(C_1,C_2)$. Then (ii) and (iv) describe the law of $X_t$ except when $T\log t$ is close to $\delta^2/8$. Under the assumptions of (ii), the law is close to $a_1u_{C_1}+a_2u_{C_2}$ for $1/T\ll t\le e^{c/T}$, $c<\delta^2/8$. It is close to $w_1u_{C_1}+w_2u_{C_2}$ for $t\ge e^{c/T}$, $c>\delta^2/8$ (Corollary \ref{cor:two}).

\emph{The control problem.} Suppose that an agent can add a drift $a$ to \eqref{eq:SDEintro} at the running cost $\frac12|a|^2$ and pays a terminal cost $g$ at the time $t$. Let $v_T(t,x)$ be the least expected cost of steering the state from $x$ over the horizon $t$. By the logarithmic transformation, $v_T(t,x)=-2T\log\E\exp(-g(X^x_t)/2T)$, where $X^x$ is the solution of \eqref{eq:SDEintro} started at $x$. We take this formula as the definition of $v_T$. For a compact convex set $C$ and every bounded Lipschitz function $g$, the companion paper \cite{DarapaneniConvex} shows that
\[
\lim_{T\to0}\lim_{t\to\infty}v_T(t,x)=\lim_{t\to\infty}\lim_{T\to0}v_T(t,x)=\min_{y\in\R^d}\{g(y)+d_C(y)^2\}
\]
(Theorem \ref{thm:hjb} below). So for a convex set the two limits commute for the value function, although they do not commute in general for the law of the state. For a disjoint union this fails for some $g$. Theorem \ref{thm:hjbdisc} gives a lower bound for $v_T(t,x)$ as $T\to0$ when $x\in C_i$. For suitable $g$ it shows that the two limits of $v_T$ differ by $\delta_i^2/4$. Hence, among finite disjoint unions of compact convex sets, the two limits commute for every bounded Lipschitz $g$ and every starting point in $C$ exactly when the set is convex, that is, when $m=1$ (Corollary \ref{cor:hjbconvex}).

\subsection{Methods and relation to earlier work}
For smooth potentials, exit times and transitions between wells are described by the theory of \cite{FreidlinWentzell2012}, which also allows wells that are compact sets. The potential theoretic approach \cite{BEGK2004,BovierdenHollander2015} gives in addition sharp prefactors. See also \cite{OlivieriVares2005,HKS1989,MenzSchlichting2014}. Degenerate saddles and wells are treated in \cite{BerglundGentz2010,AvelinJulinViitasaari2023}. Our situation is different. Here the wells are flat sets of different dimensions, the potential is not differentiable on the ridge and the basins may be unbounded with boundaries which are not smooth. The theory of \cite{FreidlinWentzell2012} is written for Lipschitz coefficients and for domains with smooth boundary, so it does not apply directly.

We do not use large deviation estimates on path space for the results on metastability. The lower bound for exit times comes from the envelope. We obtain the upper bound as follows. First we use the reflected Brownian motion limit \eqref{eq:introdiagram} to bring the state near a point of its component which is nearest to another component. Then we bound by the Girsanov theorem the probability that the state follows one explicit path to the ridge. At the ridge we use a submartingale argument. A renewal argument concludes the proof. For the convergence to equilibrium we need in addition a Doeblin argument. We obtain the exponential rates but not the prefactors.

This paper is a companion of \cite{DarapaneniConvex}, which treats a single convex set. We use several results of that paper. The envelope bounds the distance of the state from a convex set, path by path, by its initial value and the norm of the OU process. A comparison in law gives a weaker bound for every closed set. We also use a maximal inequality for the OU process, the limit diagram \eqref{eq:introdiagram}, a lower bound for the transition density at the scale $\sqrt T$ and the formula for $\cZ_T(C)$ in terms of intrinsic volumes. After a change of variable, this formula is the integral representation of the Wills functional \cite{Hadwiger1975}. We repeat the statements in Section \ref{sec:prelimB} without proofs. In a few places we use an estimate which appears inside a proof in \cite{DarapaneniConvex}. We say so where this happens.

The link between \eqref{eq:FP} and the control problem through a logarithmic change of the unknown goes back to \cite{Fleming1978}. A complete account is given in \cite[Chapter VI]{FlemingSoner2006}. See also \cite{BoueDupuis1998,EvansIshii1985}.

\subsection{Organisation}
Section \ref{sec:prelimB} contains the notation and the results of \cite{DarapaneniConvex} which we use. Section \ref{sec:meta} is the main part of the paper. It describes the geometry of the basins and the Gibbs state of a union. It then proves the freezing of the basin weights, the estimates for the exit times and the convergence to equilibrium. We treat the control problem in Section \ref{sec:hjbdisc} and collect some open questions in Section \ref{sec:open}.

\section{Preliminaries and results from \texorpdfstring{\cite{DarapaneniConvex}}{the companion paper}}\label{sec:prelimB}

In this section we fix the notation and we recall the results of the companion paper \cite{DarapaneniConvex} which we use below. We state them as they appear in that paper and we do not prove them again. After each statement we give its place in \cite{DarapaneniConvex}. Most of these results concern a single convex set.

\subsection{Metric projections and the state equation}\label{subsec:Bproj}
Let $H$ be a real separable Hilbert space with inner product $\ip{\cdot}{\cdot}$ and norm $|\cdot|$. For a nonempty closed set $C\subset H$ we put $d_C(x)=\inf\{|x-c|:c\in C\}$ and $V_C=\frac12d_C^2$. If $C$ is convex, then every $x\in H$ has a unique nearest point $P_Cx$ in $C$. We now recall some standard facts. They hold for every nonempty closed convex set $C$ and all $x,y\in H$. All of them are in \cite[Chapters 3, 4 and 12]{BauschkeCombettes2017}.
\begin{enumerate}[label=\textup{(P\arabic*)},leftmargin=3em]
\item\label{P1} $\ip{x-P_Cx}{c-P_Cx}\le0$ for every $c\in C$. The point $P_Cx$ is the only point of $C$ with this property.
\item\label{P2} $\ip{(x-P_Cx)-(y-P_Cy)}{x-y}\ge|(x-P_Cx)-(y-P_Cy)|^2$. Hence the map $x\mapsto x-P_Cx$ is monotone and Lipschitz with constant $1$. A map which is Lipschitz with constant $1$ is called nonexpansive. The map $P_C$ is also nonexpansive.
\item\label{P3} The function $d_C$ is convex and nonexpansive.
\item\label{P4} The function $V_C$ is Fr\'echet differentiable and $\nabla V_C(x)=x-P_Cx$.
\item\label{P5} For $\lambda\in[0,1]$ we have $P_C\big(P_Cx+\lambda(x-P_Cx)\big)=P_Cx$.
\end{enumerate}
These facts carry the same names in \cite[Section~2]{DarapaneniConvex}.
The Hausdorff distance of two nonempty closed sets is
\[
\haus(C,C')=\max\Big\{\sup_{c\in C}d_{C'}(c),\ \sup_{c'\in C'}d_C(c')\Big\}\in[0,\infty].
\]

\begin{definition}\label{def:path}
A \emph{moving convex set of speed $v\ge0$} is a family $\{C(t)\}_{t\ge0}$ of nonempty closed convex subsets of $H$ such that
\[
\haus(C(s),C(t))\le v\,|t-s|\qquad(s,t\ge0).
\]
\end{definition}

Fix a moving convex set $\{C(t)\}$ of speed $v$, a point $x_0\in H$ and a continuous map $w:[0,\infty)\to H$ with $w(0)=0$. We consider the integral equation
\begin{equation}\label{eq:state}
x(t)=x_0-\int_0^t\big(x(u)-P_{C(u)}x(u)\big)\,du+w(t)\qquad(t\ge0),
\end{equation}
and together with it the linear equation
\begin{equation}\label{eq:ou}
\zeta(t)=-\int_0^t\zeta(u)\,du+w(t)\qquad(t\ge0).
\end{equation}

\begin{lemma}\label{lem:wellposed}
Equation \eqref{eq:state} has a unique continuous solution $x$. Equation \eqref{eq:ou} has the unique continuous solution
\begin{equation}\label{eq:ouformula}
\zeta(t)=w(t)-\int_0^te^{-(t-u)}w(u)\,du .
\end{equation}
\end{lemma}
This lemma is \cite[Lemma~2.3]{DarapaneniConvex}. We also take Definition \ref{def:path} and the equations \eqref{eq:state} and \eqref{eq:ou} from \cite[Section~2]{DarapaneniConvex}.

In the stochastic setting $W$ is a Wiener process in $H$ with covariance operator $Q$, where $Q$ is symmetric, nonnegative and of trace class. We refer to \cite{DPZ2014} for these notions. We write $\tr Q$ for the trace of $Q$ and $\|Q\|$ for its operator norm. For $T>0$ and an initial state $X_0$ independent of $W$, the \emph{state process} $X=X^T$ and the \emph{envelope process} $\zeta=\zeta^T$ are the solutions of \eqref{eq:state} and \eqref{eq:ou} with $w=\sqrt{2T}\,W$. In differential form,
\begin{equation}\label{eq:SDE}
dX_t=-\big(X_t-P_{C(t)}X_t\big)\,dt+\sqrt{2T}\,dW_t,\qquad d\zeta_t=-\zeta_t\,dt+\sqrt{2T}\,dW_t,\quad\zeta_0=0 .
\end{equation}
The variable $\zeta_t$ is a Wiener integral of a deterministic function. Hence it is a centred Gaussian variable with covariance operator
\begin{equation}\label{eq:oucov}
2T\int_0^te^{-2(t-u)}du\,Q=T(1-e^{-2t})\,Q .
\end{equation}
Taking the trace we also get $\E|\zeta_t|^2=T(1-e^{-2t})\tr Q$. When $H=\R^d$ we always take $Q$ to be the identity. Then $W$ is a standard Brownian motion and $\zeta$ is the OU process of the Introduction. Both displays \eqref{eq:SDE} and \eqref{eq:oucov} come from \cite[Section~2]{DarapaneniConvex}, the second one being \cite[equation~(2.7)]{DarapaneniConvex}.

\emph{Sets which are not convex.} Let $C\subset\R^d$ be nonempty and closed. The set of nearest points $\Pi_C(x)=\{c\in C:|x-c|=d_C(x)\}$ is nonempty and compact for every $x$. The function $V_C$ is locally Lipschitz and it is differentiable almost everywhere. At each point $x$ of differentiability the set $\Pi_C(x)$ has exactly one point $p$, and $\nabla V_C(x)=x-p$. We fix once and for all a Borel map $P_C:\R^d\to C$ with $P_Cx\in\Pi_C(x)$ for all $x$. Such a map exists by the measurable selection theorem \cite[Chapter 14]{RockafellarWets1998}. When $C$ is convex this is the metric projection. The state equation in $\R^d$ is the equation \eqref{eq:SDEintro} of the Introduction, which we repeat for the closed set $C$:
\begin{equation}\label{eq:SDEgen}
dX_t=-\big(X_t-P_CX_t\big)\,dt+\sqrt{2T}\,dW_t .
\end{equation}
We recall one fact from \cite[Section~2.3]{DarapaneniConvex}, with the name used there.
\begin{enumerate}[label=\textup{(S\arabic*)},leftmargin=3em]
\item\label{S1} For every initial condition independent of $W$, equation \eqref{eq:SDEgen} has a unique strong solution. The drift is bounded on every ball, so strong existence and pathwise uniqueness hold up to the exit time from each ball \cite{Veretennikov1981,KrylovRockner2005}. The linear growth of the drift together with the Gronwall lemma excludes explosion. To see this, fix $S>0$. As long as the solution is defined and $t\le S$, we have
\[
|X_t|\le|X_0|+\int_0^t\big(|X_u|+d_C(0)\big)\,du+\sqrt{2T}\,|W_t|\le a_S+\int_0^t|X_u|\,du ,
\]
where $a_S=|X_0|+d_C(0)\,S+\sqrt{2T}\,\sup_{u\le S}|W_u|$. The Gronwall lemma gives $|X_t|\le a_Se^{t}\le a_Se^{S}$, which is finite almost surely. It follows that the martingale problem for \eqref{eq:SDEgen} is well posed. The solutions started from the points of $\R^d$ then form a strong Markov family with measurable transition kernels \cite[Chapter 6]{StroockVaradhan1979}. Also, solutions driven by the same Brownian motion restart consistently at stopping times. This means that after a finite stopping time $\tau$ the solution is the solution which starts at $X_\tau$ and is driven by the Brownian motion $W_{\tau+t}-W_\tau$, $t\ge0$.
\end{enumerate}
We denote by $P_t(x,\cdot)$ the transition probabilities of \eqref{eq:SDEgen}. For an initial law $\rho_0$ we write $\rho_{T,t}$ for the law of $X_t$. The solution started at the point $x$ is denoted by $X^x$.

\subsection{The Gibbs state and distances between laws}\label{subsec:Bgibbs}
Let $C\subset\R^d$ be compact and $T>0$. The number
\[
\cZ_T(C)=\int_{\R^d}e^{-d_C(x)^2/2T}\,dx
\]
is finite \cite[Section~2]{DarapaneniConvex}. The \emph{Gibbs state} of $C$ at temperature $T$ is the probability measure
\[
\gamma_T^C(dx)=\cZ_T(C)^{-1}e^{-d_C(x)^2/2T}\,dx=\cZ_T(C)^{-1}e^{-V_C(x)/T}\,dx .
\]
\begin{lemma}\label{lem:invariant}
Let $C\subset\R^d$ be compact and $T>0$. Then $\gamma_T^C$ is an invariant probability measure for \eqref{eq:SDEgen}. For every initial law $\rho_0$ and every Borel set $A\subset\R^d$,
\[
\rho_{T,t}(A)\longrightarrow\gamma_T^C(A)\qquad(t\to\infty).
\]
\end{lemma}
This is \cite[Lemma~2.4]{DarapaneniConvex}.

For probability measures $\mu,\nu$ on $\R^d$ we use the total variation distance $\TV(\mu,\nu)=\sup_A|\mu(A)-\nu(A)|$, where $A$ runs over the Borel sets. For measures with densities $g$ and $g'$ it is equal to $\frac12\int|g-g'|\,dx$. When the second moments are finite, we also use the Wasserstein distance
\[
\Wt(\mu,\nu)=\inf\big\{(\E|\xi-\xi'|^2)^{1/2}:\ \Law(\xi)=\mu,\ \Law(\xi')=\nu\big\}.
\]
By \cite[Chapter 7]{Villani2003}, $\Wt(\mu_n,\mu)\to0$ if and only if $\mu_n\to\mu$ weakly and the second moments converge. The same reference shows that $\Wt^2$ is convex in each argument.
We recall the formula $\TV(\mu,\nu)=\frac12\sup\{|\int f\,d\mu-\int f\,d\nu|:f\text{ continuous},\ |f|\le1\}$. It shows that the total variation distance is lower semicontinuous for weak convergence. That is, if $\mu_n\to\mu$ and $\nu_n\to\nu$ weakly, then
\begin{equation}\label{eq:TVlsc}
\TV(\mu,\nu)\le\liminf_{n\to\infty}\TV(\mu_n,\nu_n).
\end{equation}
The image of a measure $\mu$ under a Borel map $f$ is denoted by $f_\#\mu$. The formula and the inequality \eqref{eq:TVlsc} come from \cite[Section~2]{DarapaneniConvex}, where the inequality is \cite[equation~(2.15)]{DarapaneniConvex}.

We shall use two elementary facts about total variation. For a nonnegative measure $\mu$ and a number $a\ge0$ we write $\mu\ge a\,m$ if $\mu(A)\ge a\,m(A)$ for all Borel sets $A$.

\begin{lemma}\label{lem:doeblin}
\begin{enumerate}[label=\textup{(\roman*)}]
\item Let $\mu_1,\mu_2,m_1,m_2$ be probability measures on $\R^d$ and $a\in[0,1]$ with $\mu_1\ge a\,m_1$ and $\mu_2\ge a\,m_2$. Then $\TV(\mu_1,\mu_2)\le1-a\,(1-\TV(m_1,m_2))$.
\item Let $\Pi$ be a Markov kernel on $\R^d$, let $M\subset\R^d$ be a Borel set and let $\alpha\in[0,1]$ be such that $\TV(\Pi(y,\cdot),\Pi(y',\cdot))\le1-\alpha$ for all $y,y'\in M$. Then for all probability measures $\mu,\mu'$ on $\R^d$
\[
\TV(\mu \Pi,\mu'\Pi)\le(1-\alpha)\,\TV(\mu,\mu')+\alpha\,\big(\mu(M^c)+\mu'(M^c)\big).
\]
\end{enumerate}
\end{lemma}
Both facts are \cite[Lemma~2.5]{DarapaneniConvex}.

\subsection{The envelope and the maximal inequality}\label{subsec:Benv}
\begin{theorem}\label{thm:envelope}
Let $\{C(t)\}_{t\ge0}$ be a moving convex set of speed $v$ in $H$, let $x_0\in H$ and let $w:[0,\infty)\to H$ be continuous with $w(0)=0$. Let $x$ and $\zeta$ be the solutions of \eqref{eq:state} and \eqref{eq:ou}. Then for every $t\ge0$
\begin{align}
d_{C(t)}\big(x(t)-\zeta(t)\big)&\le e^{-t}\,d_{C(0)}(x_0)+v\,(1-e^{-t}),\label{eq:env1}\\
d_{C(t)}\big(x(t)\big)&\le e^{-t}\,d_{C(0)}(x_0)+v\,(1-e^{-t})+|\zeta(t)|.\label{eq:env2}
\end{align}
In particular, the state process and the envelope process of \eqref{eq:SDE} satisfy, almost surely and for all $t\ge0$,
\begin{equation}\label{eq:env3}
d_{C(t)}(X_t)\le e^{-t}\,d_{C(0)}(X_0)+v\,(1-e^{-t})+|\zeta_t| .
\end{equation}
\end{theorem}
This is \cite[Theorem~3.1]{DarapaneniConvex}. For an arbitrary closed set a weaker statement holds.

\begin{theorem}\label{thm:compare}
Let $C\subset\R^d$ be nonempty and closed. Let $T>0$, $x\in\R^d$ and let $X$ be the solution of \eqref{eq:SDEgen} with $X_0=x$.
\begin{enumerate}[label=\textup{(\roman*)}]
\item There is a continuous adapted increasing process $A$ with $A_0=0$ such that, for all $t\ge0$,
\begin{equation}\label{eq:Vito}
V_C(X_t)=V_C(x)+\int_0^t\big(Td-2V_C(X_u)\big)\,du+\sqrt{2T}\int_0^t\ip{X_u-P_CX_u}{dW_u}-A_t .
\end{equation}
\item Let $v\in\R^d$ with $|v|=d_C(x)$. On an extension of the probability space there is a continuous process $\varrho$ which has the law of $(|e^{-t}v+\zeta_t|)_{t\ge0}$ and satisfies
\[
d_C(X_t)\le\varrho_t\qquad\text{for all }t\ge0,\ \text{almost surely}.
\]
\item For all $t\ge0$, $N>0$ and $r>0$,
\begin{equation}\label{eq:msq}
\begin{aligned}
\E\,d_C(X_t)^2&\le e^{-2t}d_C(x)^2+dT\,(1-e^{-2t}),\\
\P\big(d_C(X_t)\ge r\big)&\le\P\big(|\zeta_t|\ge r-d_C(x)\big),\\
\P\Big(\sup_{0\le t\le N}d_C(X_t)\ge r\Big)&\le\P\Big(\sup_{0\le t\le N}|\zeta_t|\ge r-d_C(x)\Big).
\end{aligned}
\end{equation}
\end{enumerate}
\end{theorem}
This comparison is \cite[Theorem~3.7]{DarapaneniConvex}.

Below, $a_+=\max\{a,0\}$ and $\log_+a=(\log a)_+$.

\begin{proposition}\label{prop:oumax}
Let $W$ be a Wiener process in $H$ with covariance $Q\neq0$, let $T>0$ and let $\zeta$ be the envelope process of \eqref{eq:SDE}. Then for all $N>0$, $r>0$ and $\eta\in(0,1)$,
\begin{equation}\label{eq:oumax}
\begin{aligned}
\P\Big(\sup_{0\le t\le N}|\zeta_t|\ge r\Big)\le\Big(\frac{16N}{\eta^2}+1\Big)\bigg[&\exp\Big(-\frac{\big((1-\eta)r-\sqrt{T\tr Q}\big)_+^2}{2T\|Q\|}\Big)\\
&+\exp\Big(\frac{\tr Q}{2\|Q\|}-\frac{r^2}{2T\|Q\|}\Big)\bigg].
\end{aligned}
\end{equation}
\end{proposition}
This is \cite[Proposition~3.9]{DarapaneniConvex}. We also need a bound from Step 2 of the proof of \cite[Proposition~3.9]{DarapaneniConvex}, where it is \cite[equation~(3.20)]{DarapaneniConvex}. For every $a\ge0$,
\begin{equation}\label{eq:gaussconc}
\P\big(|\zeta_a|\ge\sqrt{T\tr Q}+\ell\big)\le\exp\Big(-\frac{\ell^2}{2T\|Q\|}\Big)\qquad(\ell\ge0).
\end{equation}

\subsection{Convex sets: the Gibbs state}\label{subsec:Bconvex}
In this subsection and in the next one, $C\subset\R^d$ is a nonempty compact convex set. Let $k=\dim C\in\{0,1,\dots,d\}$ be the dimension of its affine hull and let $\cH^k$ be the Hausdorff measure of dimension $k$. We write
\[
u_C=\frac{\cH^k(C\cap\,\cdot\,)}{\cH^k(C)}
\]
for the uniform probability measure on $C$, which is a Dirac mass when $k=0$. The Lebesgue measure of $\R^d$ is denoted by $\Leb$.

The notions from convex geometry which we need are in \cite[Section 4.2]{Schneider2014}. Let $\omega_n=\pi^{n/2}/\Gamma(\frac n2+1)$ be the volume of the unit ball of $\R^n$. For every Borel set $E\subset\R^d$ and every $r>0$ the \emph{local Steiner formula} holds:
\begin{equation}\label{eq:steiner}
\Leb\big\{x\in\R^d:\ 0<d_C(x)\le r,\ P_Cx\in E\big\}=\sum_{j=0}^{d-1}r^{d-j}\,\omega_{d-j}\,\Phi_j(C,E).
\end{equation}
The left side is the volume of the set of points outside $C$ which are at distance at most $r$ from $C$ and whose nearest point lies in $E$. Here $\Phi_0(C,\cdot),\dots,\Phi_{d-1}(C,\cdot)$ are finite Borel measures concentrated on the boundary of $C$, called the \emph{curvature measures} of $C$. We put $\Phi_d(C,E)=\Leb(C\cap E)$. The total masses $V_j(C)=\Phi_j(C,\R^d)$ are the \emph{intrinsic volumes} of $C$. We shall use that $V_0(C)=1$, that $V_j(C)>0$ for $0\le j\le k$, that $V_j(C)=0$ and $\Phi_j(C,\cdot)=0$ for $j>k$, and that
\begin{equation}\label{eq:Phik}
\Phi_k(C,E)=\cH^k(C\cap E),\qquad V_k(C)=\cH^k(C)>0 .
\end{equation}
For a set with nonempty interior, $2V_{d-1}(C)$ is the surface area of $C$. We denote by $\chi_n$ a random variable distributed as the Euclidean norm of a standard Gaussian vector in $\R^n$, with $\chi_0=0$. For $n\ge1$ its density is
\begin{equation}\label{eq:chi}
\frac{2}{2^{n/2}\,\Gamma(n/2)}\,r^{n-1}e^{-r^2/2}\qquad(r>0).
\end{equation}
These facts also appear in \cite[Section~4]{DarapaneniConvex}.

\begin{theorem}\label{thm:gibbs}
Let $C\subset\R^d$ be a nonempty compact convex set of dimension $k$ and let $T>0$. Then
\begin{equation}\label{eq:wills}
\cZ_T(C)=\sum_{j=0}^{k}(2\pi T)^{(d-j)/2}\,V_j(C).
\end{equation}
Moreover, if $x$ is distributed according to $\gamma_T^C$, then the law of the pair $(P_Cx,d_C(x))$ is
\begin{equation}\label{eq:mixture}
\sum_{j=0}^{k}w_j(T)\ \frac{\Phi_j(C,\cdot)}{V_j(C)}\otimes\Law\big(\sqrt T\,\chi_{d-j}\big),\qquad w_j(T)=\frac{(2\pi T)^{(d-j)/2}\,V_j(C)}{\cZ_T(C)}.
\end{equation}
\end{theorem}
The theorem is \cite[Theorem~4.2]{DarapaneniConvex}.

\begin{corollary}\label{cor:zerotemp}
In the setting of Theorem \ref{thm:gibbs}, let $D=\diam C$. Then
\begin{equation}\label{eq:zt1}
\TV\big((P_C)_\#\gamma_T^C,\ u_C\big)\le1-w_k(T)\le\sum_{j=0}^{k-1}(2\pi T)^{(k-j)/2}\,\frac{V_j(C)}{V_k(C)},
\end{equation}
\begin{equation}\label{eq:zt2}
\Wt\big(\gamma_T^C,\ u_C\big)\le\sqrt{dT}+D\,\big(1-w_k(T)\big)^{1/2}.
\end{equation}
In particular $\gamma_T^C\to u_C$ as $T\to0$ and $\Wt(\gamma_T^C,u_C)=O(T^{1/4})$.
\end{corollary}
See \cite[Corollary~4.5]{DarapaneniConvex}.

\subsection{Convex sets: the limit diagram and the mixing time}\label{subsec:Bdiagram}
Throughout this subsection, except in Lemmas \ref{lem:ballistic} and \ref{lem:density}, $C\subset\R^d$ is a nonempty compact convex set, $k=\dim C$, $L=\aff C$ is its affine hull, $L_0=L-L$ is the linear subspace parallel to $L$. We also need the numbers
\[
D=\diam C,\quad r_C=\sup\{r>0:\ \text{some ball of radius $r$ of $L$ is contained in }C\}\ \ (k\ge1).
\]
The initial law $\rho_0$ is a probability measure on $\R^d$ with finite second moment. We put
\[
m_0=\Big(\int d_C^2\,d\rho_0\Big)^{1/2}.
\]
The state process $X=X^T$ solves \eqref{eq:SDEgen} with $\Law(X_0)=\rho_0$, and $\rho_{T,t}=\Law(X^T_t)$. We define the noiseless flow by
\begin{equation}\label{eq:flow}
\varphi_t(x)=P_Cx+e^{-t}(x-P_Cx)\qquad(t\ge0,\ x\in\R^d).
\end{equation}
By \ref{P5} with $\lambda=e^{-t}$ we have $P_C\varphi_t(x)=P_Cx$, hence $\varphi_t(x)-P_C\varphi_t(x)=e^{-t}(x-P_Cx)$ and
\begin{equation}\label{eq:flowode}
\frac{d}{dt}\varphi_t(x)=-e^{-t}(x-P_Cx)=-\big(\varphi_t(x)-P_C\varphi_t(x)\big).
\end{equation}
Since $\varphi_0(x)=x$, integration of \eqref{eq:flowode} gives
\[
\varphi_t(x)=x-\int_0^t\big(\varphi_u(x)-P_C\varphi_u(x)\big)\,du ,
\]
so $\varphi_t(x)$ is the solution of \eqref{eq:state} with $w=0$, for the constant set $C$ and the initial point $x$.

\emph{Reflected Brownian motion.} Let $k\ge1$. For $x\in C$ let $N_C(x)=\{n\in L_0:\ \ip{n}{c-x}\le0\text{ for all }c\in C\}$ be the normal cone of $C$ at $x$ relative to $L$. Let $\beta$ be a standard Brownian motion in $L_0$ and let $\xi$ be a random point of $C$ independent of $\beta$. By \cite{Tanaka1979} there is a unique pair $(R,K)$ of continuous adapted processes such that $R_s\in C$ for all $s\ge0$, $K$ has locally bounded variation with total variation process $|K|$, and
\begin{equation}\label{eq:rbm}
\begin{aligned}
R_s&=\xi+\sqrt2\,\beta_s-K_s,\qquad K_s=\int_0^sn_u\,d|K|_u,\\
n_u&\in N_C(R_u),\ |n_u|=1\quad\text{for }d|K|\text{ almost every }u.
\end{aligned}
\end{equation}
Here we identify $L$ with $\R^k$. In $\R^k$ the set $C$ is bounded and contains a ball, which is the situation of \cite{Tanaka1979}. If $x$ is an interior point of $C$ relative to $L$, then $N_C(x)=\{0\}$. So the last condition in \eqref{eq:rbm} says that $K$ acts only at the times when $R$ is on the boundary of $C$ relative to $L$. There it keeps $R$ in $C$. The pair $(R,K)$ is the image of $(\xi,\beta)$ under a measurable map which is defined path by path. The process $R$ is the \emph{reflected Brownian motion} in $C$ with generator $\Delta$. For a probability measure $\nu$ on $C$ we denote by $\mu_s[\nu]$ the law of $R_s$ when $\Law(\xi)=\nu$, so that $\mu_s[\nu]=\int\mu_s[\delta_x]\,\nu(dx)$. If $k=0$, so that $C=\{p\}$, we put $\mu_s[\nu]=\delta_p$. Finally we put
\begin{equation}\label{eq:Lambda}
\Lambda_C(S)=\frac{\sqrt3}{2\,r_C}\,\big(D^4+8D^2S+4k^2S^2\big)^{1/2}\quad(k\ge1),\qquad\Lambda_C(S)=0\quad(k=0).
\end{equation}
The notation is that of \cite[Section~5]{DarapaneniConvex}.

\begin{theorem}\label{thm:diagram}
Let $C$, $\rho_0$ and $\rho_{T,t}$ be as above and let $\mu_s=\mu_s[(P_C)_\#\rho_0]$.
\begin{enumerate}[label=\textup{(\alph*)}]
\item \textup{(Short times.)} For all $T>0$ and $t\ge0$,
\[
\Wt\big(\rho_{T,t},(\varphi_t)_\#\rho_0\big)\le\sqrt{2dTt},\qquad\Wt\big(\rho_{T,t},(P_C)_\#\rho_0\big)\le\sqrt{2dTt}+e^{-t}m_0 .
\]
\item \textup{(Times of order $1/T$.)} Let $0<T\le1/e$ and $S>0$. For every $s$ with $T\log(1/T)\le s\le S$,
\[
\begin{aligned}
\Wt\big(\rho_{T,s/T},\mu_s\big)^2\le{}&4T^2m_0^2+12\,d\,T\log\frac1T\\
&+2\Big[Tm_0+11\Big(T\Big(k+\log\Big(1+\frac ST\Big)\Big)\Big)^{1/2}\Big]\Lambda_C(S).
\end{aligned}
\]
\item \textup{(Long times at fixed temperature.)} For every $T>0$, $\Wt(\rho_{T,t},\gamma_T^C)\to0$ as $t\to\infty$, and $\Wt(\gamma_T^C,u_C)\to0$ as $T\to0$.
\item \textup{(All joint limits.)} Let $T_n\to0$ and $t_n\to\infty$ with $T_nt_n\to s\in[0,\infty]$. Then
\[
\Wt\big(\rho_{T_n,t_n},\mu_s\big)\to0,\qquad\text{where }\mu_0=(P_C)_\#\rho_0\ \text{ and }\ \mu_\infty=u_C .
\]
\end{enumerate}
\end{theorem}
This theorem is \cite[Theorem~5.1]{DarapaneniConvex}.

\begin{lemma}\label{lem:ballistic}
Let $C\subset\R^d$ be a nonempty closed convex set, let $X$ solve \eqref{eq:SDEgen} and let $\varphi$ be the flow \eqref{eq:flow}. Then for all $t\ge0$
\[
\E|X_t-\varphi_t(X_0)|^2\le2dTt,\qquad\E\sup_{0\le u\le t}|X_u-\varphi_u(X_0)|^2\le8dTt .
\]
\end{lemma}
The lemma is \cite[Lemma~5.2]{DarapaneniConvex}.

There is a universal constant $c_*$ such that every logarithmically concave probability measure $\gamma$ on a Euclidean space satisfies the Poincar\'e inequality
\begin{equation}\label{eq:KLS}
\int\Big(f-\int f\,d\gamma\Big)^2d\gamma\le c_P(\gamma)\int|\nabla f|^2\,d\gamma,\qquad c_P(\gamma)\le c_*\int\Big|y-\int z\,\gamma(dz)\Big|^2\gamma(dy),
\end{equation}
for all smooth bounded $f$ with bounded gradient \cite[Theorem 1.2 and (1.8)]{Bobkov1999}. For the uniform law on a convex body the corresponding isoperimetric inequality is \cite[Theorem 5.1]{KLS1995}. Here $c_P(\gamma)$ is the smallest constant for which the first inequality holds.

\begin{lemma}\label{lem:L2}
Let $C\subset\R^d$ be a nonempty compact convex set with diameter $D$, let $T>0$ and let $\nu$ be a probability measure on $\R^d$ which has a bounded density with compact support. Let $h_0$ be the density of $\nu$ with respect to $\gamma_T^C$, and let $\nu P_t$ be the law at time $t$ of the solution of \eqref{eq:SDEgen} started from $\nu$. Then for all $t\ge0$
\begin{equation}\label{eq:L2}
\TV\big(\nu P_t,\gamma_T^C\big)\le\frac12\,\exp\Big(-\frac{Tt}{c_*M_T}\Big)\,\|h_0-1\|_{L^2(\gamma_T^C)},
\end{equation}
where, for any fixed point $c$ of $C$,
\[
M_T:=\int|y-c|^2\,\gamma_T^C(dy)\le2D^2+2dT .
\]
\end{lemma}
This lemma appears as \cite[Lemma~5.4]{DarapaneniConvex}.

\begin{lemma}\label{lem:ergodic}
Let $k\ge1$. For every probability measure $\nu$ on $C$ and every $s\ge0$,
\[
\Wt\big(\mu_s[\nu],u_C\big)\le2D\exp\Big(-\frac{s}{(k+4)\,c_*D^2}\Big).
\]
\end{lemma}
See \cite[Lemma~5.5]{DarapaneniConvex}. Step 1 of the proof of \cite[Lemma~5.5]{DarapaneniConvex} contains a volume bound which we also need. There it is \cite[equation~(5.25)]{DarapaneniConvex}. Let $k\ge1$, let $\Leb$ denote, in this display only, the Lebesgue measure of the affine hull $L$ of $C$, let $x\in C$, $r\in(0,D]$ and $A=\{y\in C:|y-x|\le r\}$. Then
\begin{equation}\label{eq:volA}
\Leb(A)\ge\Big(\frac rD\Big)^k\Leb(C) .
\end{equation}
On $L$ the measure $\cH^k$ is the Lebesgue measure. We write $\bar B(x,r)$ for the closed ball with centre $x$ and radius $r$. So \eqref{eq:volA} says that $u_C(\bar B(x,r))\ge(r/D)^k$.

\begin{corollary}\label{cor:commute}
Let $C$ and $\rho_0$ be as in Theorem \ref{thm:diagram}.
\begin{enumerate}[label=\textup{(\roman*)}]
\item The iterated limits exist in $\Wt$ and
\[
\lim_{T\to0}\lim_{t\to\infty}\rho_{T,t}=u_C,\qquad\lim_{t\to\infty}\lim_{T\to0}\rho_{T,t}=(P_C)_\#\rho_0 .
\]
\item The set of all limit points of $\rho_{T,t}$ as $T\to0$ and $t\to\infty$ is the curve $\{\mu_s:0\le s\le\infty\}$.
\item The following three statements are equivalent. First, the two iterated limits in \textup{(i)} coincide. Second, the joint limit of $\rho_{T,t}$ as $T\to0$, $t\to\infty$ exists. Third, $(P_C)_\#\rho_0=u_C$.
\item If $C$ has nonempty interior, then $(P_C)_\#\rho_0=u_C$ holds only for $\rho_0=u_C$.
\end{enumerate}
\end{corollary}
The corollary is \cite[Corollary~5.6]{DarapaneniConvex}.

Recall that $\Pi_C(x)$ is the set of nearest points of $x$ in $C$ and that $\omega_d$ is the volume of the unit ball.

\begin{lemma}\label{lem:density}
Let $C\subset\R^d$ be nonempty and closed, let $C'\subset C$ be nonempty, closed and convex, and let $\vartheta\in(0,\infty]$ be such that $\Pi_C(x)=\{P_{C'}x\}$ for all $x$ with $d_{C'}(x)<\vartheta$. Let $\lambda\ge0$ and let $T>0$ satisfy $(\lambda+1+8\sqrt d)\sqrt T<\vartheta$. Let $y\in\R^d$ with $d_{C'}(y)\le\lambda\sqrt T$, let $p=P_{C'}y$, and let $X^y$ be the solution of \eqref{eq:SDEgen} with $X^y_0=y$. Then for every Borel set $A$ contained in the ball $\bar B(p,\sqrt T)$,
\[
\P\big(X^y_1\in A\big)\ge c_{d,\lambda}\,T^{-d/2}\,\Leb(A),\qquad c_{d,\lambda}:=\frac12\,(4\pi)^{-d/2}\exp\Big(-\frac14-\frac34(\lambda+1)^2-\frac d3\Big).
\]
\end{lemma}

This is \cite[Lemma~5.9]{DarapaneniConvex}. For a closed convex set $C$ we take $C'=C$ and $\vartheta=\infty$. Then there is no condition on $T$.

Let
\[
\lambda_d=\tfrac12+\sqrt{2d},\qquad\alpha_d=\tfrac14\,c_{d,\lambda_d}\,\omega_d\in(0,1),
\]
which depend only on the dimension. Lemma \ref{lem:density}, applied to a closed convex set with $A=\bar B(p,\sqrt T)$, gives $1\ge c_{d,\lambda}\,T^{-d/2}\,\omega_dT^{d/2}=c_{d,\lambda}\,\omega_d$ for every $\lambda\ge0$. So $\alpha_d\le\frac14$.

\begin{theorem}\label{thm:tvmix}
Let $C\subset\R^d$ be a nonempty compact convex set with diameter $D$, let $\ell_0>0$ and $0<T\le1$, and put
\[
r_T=\log_+\frac{4\ell_0}{\sqrt T}+1+\frac{c_*(2D^2+2d)}{T}\,\log\Big(2\Big(\frac{\sqrt e\,\cZ_1(C)}{\omega_d}\Big)^{1/2}T^{-d/4}\Big).
\]
Then for every $x$ with $d_C(x)\le\ell_0$, every integer $n\ge1$ and every $t\ge n\,r_T$,
\[
\TV\big(P_t(x,\cdot),\gamma_T^C\big)\le(1-\alpha_d)^n+2\exp\Big(-\frac{(\ell_0-\sqrt{dT})_+^2}{2T}\Big).
\]
\end{theorem}
We quote this theorem from \cite[Theorem~5.10]{DarapaneniConvex}. Here $\cZ_1(C)$ is the number $\cZ_T(C)$ for $T=1$.

\subsection{The control problem}\label{subsec:Bcontrol}
In this subsection $C\subset\R^d$ is compact, $T>0$ and $g:\R^d\to\R$ is bounded and Lipschitz. Consider an agent who can add a drift $a_u$ to the state equation \eqref{eq:SDEgen}. The agent pays $\frac12|a_u|^2$ per unit time for this drift and pays $g$ at the final time $t$. Let $v_T(t,x)$ be the least expected cost when the state starts at $x$. The logarithmic transformation \cite{Fleming1978}, \cite[Chapter VI]{FlemingSoner2006} expresses $v_T$ by the state process without control. As in \cite[Section~6]{DarapaneniConvex}, we take this expression as the \emph{definition} of $v_T$:
\begin{equation}\label{eq:log}
v_T(t,x)=-2T\log\E\exp\Big(-\frac{g(X^x_t)}{2T}\Big),
\end{equation}
where $X^x$ is the solution of \eqref{eq:SDEgen} started at $x$. All the results on $v_T$ below are statements about the right side of \eqref{eq:log}. For a curve $\phi:[0,t]\to\R^d$ which is absolutely continuous with square integrable derivative, and for a closed convex set $C$, we put
\[
J_t(\phi)=\frac12\int_0^t\big|\phi'(u)+\phi(u)-P_C\phi(u)\big|^2\,du .
\]
This is the energy of the control $a=\phi'+\phi-P_C\phi$ which makes $\phi$ a solution of the equation $\phi'=-(\phi-P_C\phi)+a$.
For a bounded function $g$ we write $\osc g=\sup g-\inf g$, and for a compact set $C$
\[
g_C=\min_{y\in\R^d}\big\{g(y)+d_C(y)^2\big\}.
\]
The minimum exists because $g$ is bounded and continuous and $d_C(y)\to\infty$ as $|y|\to\infty$. See \cite[Section~6]{DarapaneniConvex}.

\begin{theorem}\label{thm:hjb}
Let $g:\R^d\to\R$ be bounded and Lipschitz and let $v_T$ be defined by \eqref{eq:log}.
\begin{enumerate}[label=\textup{(\alph*)}]
\item Let $C\subset\R^d$ be compact. Then for every $x\in\R^d$
\[
\lim_{T\to0}\ \lim_{t\to\infty}v_T(t,x)=g_C .
\]
\item Let $C\subset\R^d$ be compact and convex. Then for all $t>0$ and $x\in\R^d$ the limit
\begin{equation}\label{eq:v0}
v_0(t,x):=\lim_{T\to0}v_T(t,x)=\inf\big\{g(\phi(t))+J_t(\phi):\ \phi(0)=x\big\}
\end{equation}
exists, and
\begin{equation}\label{eq:v0bounds}
\begin{aligned}
v_0(t,x)&\ge g_C-2e^{-t}d_C(x)\big(e^{-t}d_C(x)+\sqrt{\osc g}\big)&&(t>0),\\
v_0(t,x)&\le g_C+\frac{3\,(D+e^{-t/3}\varrho)^2}{t}+\frac t3\,e^{-2t/3}\varrho^2&&(t\ge3),
\end{aligned}
\end{equation}
where $D=\diam C$ and $\varrho=d_C(x)+\sqrt{\osc g}$. In particular
\[
\lim_{t\to\infty}\ \lim_{T\to0}v_T(t,x)=g_C=\lim_{T\to0}\ \lim_{t\to\infty}v_T(t,x).
\]
\end{enumerate}
\end{theorem}
This is \cite[Theorem~6.3]{DarapaneniConvex}.

\section{Disconnected sets: metastability}\label{sec:meta}

Suppose now that the feasible set has several components:
\begin{equation}\label{eq:union}
C=C_1\cup\dots\cup C_m,\qquad m\ge2,
\end{equation}
where $C_1,\dots,C_m$ are nonempty compact convex subsets of $\R^d$ which are pairwise disjoint. The potential $V_C=\min_iV_{C_i}$ has $m$ flat wells. Its gradient jumps across the set of points which are equally close to two of the components. A state which starts near one component stays near it for a very long time before it moves to another one. This is the metastability described in the Introduction. We prove three things. The description of \cite[Section~5]{DarapaneniConvex} holds separately in each well for a time which is exponentially long in $1/T$. We find the exponential rate of the exit time from a well. Finally we identify the limits before the first exit time and after the time of Theorem \ref{thm:equil}, at which equilibrium is reached.

\subsection{Geometry}
We fix the following notation for the set \eqref{eq:union}:
\begin{align*}
\delta_i&=\min_{j\ne i}\dist(C_i,C_j)>0,\qquad\delta_*=\min_i\delta_i,\qquad k_i=\dim C_i,\qquad k^*=\max_ik_i,\\
B_i&=\{x\in\R^d:\ d_{C_i}(x)<d_{C_j}(x)\text{ for all }j\ne i\},\qquad\cR=\R^d\setminus(B_1\cup\dots\cup B_m),\\
U_i&=\{x\in\R^d:\ d_{C_i}(x)<\delta_i/2\},\qquad\widetilde U_i=\{x\in\R^d:\ d_{C_i}(x)\le\delta_*/4\}.
\end{align*}
We call the sets $B_i$ the \emph{basins} (the Voronoi cells of the components) and the set $\cR$ the \emph{ridge}. The functions $d_{C_j}$ are continuous, so the basins are open and the ridge is closed. Two different basins do not meet, since $d_{C_i}(x)<d_{C_j}(x)$ and $d_{C_j}(x)<d_{C_i}(x)$ cannot hold together. The sets $U_i$ and $\widetilde U_i$ are neighbourhoods of $C_i$, and $\widetilde U_i\subset U_i$ because $\delta_*/4\le\delta_i/4<\delta_i/2$. We write $P_i=P_{C_i}$ and $\varphi^i_t(x)=P_ix+e^{-t}(x-P_ix)$.

\begin{lemma}\label{lem:geometry}
\begin{enumerate}[label=\textup{(\roman*)}]
\item $U_i\subset B_i$. For $x\in B_i$ we have $\Pi_C(x)=\{P_ix\}$ and $d_C(x)=d_{C_i}(x)$.
\item If $x\in B_i$, then $\varphi^i_t(x)\in B_i$ for all $t\ge0$.
\item The ridge $\cR$ is contained in a finite union of hypersurfaces of class $C^1$. In particular it has Lebesgue measure zero.
\item The lowest value of $V_C$ on the boundary of $B_i$ is $\delta_i^2/8$.
\end{enumerate}
\end{lemma}

\begin{proof}
(i) Let $x\in U_i$, let $j\ne i$ and let $c\in C_j$. Since $P_ix\in C_i$ and $c\in C_j$, we have $|P_ix-c|\ge\dist(C_i,C_j)\ge\delta_i$. Also $|x-P_ix|=d_{C_i}(x)<\delta_i/2$. By the triangle inequality,
\[
|x-c|\ge|P_ix-c|-|x-P_ix|\ge\delta_i-d_{C_i}(x).
\]
Taking the minimum over $c\in C_j$ we get
\[
d_{C_j}(x)\ge\delta_i-d_{C_i}(x)>\delta_i-\frac{\delta_i}2=\frac{\delta_i}2>d_{C_i}(x).
\]
This holds for every $j\ne i$, hence $x\in B_i$.

Now let $x\in B_i$. Since $d_{C_i}(x)<d_{C_j}(x)$ for all $j\ne i$,
\[
d_C(x)=\min_{1\le j\le m}d_{C_j}(x)=d_{C_i}(x).
\]
Let $c\in\Pi_C(x)$, so that $c\in C$ and $|x-c|=d_C(x)=d_{C_i}(x)$. If $c$ were a point of $C_j$ with $j\ne i$, we would have $|x-c|\ge d_{C_j}(x)>d_{C_i}(x)$, which is false. So $c\in C_i$ and $|x-c|=d_{C_i}(x)$. The nearest point of $x$ in the convex set $C_i$ is unique, hence $c=P_ix$. Conversely, $P_ix\in C$ and $|x-P_ix|=d_C(x)$, so $P_ix\in\Pi_C(x)$. This proves $\Pi_C(x)=\{P_ix\}$.

(ii) Let $x\in B_i$, $t\ge0$ and $z=\varphi^i_t(x)=P_ix+e^{-t}(x-P_ix)$. Then
\[
x-z=(1-e^{-t})(x-P_ix),\qquad\text{so}\qquad|x-z|=(1-e^{-t})\,d_{C_i}(x).
\]
By \ref{P5} with $\lambda=e^{-t}$ we have $P_iz=P_ix$. Therefore
\[
d_{C_i}(z)=|z-P_ix|=e^{-t}|x-P_ix|=e^{-t}d_{C_i}(x).
\]
Let $j\ne i$. By \ref{P3} the function $d_{C_j}$ is nonexpansive. Also $d_{C_j}(x)>d_{C_i}(x)$, because $x\in B_i$. Hence
\[
d_{C_j}(z)\ge d_{C_j}(x)-|x-z|>d_{C_i}(x)-(1-e^{-t})\,d_{C_i}(x)=e^{-t}d_{C_i}(x)=d_{C_i}(z).
\]
Thus $z\in B_i$.

(iii) Let $x\in\cR$ and let $i$ be an index with $d_{C_i}(x)=\min_ld_{C_l}(x)$. Since $x\notin B_i$, there is $j\ne i$ with $d_{C_j}(x)\le d_{C_i}(x)$, and then $d_{C_j}(x)=d_{C_i}(x)$ by the choice of $i$. So $x$ belongs to one of the sets
\[
E_{ij}=\{x\in\R^d:\ d_{C_i}(x)=d_{C_j}(x)\},\qquad i\ne j .
\]
A point $x$ of $E_{ij}$ does not belong to $C_i\cup C_j$. Indeed, if $x\in C_i$, then $d_{C_j}(x)=d_{C_i}(x)=0$, so $x\in C_i\cap C_j$, which is empty. The same holds if $x\in C_j$.

We show that $E_{ij}$ is a hypersurface. By \ref{P4} the function $V_{C_i}$ is differentiable with gradient $x-P_ix$. This gradient is continuous by \ref{P2}. Outside $C_i$ we have $V_{C_i}>0$ and $d_{C_i}=(2V_{C_i})^{1/2}$. By the chain rule $d_{C_i}$ is of class $C^1$ outside $C_i$ and
\[
\nabla d_{C_i}(x)=\frac{\nabla V_{C_i}(x)}{(2V_{C_i}(x))^{1/2}}=\frac{x-P_ix}{d_{C_i}(x)}\qquad(x\notin C_i).
\]
Hence the function $g=d_{C_i}-d_{C_j}$ is of class $C^1$ on the open set $\R^d\setminus(C_i\cup C_j)$, with
\[
\nabla g(x)=\frac{x-P_ix}{d_{C_i}(x)}-\frac{x-P_jx}{d_{C_j}(x)} .
\]
Suppose that $\nabla g(x)=0$ at a point $x\in E_{ij}$ and let $r=d_{C_i}(x)=d_{C_j}(x)>0$. Then
\[
\frac{x-P_ix}{r}=\frac{x-P_jx}{r},\qquad\text{that is,}\qquad P_ix=P_jx .
\]
This point belongs to $C_i\cap C_j$, which is impossible. So $\nabla g\ne0$ at every point of $E_{ij}=\{g=0\}$. By the implicit function theorem $E_{ij}$ is a hypersurface of class $C^1$. This proves the first statement, because $\cR$ is contained in the union of the sets $E_{ij}$ and there are finitely many of them.

Near each of its points, $E_{ij}$ is the graph of a function of class $C^1$ of $d-1$ of the coordinates. By the Fubini theorem such a graph has Lebesgue measure zero. Countably many of these graphs cover $E_{ij}$. Hence each $E_{ij}$ has measure zero, and so has $\cR$.

(iv) Let $x$ be a point of the boundary of $B_i$. Since $B_i$ is open, $x\notin B_i$, so there is $j\ne i$ with $d_{C_j}(x)\le d_{C_i}(x)$. Also $x$ is the limit of a sequence of points $x_n\in B_i$. For every $l\ne i$ we have $d_{C_i}(x_n)<d_{C_l}(x_n)$. The distance functions are continuous, so in the limit $d_{C_i}(x)\le d_{C_l}(x)$. Therefore
\[
d_{C_i}(x)=d_{C_j}(x)=d_C(x).
\]
The points $P_ix\in C_i$ and $P_jx\in C_j$ satisfy
\[
\delta_i\le\dist(C_i,C_j)\le|P_ix-P_jx|\le|P_ix-x|+|x-P_jx|=d_{C_i}(x)+d_{C_j}(x)=2\,d_C(x).
\]
So $d_C(x)\ge\delta_i/2$ and
\[
V_C(x)=\frac12\,d_C(x)^2\ge\frac12\cdot\frac{\delta_i^2}4=\frac{\delta_i^2}8 .
\]

We now show that the value $\delta_i^2/8$ is attained. Let $j$ be an index with $\dist(C_i,C_j)=\delta_i$. Since $C_i$ and $C_j$ are compact, there are $a\in C_i$ and $b\in C_j$ with $|a-b|=\delta_i$. Let $z=(a+b)/2$ be the midpoint of the segment $[a,b]$. Then
\[
d_{C_i}(z)\le|z-a|=\frac{\delta_i}2,\qquad d_{C_j}(z)\le|z-b|=\frac{\delta_i}2 .
\]
Let $l\ne i$ and $c\in C_l$. As in the proof of (i),
\[
|z-c|\ge|P_iz-c|-|z-P_iz|\ge\dist(C_i,C_l)-d_{C_i}(z)\ge\delta_i-\frac{\delta_i}2=\frac{\delta_i}2 ,
\]
and so $d_{C_l}(z)\ge\delta_i/2$ for every $l\ne i$. In the same way, for $c\in C_i$,
\[
|z-c|\ge|P_jz-c|-|z-P_jz|\ge\dist(C_i,C_j)-d_{C_j}(z)\ge\delta_i-\frac{\delta_i}2=\frac{\delta_i}2 ,
\]
and so $d_{C_i}(z)\ge\delta_i/2$. Putting these inequalities together,
\[
d_{C_i}(z)=d_{C_j}(z)=\frac{\delta_i}2,\qquad d_C(z)=\min_ld_{C_l}(z)=\frac{\delta_i}2 .
\]
Hence $z\notin B_i$ and $V_C(z)=\delta_i^2/8$. It remains to see that $z$ is on the boundary of $B_i$. Let $0\le s<1$ and $y_s=a+s(z-a)$. Then $d_{C_i}(y_s)\le|y_s-a|=s\,\delta_i/2<\delta_i/2$, so $y_s\in U_i\subset B_i$ by (i). Since $y_s\to z$ as $s\to1$, the point $z$ is a limit of points of $B_i$. As $z\notin B_i$, it is on the boundary of $B_i$.
\end{proof}

For $i\ne j$ we put $\Delta_{ij}=\dist(C_i,C_j)$, so that $\delta_i=\min_{j\ne i}\Delta_{ij}$. A \emph{nearest pair} of $(i,j)$ is a pair of points $a\in C_i$, $b\in C_j$ with $|a-b|=\Delta_{ij}$. Such a pair exists because $C_i$ and $C_j$ are compact. For a nearest pair we have $a=P_ib$ and $b=P_ja$. Indeed,
\[
\Delta_{ij}\le d_{C_i}(b)\le|b-a|=\Delta_{ij},
\]
so $a$ is the nearest point of $b$ in $C_i$, and the same argument gives $b=P_ja$. Not every pair of components communicates directly, because a third component may lie between them. The following definition describes the pairs for which this does not happen.

\begin{definition}\label{def:direct}
The pair $\{i,j\}$ is \emph{direct} if it has a nearest pair $(a,b)$ such that
\begin{equation}\label{eq:direct}
d_{C_l}(w)>\min\{|w-a|,|w-b|\}\qquad\text{for all }w\in[a,b]\text{ and all }l\notin\{i,j\}.
\end{equation}
We denote by $\Delta_*$ the smallest number $\Delta$ such that the graph with the vertices $1,\dots,m$ and the edges $\{i,j\}$, $\Delta_{ij}\le\Delta$, is connected.
\end{definition}

\begin{lemma}\label{lem:direct}
\begin{enumerate}[label=\textup{(\roman*)}]
\item If $\Delta_{ij}=\delta_i$, then $\{i,j\}$ is direct, and \eqref{eq:direct} holds for every nearest pair of $(i,j)$.
\item If $\{i,j\}$ is not direct, there is $l\notin\{i,j\}$ with $\Delta_{il}<\Delta_{ij}$ and $\Delta_{lj}<\Delta_{ij}$.
\item The graph whose edges are the direct pairs $\{i,j\}$ with $\Delta_{ij}\le\Delta_*$ is connected. Moreover $\delta_i\le\Delta_*$ for all $i$, and $\Delta_*=\delta_1=\delta_2$ if $m=2$.
\item Let $(a,b)$ be a nearest pair of $(i,j)$ satisfying \eqref{eq:direct}, let $\Delta=\Delta_{ij}$, $e=(b-a)/\Delta$ and $w_\lambda=a+\lambda e$. Then $P_iw_\lambda=a$ and $P_jw_\lambda=b$ for $0\le\lambda\le\Delta$. Moreover $w_\lambda\in B_i$ for $0\le\lambda<\Delta/2$, $w_\lambda\in B_j$ for $\Delta/2<\lambda\le\Delta$, and the midpoint $z=w_{\Delta/2}$ satisfies
\[
d_{C_i}(z)=d_{C_j}(z)=\frac\Delta2<d_{C_l}(z)\qquad(l\notin\{i,j\}).
\]
\end{enumerate}
\end{lemma}

\begin{proof}
Let $(a,b)$ be a nearest pair of $(i,j)$, let $\Delta=\Delta_{ij}$, $e=(b-a)/\Delta$ and $w_\lambda=a+\lambda e$ for $0\le\lambda\le\Delta$. Then $e$ is a unit vector, $w_0=a$, $w_\Delta=b$ and the points $w_\lambda$ fill the segment $[a,b]$. Since $a=P_ib$ and $b=P_ja$, we can write
\[
w_\lambda=P_ib+\frac\lambda\Delta\,(b-P_ib),\qquad w_\lambda=P_ja+\frac{\Delta-\lambda}\Delta\,(a-P_ja).
\]
The numbers $\lambda/\Delta$ and $(\Delta-\lambda)/\Delta$ belong to $[0,1]$. So property \ref{P5} gives $P_iw_\lambda=P_ib=a$ and $P_jw_\lambda=P_ja=b$. Since $|w_\lambda-a|=\lambda$ and $|w_\lambda-b|=\Delta-\lambda$, it follows that
\begin{equation}\label{eq:wlambda}
d_{C_i}(w_\lambda)=\lambda,\qquad d_{C_j}(w_\lambda)=\Delta-\lambda\qquad(0\le\lambda\le\Delta).
\end{equation}
We shall use the following fact. Let $p\in\R^d$ and $0<\lambda\le\Delta$ be such that
\begin{equation}\label{eq:xs5aray}
|p-a|=|p-w_\lambda|+|w_\lambda-a|\qquad\text{and}\qquad|p-a|=\Delta .
\end{equation}
Then $p=b$. Indeed, the first relation is the case of equality in the triangle inequality, so $w_\lambda$ lies on the segment $[a,p]$. Since $w_\lambda\ne a$, the point $p$ lies on the half line which starts at $a$ and passes through $w_\lambda$, that is, $p=a+|p-a|\,e$. With $|p-a|=\Delta$ this gives $p=a+\Delta e=b$.

(i) Let $\Delta=\delta_i$, let $l\notin\{i,j\}$, $0\le\lambda\le\Delta$ and $p_l=P_lw_\lambda$. Since $p_l\in C_l$ and $a\in C_i$, we have $|p_l-a|\ge\Delta_{il}\ge\delta_i=\Delta$. By the triangle inequality and \eqref{eq:wlambda},
\[
d_{C_l}(w_\lambda)=|w_\lambda-p_l|\ge|p_l-a|-|a-w_\lambda|\ge\Delta-\lambda=|w_\lambda-b|\ge\min\{|w_\lambda-a|,|w_\lambda-b|\}.
\]
We have to show that the first term is strictly larger than the last. Suppose that all the inequalities in this chain are equalities. If $\lambda=0$, the last one reads $\Delta=\min\{0,\Delta\}=0$, which is false. If $\lambda>0$, the first two equalities say that $|p_l-a|=|p_l-w_\lambda|+|w_\lambda-a|$ and $|p_l-a|=\Delta$. This is \eqref{eq:xs5aray} for $p=p_l$, hence $p_l=b$. Then $b\in C_l\cap C_j$, which is impossible. So
\[
d_{C_l}(w_\lambda)>\min\{|w_\lambda-a|,|w_\lambda-b|\}\qquad(0\le\lambda\le\Delta,\ l\notin\{i,j\}).
\]
This is \eqref{eq:direct} for the pair $(a,b)$. Since $(a,b)$ was an arbitrary nearest pair and a nearest pair exists, the pair $\{i,j\}$ is direct.

(ii) Fix a nearest pair $(a,b)$. Since $\{i,j\}$ is not direct, \eqref{eq:direct} fails for $(a,b)$. So there are $\lambda\in[0,\Delta]$ and $l\notin\{i,j\}$ with
\[
d_{C_l}(w_\lambda)\le\min\{|w_\lambda-a|,|w_\lambda-b|\}=\min\{\lambda,\Delta-\lambda\}.
\]
Let $p_l=P_lw_\lambda$, so that $|w_\lambda-p_l|=d_{C_l}(w_\lambda)$. Since $a\in C_i$ and $p_l\in C_l$,
\[
\Delta_{il}\le|a-p_l|\le|a-w_\lambda|+|w_\lambda-p_l|\le\lambda+\min\{\lambda,\Delta-\lambda\}\le\lambda+(\Delta-\lambda)=\Delta .
\]
Suppose that $\Delta_{il}=\Delta$. Then all the inequalities in this chain are equalities. If $\lambda=0$, then $d_{C_l}(a)\le\min\{0,\Delta\}=0$, so $a\in C_l\cap C_i$, which is impossible. If $\lambda>0$, then $|a-p_l|=|a-w_\lambda|+|w_\lambda-p_l|$ and $|a-p_l|=\Delta$, so $p_l=b$ by \eqref{eq:xs5aray}. Then $b\in C_l\cap C_j$, which is again impossible. Hence $\Delta_{il}<\Delta=\Delta_{ij}$.

We prove the second inequality in the same way, with the roles of $a$ and $b$ exchanged. Since $b\in C_j$ and $p_l\in C_l$,
\[
\Delta_{lj}\le|b-p_l|\le|b-w_\lambda|+|w_\lambda-p_l|\le(\Delta-\lambda)+\min\{\lambda,\Delta-\lambda\}\le(\Delta-\lambda)+\lambda=\Delta .
\]
Suppose that $\Delta_{lj}=\Delta$, so that all these inequalities are equalities. If $\lambda=\Delta$, then $d_{C_l}(b)\le\min\{\Delta,0\}=0$, so $b\in C_l\cap C_j$, which is impossible. If $\lambda<\Delta$, then $w_\lambda\ne b$, $|b-p_l|=|b-w_\lambda|+|w_\lambda-p_l|$ and $|b-p_l|=\Delta$. As in the proof of \eqref{eq:xs5aray}, the point $p_l$ lies on the half line which starts at $b$ and passes through $w_\lambda$, at distance $\Delta$ from $b$. So $p_l=b-\Delta e=a$ and $a\in C_l\cap C_i$, which is impossible. Hence $\Delta_{lj}<\Delta_{ij}$.

(iii) We claim that any two indices $i\ne j$ are joined by a path of direct pairs whose distances are at most $\Delta_{ij}$. Let $r_1<r_2<\dots<r_q$ be the different values taken by the finitely many numbers $\Delta_{i'j'}$. We prove by induction on $n$ that the claim holds for all pairs with $\Delta_{ij}=r_n$.

Let $n=1$ and $\Delta_{ij}=r_1$. If $\{i,j\}$ were not direct, part (ii) would give an index $l$ with $\Delta_{il}<r_1$, which contradicts the choice of $r_1$. So $\{i,j\}$ is direct, and the path consists of this one pair.

Now let $n\ge2$, suppose that the claim holds for all pairs with distance $r_1,\dots,r_{n-1}$, and let $\Delta_{ij}=r_n$. If $\{i,j\}$ is direct, there is nothing to prove. If it is not direct, let $l$ be as in (ii). The numbers $\Delta_{il}$ and $\Delta_{lj}$ are smaller than $r_n$, so they are among $r_1,\dots,r_{n-1}$. By the induction hypothesis $i$ and $l$ are joined by a path of direct pairs with distances at most $\Delta_{il}<\Delta_{ij}$. In the same way $l$ and $j$ are joined by a path of direct pairs with distances at most $\Delta_{lj}<\Delta_{ij}$. The first path followed by the second joins $i$ and $j$. This proves the claim.

The graph with the edges $\{i,j\}$, $\Delta_{ij}\le\Delta_*$, is connected by the definition of $\Delta_*$. By the claim each of its edges $\{i,j\}$ can be replaced by a path of direct pairs with distances at most $\Delta_{ij}\le\Delta_*$. Hence the graph whose edges are the direct pairs with distance at most $\Delta_*$ is connected.

For the second statement, fix $i$. In a connected graph on $m\ge2$ vertices the vertex $i$ has a neighbour $j$. For this $j$ we have $\Delta_{ij}\le\Delta_*$, and so
\[
\delta_i=\min_{j'\ne i}\Delta_{ij'}\le\Delta_{ij}\le\Delta_* .
\]
Finally let $m=2$. There is only one pair, so $\delta_1=\delta_2=\Delta_{12}$. The graph on the two vertices is connected exactly when it contains the edge $\{1,2\}$, that is, when $\Delta\ge\Delta_{12}$. So $\Delta_*=\Delta_{12}=\delta_1=\delta_2$.

(iv) We proved the statement about the projections before \eqref{eq:wlambda}. Let $l\notin\{i,j\}$. By \eqref{eq:direct} with $w=w_\lambda$,
\[
d_{C_l}(w_\lambda)>\min\{\lambda,\Delta-\lambda\}\qquad(0\le\lambda\le\Delta).
\]
Let $0\le\lambda<\Delta/2$. Then $\min\{\lambda,\Delta-\lambda\}=\lambda$. By \eqref{eq:wlambda}
\[
d_{C_i}(w_\lambda)=\lambda<\Delta-\lambda=d_{C_j}(w_\lambda),\qquad d_{C_i}(w_\lambda)=\lambda<d_{C_l}(w_\lambda).
\]
So $w_\lambda\in B_i$. Let $\Delta/2<\lambda\le\Delta$. Then $\min\{\lambda,\Delta-\lambda\}=\Delta-\lambda$. In this case
\[
d_{C_j}(w_\lambda)=\Delta-\lambda<\lambda=d_{C_i}(w_\lambda),\qquad d_{C_j}(w_\lambda)=\Delta-\lambda<d_{C_l}(w_\lambda).
\]
So $w_\lambda\in B_j$. For $\lambda=\Delta/2$ we have $\min\{\lambda,\Delta-\lambda\}=\Delta/2$. Now \eqref{eq:wlambda} gives $d_{C_i}(z)=d_{C_j}(z)=\Delta/2<d_{C_l}(z)$. This completes the proof.
\end{proof}

\begin{remark}\label{rem:connected}
For finite sets of points, part (ii) says that the graph of direct pairs contains the relative neighbourhood graph. Part (iii) then follows from the known fact that this graph contains a minimal spanning tree. This fact is stated in \cite{Toussaint1980} for points of the plane. The proof is the same in every dimension. The number $\Delta_*$ is the longest edge of such a tree. The sets $\{d_{C_i}\le r\}$ are convex. Two of them meet if and only if $\Delta_{ij}\le2r$. Hence $\{d_C\le r\}$ is connected if and only if $r\ge\Delta_*/2$. In terms of the potential, $\Delta_*^2/8$ is the smallest level $h$ such that $\{V_C\le h\}$ is connected. In the same way $\delta_i^2/8$ is the smallest level $h$ such that the connected component of $\{V_C\le h\}$ which contains $C_i$ contains another component of $C$.
\end{remark}

\subsection{The Gibbs state of a union}
The zero temperature limit of the Gibbs state is again explicit. No energy distinguishes the components, since $V_C=0$ on all of them. The selection is made first by dimension and then by volume.

\begin{proposition}\label{prop:gibbsunion}
Let $C$ be as in \eqref{eq:union}. Then, as $T\to0$,
\[
\gamma_T^C\longrightarrow\sum_{i=1}^mw_i\,u_{C_i}\ \text{ weakly},\qquad w_i=\frac{\cH^{k^*}(C_i)\,\one_{\{k_i=k^*\}}}{\sum_{j:\,k_j=k^*}\cH^{k^*}(C_j)} .
\]
\end{proposition}

\begin{proof}
Let $f$ be a bounded continuous function on $\R^d$. The space $\R^d$ is the union of the pairwise disjoint sets $B_1,\dots,B_m$ and $\cR$. The ridge $\cR$ has measure zero by Lemma \ref{lem:geometry}(iii). Also $d_C=d_{C_i}$ on $B_i$ by Lemma \ref{lem:geometry}(i). For each $i$ we write the integral over $B_i$ as the integral over $\R^d$ minus the integral over $\R^d\setminus B_i$. We also recall that $e^{-d_{C_i}^2/2T}dx=\cZ_T(C_i)\,\gamma_T^{C_i}(dx)$. This gives
\begin{equation}\label{eq:unionsplit}
\begin{aligned}
\int_{\R^d}f\,e^{-d_C^2/2T}\,dx&=\sum_{i=1}^m\int_{B_i}f\,e^{-d_{C_i}^2/2T}\,dx\\
&=\sum_{i=1}^m\Big[\cZ_T(C_i)\int f\,d\gamma_T^{C_i}-\int_{\R^d\setminus B_i}f\,e^{-d_{C_i}^2/2T}\,dx\Big].
\end{aligned}
\end{equation}

\emph{The remainder.} Since $U_i\subset B_i$, we have $\R^d\setminus B_i\subset\R^d\setminus U_i$. On this set $d_{C_i}\ge\delta_i/2$. So, at every point of it,
\[
\frac{d_{C_i}^2}{2T}=\frac{d_{C_i}^2}{4T}+\frac{d_{C_i}^2}{4T}\ge\frac{\delta_i^2}{16T}+\frac{d_{C_i}^2}{4T},\qquad\text{so}\qquad e^{-d_{C_i}^2/2T}\le e^{-\delta_i^2/16T}\,e^{-d_{C_i}^2/4T}.
\]
Hence
\[
\Big|\int_{\R^d\setminus B_i}f\,e^{-d_{C_i}^2/2T}\,dx\Big|\le\sup|f|\,e^{-\delta_i^2/16T}\int_{\R^d}e^{-d_{C_i}^2/4T}\,dx=\sup|f|\,e^{-\delta_i^2/16T}\,\cZ_{2T}(C_i).
\]
For $T\le1$ we have $e^{-d_{C_i}^2/4T}\le e^{-d_{C_i}^2/4}$. Integrating this inequality we get $\cZ_{2T}(C_i)\le\cZ_2(C_i)$. Since $d-k^*\ge0$ and $T^{-n}e^{-\kappa/T}\to0$ as $T\to0$ for all $n\ge0$ and $\kappa>0$,
\begin{equation}\label{eq:xs5arem}
(2\pi T)^{-(d-k^*)/2}\,\Big|\int_{\R^d\setminus B_i}f\,e^{-d_{C_i}^2/2T}\,dx\Big|\le\sup|f|\,\cZ_2(C_i)\,(2\pi T)^{-(d-k^*)/2}\,e^{-\delta_i^2/16T}\longrightarrow0 .
\end{equation}

\emph{The main terms.} By Theorem \ref{thm:gibbs}, more precisely by \eqref{eq:wills} and \eqref{eq:Phik} for the set $C_i$ of dimension $k_i$,
\[
\cZ_T(C_i)=\sum_{j=0}^{k_i}(2\pi T)^{(d-j)/2}\,V_j(C_i)=(2\pi T)^{(d-k_i)/2}\Big[\cH^{k_i}(C_i)+\sum_{j<k_i}(2\pi T)^{(k_i-j)/2}\,V_j(C_i)\Big].
\]
For $T\le1$ and $j<k_i$ we have $(2\pi T)^{(k_i-j)/2}\le(2\pi)^{d/2}\sqrt T$. Since also $\cH^{k_i}(C_i)>0$, we get
\[
\cZ_T(C_i)=(2\pi T)^{(d-k_i)/2}\,\cH^{k_i}(C_i)\,\big(1+O(\sqrt T)\big).
\]
We divide by $(2\pi T)^{(d-k^*)/2}$. Since $(d-k_i)-(d-k^*)=k^*-k_i\ge0$,
\begin{equation}\label{eq:xs5amain}
\frac{\cZ_T(C_i)}{(2\pi T)^{(d-k^*)/2}}=(2\pi T)^{(k^*-k_i)/2}\,\cH^{k_i}(C_i)\,\big(1+O(\sqrt T)\big)\longrightarrow
\begin{cases}\cH^{k^*}(C_i),&k_i=k^*,\\[2pt]0,&k_i<k^* .\end{cases}
\end{equation}
By Corollary \ref{cor:zerotemp}, $\gamma_T^{C_i}\to u_{C_i}$ weakly, so $\int f\,d\gamma_T^{C_i}\to\int f\,du_{C_i}$. These integrals are bounded by $\sup|f|$.

\emph{Conclusion.} We divide \eqref{eq:unionsplit} by $(2\pi T)^{(d-k^*)/2}$ and let $T\to0$. By \eqref{eq:xs5amain} the terms with $k_i<k^*$ disappear. The remainders disappear by \eqref{eq:xs5arem}. Thus
\[
(2\pi T)^{-(d-k^*)/2}\int_{\R^d}f\,e^{-d_C^2/2T}\,dx\longrightarrow\sum_{i:\,k_i=k^*}\cH^{k^*}(C_i)\int f\,du_{C_i} .
\]
For $f=1$ the left side is $(2\pi T)^{-(d-k^*)/2}\cZ_T(C)$ and the limit is $\sum_{j:\,k_j=k^*}\cH^{k^*}(C_j)$, which is positive. Taking the quotient,
\[
\int f\,d\gamma_T^C=\frac{\int_{\R^d}f\,e^{-d_C^2/2T}\,dx}{\cZ_T(C)}\longrightarrow\frac{\sum_{i:\,k_i=k^*}\cH^{k^*}(C_i)\int f\,du_{C_i}}{\sum_{j:\,k_j=k^*}\cH^{k^*}(C_j)}=\sum_{i=1}^mw_i\int f\,du_{C_i} .
\]
This is the claim.
\end{proof}

\subsection{Freezing of the basin weights}
Let $\rho_0$ be a probability measure on $\R^d$ with finite second moment and with $\rho_0(\cR)=0$. By Lemma \ref{lem:geometry}(iii) this holds for instance when $\rho_0$ has a density. Put
\[
a_i=\rho_0(B_i),\qquad\rho_0^i=\rho_0(\,\cdot\,|B_i)\quad\text{when }a_i>0 .
\]
Since the basins are pairwise disjoint and cover $\R^d\setminus\cR$, we have $\sum_ia_i=1$. For a probability measure $\nu$ on $C_i$ let $\mu^i_s[\nu]$ be the law at time $s$ of the reflected Brownian motion in $C_i$ started from $\nu$, as in Section \ref{subsec:Bdiagram}, with $\mu^i_\infty[\nu]=u_{C_i}$. By Lemma \ref{lem:geometry}(i) the selection $P_C$ coincides with $P_i$ on $B_i$. So the noiseless flow $\varphi_t(x)=P_Cx+e^{-t}(x-P_Cx)$ is equal to $\varphi^i_t$ on $B_i$.

\begin{theorem}\label{thm:meta}
Let $C$ be as in \eqref{eq:union}, let $\rho_0$ be as above and let $\rho_{T,t}$ be the law at time $t$ of the solution of \eqref{eq:SDEgen} started from $\rho_0$.
\begin{enumerate}[label=\textup{(\alph*)}]
\item For every $t\ge0$, $\rho_{T,t}\to(\varphi_t)_\#\rho_0$ weakly as $T\to0$, and $(\varphi_t)_\#\rho_0\to(P_C)_\#\rho_0$ weakly as $t\to\infty$.
\item Let $0<c<\delta_*^2/8$. Let $T_n\to0$ and $t_n\to\infty$ with $t_n\le\exp(c/T_n)$ and $T_nt_n\to s\in[0,\infty]$. Then
\[
\rho_{T_n,t_n}\longrightarrow\sum_{i=1}^ma_i\,\mu^i_s\big[(P_i)_\#\rho_0^i\big]\ \text{ weakly}.
\]
In particular, if $T_nt_n\to\infty$ and $t_n\le\exp(c/T_n)$, the limit is $\sum_ia_iu_{C_i}$.
\item For every $T>0$, $\rho_{T,t}(A)\to\gamma_T^C(A)$ for every Borel set $A$ as $t\to\infty$, and $\gamma_T^C\to\sum_iw_iu_{C_i}$ weakly as $T\to0$, with $w_i$ as in Proposition \ref{prop:gibbsunion}.
\end{enumerate}
\end{theorem}

Most of the work lies in showing that a state started in a basin does not see the other components for a long time. Let $X^x$ denote the solution of \eqref{eq:SDEgen} started at $x$. We write $X^{i,x}$ for the solution of the equation $dX=-(X-P_iX)dt+\sqrt{2T}dW$ for the convex set $C_i$, with the same Brownian motion and the same starting point.

\begin{lemma}\label{lem:confine}
\begin{enumerate}[label=\textup{(\roman*)}]
\item Let $x\in U_i$ and $N>0$. Then
\[
\P\big(X^x_t\ne X^{i,x}_t\text{ for some }t\le N\big)\le\P\Big(\sup_{0\le t\le N}|\zeta_t|\ge\frac{\delta_i}2-d_{C_i}(x)\Big).
\]
\item Let $x\in B_i$ and $0<c<\delta_i^2/8$. Then
\[
\P\big(X^x_t=X^{i,x}_t\text{ for all }t\le e^{c/T}\big)\longrightarrow1\qquad(T\to0).
\]
\item Let $x\in B_i$ and $t\ge0$. Then $\P(X^x_u=X^{i,x}_u\text{ for all }u\le t)\to1$ as $T\to0$.
\end{enumerate}
\end{lemma}

\begin{proof}
Let $x\in B_i$. We first show that the two processes coincide as long as they are in $B_i$. Let
\[
\sigma=\inf\{t\ge0:\ X^{i,x}_t\notin B_i\},\qquad\tau=\inf\{t\ge0:\ X^x_t\notin B_i\}
\]
be the first times at which $X^{i,x}$ and $X^x$ leave $B_i$, with the value $\infty$ if the process never leaves. For $u<\tau$ we have $X^x_u\in B_i$. So $P_CX^x_u=P_iX^x_u$ by Lemma \ref{lem:geometry}(i). Hence, for $t\le\tau$,
\begin{align*}
X^x_t&=x-\int_0^t\big(X^x_u-P_iX^x_u\big)\,du+\sqrt{2T}\,W_t ,\\
X^{i,x}_t&=x-\int_0^t\big(X^{i,x}_u-P_iX^{i,x}_u\big)\,du+\sqrt{2T}\,W_t .
\end{align*}
We subtract the second line from the first. The starting points and the terms with $W_t$ cancel. So the difference $\Xi_t=X^x_t-X^{i,x}_t$ satisfies
\[
\Xi_t=-\int_0^t\big[(X^x_u-P_iX^x_u)-(X^{i,x}_u-P_iX^{i,x}_u)\big]\,du\qquad(t\le\tau).
\]
By \ref{P2} the map $y\mapsto y-P_iy$ is nonexpansive. Hence the norm of the integrand is at most $|\Xi_u|$, which gives
\begin{equation}\label{eq:xs5agron}
|\Xi_t|\le\int_0^t|\Xi_u|\,du\qquad(t\le\tau).
\end{equation}
We now use the argument of the Gronwall lemma. Put $g(t)=\int_0^t|\Xi_u|\,du$ for $t\le\tau$. Then $g\ge0$, $g(0)=0$ and $g'(t)=|\Xi_t|\le g(t)$ by \eqref{eq:xs5agron}. Hence
\[
\frac{d}{dt}\big(e^{-t}g(t)\big)=e^{-t}\big(g'(t)-g(t)\big)\le0 ,\qquad\text{so}\qquad e^{-t}g(t)\le g(0)=0 .
\]
Thus $g=0$. Now \eqref{eq:xs5agron} gives $\Xi_t=0$ for $t\le\tau$. It follows that $X^{i,x}_t=X^x_t\in B_i$ for $t<\tau$. Suppose that $\tau<\infty$. The paths are continuous and the complement of $B_i$ is closed, so $X^x_\tau\notin B_i$. Therefore $X^{i,x}_\tau=X^x_\tau\notin B_i$. This means that $\sigma=\tau$. We have proved that
\begin{equation}\label{eq:coincide}
X^x_t=X^{i,x}_t\qquad\text{for all }t\le\sigma .
\end{equation}

(i) Suppose now that $x\in U_i$, so that $\delta_i/2-d_{C_i}(x)>0$. The process $X^{i,x}$ is the state process of the convex set $C_i$, which does not move. By Theorem \ref{thm:envelope} with $v=0$, for all $t\ge0$,
\[
d_{C_i}(X^{i,x}_t)\le e^{-t}d_{C_i}(x)+|\zeta_t|\le d_{C_i}(x)+|\zeta_t| .
\]
Suppose that $\sup_{t\le N}|\zeta_t|<\delta_i/2-d_{C_i}(x)$. Then for all $t\le N$
\[
d_{C_i}(X^{i,x}_t)<d_{C_i}(x)+\frac{\delta_i}2-d_{C_i}(x)=\frac{\delta_i}2 ,
\]
that is, $X^{i,x}_t\in U_i\subset B_i$. Hence $\sigma>N$. By \eqref{eq:coincide}, $X^x_t=X^{i,x}_t$ for all $t\le N$. So the event that $X^x_t\ne X^{i,x}_t$ for some $t\le N$ is contained in the event $\{\sup_{t\le N}|\zeta_t|\ge\delta_i/2-d_{C_i}(x)\}$. This is (i).

We prove (iii) before (ii), because the proof of (ii) uses \eqref{eq:tube}.

(iii) Let $x\in B_i$ and $t\ge0$. The curve $\Gamma_t=\{\varphi^i_u(x):0\le u\le t\}$ is compact, being the continuous image of the interval $[0,t]$. By Lemma \ref{lem:geometry}(ii) it is contained in the open set $B_i$. The distance to the complement of $B_i$ is a continuous function which is positive on $B_i$, so it has a positive minimum on $\Gamma_t$. Let $\mathsf d>0$ be smaller than this minimum. Lemma \ref{lem:ballistic} for the convex set $C_i$ and the starting point $x$ gives
\[
\E\sup_{0\le u\le t}|X^{i,x}_u-\varphi^i_u(x)|^2\le8dTt .
\]
By the Chebyshev inequality,
\begin{equation}\label{eq:tube}
\P\Big(\sup_{0\le u\le t}|X^{i,x}_u-\varphi^i_u(x)|\ge\mathsf d\Big)\le\frac1{\mathsf d^2}\,\E\sup_{0\le u\le t}|X^{i,x}_u-\varphi^i_u(x)|^2\le\frac{8dTt}{\mathsf d^2}\longrightarrow0\qquad(T\to0).
\end{equation}
On the complementary event, for every $u\le t$ the point $X^{i,x}_u$ is at distance less than $\mathsf d$ from the point $\varphi^i_u(x)$ of $\Gamma_t$. By the choice of $\mathsf d$ it belongs to $B_i$. Thus $X^{i,x}$ stays in $B_i$ up to time $t$, so $\sigma>t$ and \eqref{eq:coincide} applies. Therefore
\[
\P\big(X^x_u=X^{i,x}_u\text{ for all }u\le t\big)\ge1-\frac{8dTt}{\mathsf d^2}\longrightarrow1\qquad(T\to0).
\]

(ii) Fix $x\in B_i$ and $0<c<\delta_i^2/8$. As $\eta\to0$ the number $(1-\eta)^2(1-4\eta)^2\delta_i^2/8$ tends to $\delta_i^2/8>c$. So we can choose $\eta\in(0,1/8)$ so small that
\begin{equation}\label{eq:etachoice}
(1-\eta)^2(1-4\eta)^2\,\frac{\delta_i^2}8>c .
\end{equation}
The proof has two phases. In the first phase, whose length does not depend on $T$, the state comes close to $C_i$. In the second phase it stays in $U_i$.

\emph{First phase.} Define $t_1=\max\{0,\log(d_{C_i}(x)/(\eta\delta_i))\}$, where the logarithm is $-\infty$ when $d_{C_i}(x)=0$. If $d_{C_i}(x)\le\eta\delta_i$, then $t_1=0$. If $d_{C_i}(x)>\eta\delta_i$, then $e^{-t_1}=\eta\delta_i/d_{C_i}(x)$. In both cases, as in the proof of Lemma \ref{lem:geometry}(ii),
\[
d_{C_i}(\varphi^i_{t_1}(x))=e^{-t_1}d_{C_i}(x)\le\eta\delta_i .
\]
Choose $\mathsf d\in(0,\eta\delta_i]$ smaller than the distance of $\Gamma_{t_1}$ from the complement of $B_i$. Let $G_1$ be the event that $\sup_{u\le t_1}|X^{i,x}_u-\varphi^i_u(x)|<\mathsf d$. By \eqref{eq:tube} with $t=t_1$,
\[
\P(G_1)\ge1-\frac{8dTt_1}{\mathsf d^2} .
\]
On $G_1$ the process $X^{i,x}$ stays in $B_i$ on $[0,t_1]$. As in the proof of (iii), this gives $X^x=X^{i,x}$ on $[0,t_1]$. Since $d_{C_i}$ is nonexpansive and $\mathsf d\le\eta\delta_i$, on $G_1$
\[
d_{C_i}(X^x_{t_1})\le d_{C_i}(\varphi^i_{t_1}(x))+|X^{i,x}_{t_1}-\varphi^i_{t_1}(x)|\le d_{C_i}(\varphi^i_{t_1}(x))+\mathsf d\le2\eta\delta_i<\frac{\delta_i}2 ,
\]
where the last inequality holds because $\eta<1/8$. So $X^x_{t_1}\in U_i$ on $G_1$.

\emph{Second phase.} We now apply the argument of (i) after the time $t_1$. Put $W'_t=W_{t_1+t}-W_{t_1}$ and let $\zeta'$ be the solution of $d\zeta'=-\zeta'\,dt+\sqrt{2T}\,dW'$, $\zeta'_0=0$. This is the envelope process restarted from zero at the time $t_1$. It is driven by the increments of $W$ after $t_1$. Since $W'$ is again a standard Brownian motion, $\zeta'$ has the same law as $\zeta$. The equation of $X^{i,x}$ gives, for $t\ge0$,
\[
X^{i,x}_{t_1+t}=X^{i,x}_{t_1}-\int_0^t\big(X^{i,x}_{t_1+u}-P_iX^{i,x}_{t_1+u}\big)\,du+\sqrt{2T}\,W'_t .
\]
So $t\mapsto X^{i,x}_{t_1+t}$ is the solution of \eqref{eq:state} for the set $C_i$, with the starting point $X^{i,x}_{t_1}$ and the path $w=\sqrt{2T}\,W'$. By Theorem \ref{thm:envelope},
\[
d_{C_i}(X^{i,x}_{t_1+t})\le e^{-t}d_{C_i}(X^{i,x}_{t_1})+|\zeta'_t|\le d_{C_i}(X^{i,x}_{t_1})+|\zeta'_t|\qquad(t\ge0).
\]
Let $N>0$. On the event $G_1\cap\{\sup_{t\le N}|\zeta'_t|<\delta_i/2-2\eta\delta_i\}$ we have $d_{C_i}(X^{i,x}_{t_1})\le2\eta\delta_i$ and therefore
\[
d_{C_i}(X^{i,x}_{t_1+t})<2\eta\delta_i+\frac{\delta_i}2-2\eta\delta_i=\frac{\delta_i}2\qquad(0\le t\le N).
\]
On this event $X^{i,x}$ stays in $B_i$ on $[0,t_1]$ and in $U_i\subset B_i$ on $[t_1,t_1+N]$. Hence $\sigma>t_1+N$. By \eqref{eq:coincide} the two processes then coincide on $[0,t_1+N]$. Since $\zeta'$ and $\zeta$ have the same law, we get for every $N>0$
\[
\P\big(G_1\cap\{X^x_t\ne X^{i,x}_t\text{ for some }t\in[t_1,t_1+N]\}\big)\le\P\Big(\sup_{0\le t\le N}|\zeta_t|\ge\frac{\delta_i}2-2\eta\delta_i\Big).
\]

\emph{The two phases together.} On $G_1$ the two processes coincide on $[0,t_1]$. Suppose that $G_1$ holds and that the processes differ at some time of $[0,N]$. Then this time belongs to $[t_1,t_1+N]$. So the event that they differ somewhere on $[0,N]$ is contained in the union of the complement of $G_1$ and of the event in the last display. Hence
\[
\P\big(X^x_t\ne X^{i,x}_t\text{ for some }t\le N\big)\le\frac{8dTt_1}{\mathsf d^2}+\P\Big(\sup_{0\le t\le N}|\zeta_t|\ge r\Big),\qquad r=(1-4\eta)\,\frac{\delta_i}2 .
\]
We bound the last probability by Proposition \ref{prop:oumax}, applied with this $r$, the same $\eta$, $\tr Q=d$ and $\|Q\|=1$. Since $(1-\eta)r=(1-\eta)(1-4\eta)\delta_i/2$ and $r^2/2T=(1-4\eta)^2\delta_i^2/8T$, the probability that $X^x$ and $X^{i,x}$ differ somewhere on $[0,N]$ is at most
\[
\frac{8dTt_1}{\mathsf d^2}+\Big(\frac{16N}{\eta^2}+1\Big)\bigg[\exp\Big(-\frac{\big((1-\eta)(1-4\eta)\frac{\delta_i}2-\sqrt{dT}\big)_+^2}{2T}\Big)+\exp\Big(\frac d2-\frac{(1-4\eta)^2\delta_i^2}{8T}\Big)\bigg].
\]

\emph{The limit.} Let $N=e^{c/T}$. We show that each of the three terms tends to zero as $T\to0$. Put
\[
r_1=(1-\eta)(1-4\eta)\,\frac{\delta_i}2,\qquad\kappa_1=\frac{r_1^2}2=(1-\eta)^2(1-4\eta)^2\,\frac{\delta_i^2}8,\qquad\kappa_2=(1-4\eta)^2\,\frac{\delta_i^2}8 .
\]
By \eqref{eq:etachoice} we have $\kappa_1>c$. Since $(1-\eta)^2\le1$, also $\kappa_2\ge\kappa_1>c$. The first term $8dTt_1/\mathsf d^2$ tends to zero because $t_1$ and $\mathsf d$ do not depend on $T$. For the other two terms we use $N\ge1$, which gives
\[
\frac{16N}{\eta^2}+1\le\Big(\frac{16}{\eta^2}+1\Big)e^{c/T} .
\]
Let $T<r_1^2/d$, so that $\sqrt{dT}<r_1$. Then
\[
\frac{(r_1-\sqrt{dT})^2}{2T}=\frac{r_1^2}{2T}-\frac{r_1\sqrt{dT}}{T}+\frac{dT}{2T}=\frac{\kappa_1}T-\frac{r_1\sqrt d}{\sqrt T}+\frac d2 ,
\]
so the second term is at most
\[
\Big(\frac{16}{\eta^2}+1\Big)\exp\Big(\frac cT-\frac{\kappa_1}T+\frac{r_1\sqrt d}{\sqrt T}-\frac d2\Big)=\Big(\frac{16}{\eta^2}+1\Big)\exp\Big(-\frac1{\sqrt T}\Big(\frac{\kappa_1-c}{\sqrt T}-r_1\sqrt d\Big)-\frac d2\Big).
\]
This tends to zero, because $\kappa_1-c>0$ and so $(\kappa_1-c)/\sqrt T-r_1\sqrt d\to\infty$. The third term is at most
\[
\Big(\frac{16}{\eta^2}+1\Big)\exp\Big(\frac cT+\frac d2-\frac{\kappa_2}T\Big)=\Big(\frac{16}{\eta^2}+1\Big)\exp\Big(\frac d2-\frac{\kappa_2-c}T\Big),
\]
which tends to zero because $\kappa_2-c>0$. Hence the probability that $X^x$ and $X^{i,x}$ differ somewhere on $[0,e^{c/T}]$ tends to zero. This proves (ii).
\end{proof}

\begin{proof}[Proof of Theorem \ref{thm:meta}]
Let $f$ be a bounded Lipschitz function on $\R^d$ with Lipschitz constant $\ell$. Such functions determine weak convergence. By \ref{S1} the solutions started from the points of $\R^d$ form a Markov family, so the law of the solution started from $\rho_0$ satisfies $\int f\,d\rho_{T,t}=\int_{\R^d}\E f(X^x_t)\,\rho_0(dx)$. Since $\rho_0(\cR)=0$ and $\R^d\setminus\cR$ is the union of the pairwise disjoint basins,
\begin{equation}\label{eq:decomp}
\int f\,d\rho_{T,t}=\sum_{i=1}^m\int_{B_i}\E f(X^x_t)\,\rho_0(dx).
\end{equation}

(a) Fix $x\in B_i$ and $t\ge0$. The difference $f(X^x_t)-f(X^{i,x}_t)$ is zero when $X^x_t=X^{i,x}_t$. Otherwise its absolute value is at most $2\sup|f|$. Hence
\[
|\E f(X^x_t)-f(\varphi^i_t(x))|\le2\sup|f|\,\P(X^x_t\ne X^{i,x}_t)+\E|f(X^{i,x}_t)-f(\varphi^i_t(x))| .
\]
The first term on the right tends to zero as $T\to0$ by Lemma \ref{lem:confine}(iii). For the second term we use the Lipschitz property of $f$, the Schwarz inequality and Lemma \ref{lem:ballistic}:
\[
\E|f(X^{i,x}_t)-f(\varphi^i_t(x))|\le\ell\,\E|X^{i,x}_t-\varphi^i_t(x)|\le\ell\,\big(\E|X^{i,x}_t-\varphi^i_t(x)|^2\big)^{1/2}\le\ell\,\sqrt{2dTt}\longrightarrow0 .
\]
So $\E f(X^x_t)\to f(\varphi^i_t(x))$ as $T\to0$, for every $x\in B_i$. These functions of $x$ are bounded by $\sup|f|$. By dominated convergence in \eqref{eq:decomp},
\[
\int f\,d\rho_{T,t}\longrightarrow\sum_{i=1}^m\int_{B_i}f(\varphi^i_t(x))\,\rho_0(dx)=\int_{\R^d}f(\varphi_t(x))\,\rho_0(dx)=\int f\,d(\varphi_t)_\#\rho_0 ,
\]
where we used that $\varphi_t=\varphi^i_t$ on $B_i$ and $\rho_0(\cR)=0$. This is the first statement. For the second, let $x\in B_i$. Then
\[
\varphi_t(x)=P_ix+e^{-t}(x-P_ix)\longrightarrow P_ix=P_Cx\qquad(t\to\infty).
\]
So $f(\varphi_t(x))\to f(P_Cx)$ for every $x\notin\cR$. Dominated convergence gives
\[
\int f\,d(\varphi_t)_\#\rho_0=\int f(\varphi_t(x))\,\rho_0(dx)\longrightarrow\int f(P_Cx)\,\rho_0(dx)=\int f\,d(P_C)_\#\rho_0 .
\]

(b) Let $x\in B_i$. We take the processes at temperature $T_n$. Since $t_n\le e^{c/T_n}$ and $c<\delta_*^2/8\le\delta_i^2/8$, Lemma \ref{lem:confine}(ii) gives
\[
\P\big(X^x_{t_n}\ne X^{i,x}_{t_n}\big)\le1-\P\big(X^x_t=X^{i,x}_t\text{ for all }t\le e^{c/T_n}\big)\longrightarrow0 ,
\]
and as in (a)
\[
|\E f(X^x_{t_n})-\E f(X^{i,x}_{t_n})|\le2\sup|f|\,\P\big(X^x_{t_n}\ne X^{i,x}_{t_n}\big)\longrightarrow0 .
\]
The left side is bounded by $2\sup|f|$. By dominated convergence in \eqref{eq:decomp},
\[
\int f\,d\rho_{T_n,t_n}-\sum_{i=1}^m\int_{B_i}\E f(X^{i,x}_{t_n})\,\rho_0(dx)\longrightarrow0 .
\]
The term with index $i$ is zero when $a_i=0$. When $a_i>0$, the definition of $\rho_0^i$ shows that it is equal to $a_i\int\E f(X^{i,x}_{t_n})\,\rho_0^i(dx)$. So
\[
\int f\,d\rho_{T_n,t_n}-\sum_{i:\,a_i>0}a_i\int\E f(X^{i,x}_{t_n})\,\rho_0^i(dx)\longrightarrow0 .
\]
Let $a_i>0$. The integral $\int\E f(X^{i,x}_{t_n})\rho_0^i(dx)$ is the integral of $f$ with respect to the law at time $t_n$ of the state process for the convex set $C_i$ started from $\rho_0^i$. Since $\rho_0^i\le\rho_0/a_i$, the measure $\rho_0^i$ has finite second moment. By Theorem \ref{thm:diagram}(d) for the set $C_i$ this law converges to $\mu^i_s[(P_i)_\#\rho_0^i]$ in the distance $\Wt$, hence weakly. Therefore
\[
\int f\,d\rho_{T_n,t_n}\longrightarrow\sum_{i:\,a_i>0}a_i\int f\,d\mu^i_s\big[(P_i)_\#\rho_0^i\big] ,
\]
which is the claim. If $s=\infty$, then $\mu^i_\infty[(P_i)_\#\rho_0^i]=u_{C_i}$, so the limit is $\sum_ia_iu_{C_i}$.

(c) The first statement is Lemma \ref{lem:invariant} and the second is Proposition \ref{prop:gibbsunion}.
\end{proof}

\begin{corollary}\label{cor:three}
In the setting of Theorem \ref{thm:meta} there are three limits,
\[
\lim_{t\to\infty}\lim_{T\to0}\rho_{T,t}=(P_C)_\#\rho_0,\quad\lim_{\substack{T\to0,\ tT\to\infty\\ t\le e^{c/T}}}\rho_{T,t}=\sum_{i=1}^ma_i\,u_{C_i},\quad\lim_{T\to0}\lim_{t\to\infty}\rho_{T,t}=\sum_{i=1}^mw_i\,u_{C_i} .
\]
Moreover:
\begin{enumerate}[label=\textup{(\roman*)}]
\item the first and the second limit coincide if and only if $(P_i)_\#\rho_0^i=u_{C_i}$ for every $i$ with $a_i>0$,
\item the second and the third limit coincide if and only if $a_i=w_i$ for all $i$,
\item if $\rho_0$ gives positive mass to the basin of a component whose dimension is smaller than $k^*$, then the third limit is different from the first and from the second.
\end{enumerate}
\end{corollary}

\begin{proof}
Parts (a), (b) and (c) of Theorem \ref{thm:meta} give the three limits, in this order. We first write the first limit as a sum over the components. Let $A$ be a Borel set. Since $\rho_0(\cR)=0$ and $P_C=P_i$ on $B_i$,
\[
(P_C)_\#\rho_0(A)=\rho_0\big(P_C^{-1}(A)\big)=\sum_{i=1}^m\rho_0\big(B_i\cap P_i^{-1}(A)\big)=\sum_{i:\,a_i>0}a_i\,\rho_0^i\big(P_i^{-1}(A)\big).
\]
So the first limit is $(P_C)_\#\rho_0=\sum_ia_i(P_i)_\#\rho_0^i$, where the terms with $a_i=0$ are absent. Thus each of the three limits is a sum $\sum_i\nu_i$, where $\nu_i$ is a measure carried by $C_i$:
\[
\nu_i=a_i\,(P_i)_\#\rho_0^i,\qquad\nu_i=a_i\,u_{C_i},\qquad\nu_i=w_i\,u_{C_i}
\]
for the first, the second and the third limit respectively. The sets $C_1,\dots,C_m$ are pairwise disjoint, so $\nu_i$ is the restriction of the sum to $C_i$. Hence two such sums are equal if and only if their parts on each $C_i$ are equal.

(i) The first and the second limit coincide if and only if $a_i(P_i)_\#\rho_0^i=a_iu_{C_i}$ for all $i$. This holds for the indices with $a_i=0$. For $a_i>0$ it means that $(P_i)_\#\rho_0^i=u_{C_i}$.

(ii) By the same argument, the second and the third limit coincide if and only if $a_iu_{C_i}=w_iu_{C_i}$ for all $i$. Since $u_{C_i}$ is a probability measure, this means that $a_i=w_i$.

(iii) Suppose that $a_i>0$ and $k_i<k^*$. Then $w_i=0$ by the formula of Proposition \ref{prop:gibbsunion}. So the third limit gives no mass to $C_i$, while the first and the second limit give it the mass $a_i>0$.
\end{proof}

\subsection{Exit and transition times}
By Lemma \ref{lem:geometry}(iv) the number $\delta_i^2/8$ is the height of the potential barrier around the well $C_i$. Lemma \ref{lem:confine}(ii) says that, with probability tending to one, the state does not leave the basin $B_i$ before the time $\exp(c/T)$ when $c<\delta_i^2/8$. We now prove the matching upper bound. For $x\in B_i$ let
\[
\tau_i=\inf\{t\ge0:\ X^x_t\notin B_i\},\qquad\theta_i=\inf\Big\{t\ge0:\ X^x_t\in\bigcup_{l\ne i}\widetilde U_l\Big\}
\]
be the exit time from the basin and the first time at which the state is close to another component. Since $\widetilde U_l\subset B_l$ and $B_l$ does not meet $B_i$ for $l\ne i$, the state is outside $B_i$ at the time $\theta_i$. Hence $\tau_i\le\theta_i$.

\begin{theorem}\label{thm:exit}
Let $C$ be as in \eqref{eq:union} and let $x\in B_i$. Then for every $\eta>0$
\[
\P\Big(\exp\Big(\frac{\delta_i^2/8-\eta}{T}\Big)\le\tau_i\le\theta_i\le\exp\Big(\frac{\delta_i^2/8+\eta}{T}\Big)\Big)\longrightarrow1\qquad(T\to0).
\]
In other words, $T\log\tau_i$ and $T\log\theta_i$ converge to $\delta_i^2/8$ in probability.
\end{theorem}

For smooth potentials this is the Arrhenius law \cite{FreidlinWentzell2012}. When the wells and the saddles are nondegenerate, \cite{BEGK2004} give in addition the prefactor. To leading order, the logarithm of the exit time is the height of the barrier divided by $T$. We described the differences with our setting in the Introduction. The proof of the upper bound uses two estimates. One is a lower bound for the probability of following a given path (Lemma \ref{lem:tube}). The other is a bound for the probability that the state leaves the ridge on the side of the other component (Lemma \ref{lem:ridge}). Lemma \ref{lem:attempt} combines them into one attempt of transition. We then repeat this attempt.

\emph{The energy of a path.} Let $C'\subset\R^d$ be a nonempty closed convex set. For an absolutely continuous curve $\phi:[0,t]\to\R^d$ with square integrable derivative we put
\begin{equation}\label{eq:energy}
J^{C'}_t(\phi)=\frac12\int_0^t\big|\phi'(u)+\phi(u)-P_{C'}\phi(u)\big|^2\,du .
\end{equation}
This is the energy of the additional drift $a=\phi'+(\phi-P_{C'}\phi)$ which is needed to make $\phi$ a solution of the noiseless equation. Indeed, $\phi'=-(\phi-P_{C'}\phi)+a$.

\begin{lemma}\label{lem:tube}
Let $C'\subset\R^d$ be a nonempty closed convex set, $T>0$, $y\in\R^d$, and let $X^{C'}$ be the solution of $dX=-(X-P_{C'}X)\,dt+\sqrt{2T}\,dW$ with $X^{C'}_0=y$. Let $\phi:[0,t]\to\R^d$ be absolutely continuous with square integrable derivative. Then for every $r>0$
\[
\P\Big(\sup_{0\le u\le t}|X^{C'}_u-\phi(u)|\le\big(|y-\phi(0)|+r\big)e^{t}\Big)\ge\P\Big(\sup_{0\le u\le t}|W_u|<\frac{r}{\sqrt{2T}}\Big)\exp\Big(-\frac{J^{C'}_t(\phi)}{2T}\Big).
\]
\end{lemma}

\begin{proof}
Let $b(x)=-(x-P_{C'}x)$. By \ref{P2} the map $b$ is nonexpansive. We put
\[
h(u)=\frac1{\sqrt{2T}}\int_0^u\big(\phi'(s)-b(\phi(s))\big)\,ds\qquad(0\le u\le t).
\]
Then $h(0)=0$ and, for almost every $u$,
\[
h'(u)=\frac{\phi'(u)-b(\phi(u))}{\sqrt{2T}}=\frac{\phi'(u)+\phi(u)-P_{C'}\phi(u)}{\sqrt{2T}} .
\]
Comparing with \eqref{eq:energy} we get
\begin{equation}\label{eq:xs5bh}
\frac12\int_0^t|h'(u)|^2\,du=\frac1{2T}\cdot\frac12\int_0^t\big|\phi'(u)+\phi(u)-P_{C'}\phi(u)\big|^2\,du=\frac{J^{C'}_t(\phi)}{2T} .
\end{equation}
The curve $\phi$ satisfies an equation of the same form as $X^{C'}$, with $\sqrt{2T}\,h$ in the place of the noise:
\begin{equation}\label{eq:xs5bphi}
\begin{aligned}
\phi(u)&=\phi(0)+\int_0^u\phi'(s)\,ds=\phi(0)+\int_0^ub(\phi(s))\,ds+\int_0^u\big(\phi'(s)-b(\phi(s))\big)\,ds\\
&=\phi(0)+\int_0^ub(\phi(s))\,ds+\sqrt{2T}\,h(u).
\end{aligned}
\end{equation}

\emph{Change of measure.} The function $h'$ is deterministic and square integrable. So $\int_0^t\ip{h'}{dW}$ is a Gaussian variable with mean zero and variance $\int_0^t|h'|^2du$. Hence
\[
\E\exp\Big(\int_0^t\ip{h'}{dW}-\frac12\int_0^t|h'|^2\,du\Big)=1 .
\]
Let $\mathbf Q$ be the probability measure with density $\exp(\int_0^t\ip{h'}{dW}-\frac12\int_0^t|h'|^2du)$ with respect to $\P$. By the Girsanov theorem, the process $\widetilde W_u=W_u-h(u)$ is a Brownian motion on $[0,t]$ under $\mathbf Q$. We can write the equation of $X^{C'}$ with $\widetilde W$:
\[
X^{C'}_u=y+\int_0^ub(X^{C'}_s)\,ds+\sqrt{2T}\,W_u=y+\int_0^ub(X^{C'}_s)\,ds+\sqrt{2T}\,h(u)+\sqrt{2T}\,\widetilde W_u .
\]
We subtract \eqref{eq:xs5bphi}. The terms with $h$ cancel, so that
\[
X^{C'}_u-\phi(u)=y-\phi(0)+\int_0^u\big(b(X^{C'}_s)-b(\phi(s))\big)\,ds+\sqrt{2T}\,\widetilde W_u .
\]

\emph{The Gronwall lemma.} Put
\[
D(u)=|X^{C'}_u-\phi(u)|,\qquad A=|y-\phi(0)|+\sqrt{2T}\sup_{u\le t}|\widetilde W_u| .
\]
Since $|b(X^{C'}_s)-b(\phi(s))|\le D(s)$, the last identity gives
\[
D(u)\le A+\int_0^uD(s)\,ds\qquad(0\le u\le t).
\]
Let $F(u)=\int_0^uD(s)\,ds$. Then $F(0)=0$ and $F'=D\le A+F$, so that
\[
\frac{d}{du}\big(e^{-u}F(u)\big)=e^{-u}\big(F'(u)-F(u)\big)\le A\,e^{-u} .
\]
Integrating from $0$ to $u$ we get
\[
e^{-u}F(u)\le A\int_0^ue^{-s}\,ds=A\,(1-e^{-u}).
\]
Thus $F(u)\le A\,(e^u-1)$ and $D(u)\le A+F(u)\le A\,e^u$. For $u\le t$ this gives
\[
\sup_{u\le t}|X^{C'}_u-\phi(u)|\le\Big(|y-\phi(0)|+\sqrt{2T}\sup_{u\le t}|\widetilde W_u|\Big)e^{t} .
\]
Let $G=\{\sqrt{2T}\sup_{u\le t}|\widetilde W_u|<r\}$. On $G$ the right side is smaller than $(|y-\phi(0)|+r)e^t$. So the event in the statement of the lemma contains $G$. It is therefore enough to bound $\P(G)$ from below.

\emph{The density of $\P$ with respect to $\mathbf Q$.} Since $dW_u=d\widetilde W_u+h'(u)\,du$, we have
\[
\int_0^t\ip{h'}{dW}=\int_0^t\ip{h'}{d\widetilde W}+\int_0^t|h'|^2\,du .
\]
Hence
\[
\frac{d\P}{d\mathbf Q}=\exp\Big(-\int_0^t\ip{h'}{dW}+\frac12\int_0^t|h'|^2\,du\Big)=\exp\Big(-\int_0^t\ip{h'}{d\widetilde W}-\frac12\int_0^t|h'|^2\,du\Big).
\]

\emph{The Jensen inequality and symmetry.} Write $I=\int_0^t\ip{h'}{d\widetilde W}$. Since $\widetilde W$ is a Brownian motion under $\mathbf Q$, we have $\mathbf Q(G)=\P(\sup_{u\le t}|W_u|<r/\sqrt{2T})$. This number is positive. The function $s\mapsto e^{-s}$ is convex. By the Jensen inequality for the conditional law $\mathbf Q(\,\cdot\,|G)$,
\[
\E_{\mathbf Q}\big[\one_G\,e^{-I}\big]=\mathbf Q(G)\,\E_{\mathbf Q}\big[e^{-I}\,\big|\,G\big]\ge\mathbf Q(G)\,\exp\big(-\E_{\mathbf Q}[\,I\,|\,G\,]\big).
\]
Therefore
\[
\P(G)=\E_{\mathbf Q}\Big[\one_G\,\frac{d\P}{d\mathbf Q}\Big]\ge\mathbf Q(G)\,\exp\Big(-\frac12\int_0^t|h'|^2\,du\Big)\exp\Big(-\E_{\mathbf Q}\Big[\int_0^t\ip{h'}{d\widetilde W}\,\Big|\,G\Big]\Big).
\]
The conditional expectation vanishes. Indeed, under $\mathbf Q$ the process $-\widetilde W$ is again a Brownian motion. If we replace $\widetilde W$ by $-\widetilde W$, the event $G$ does not change, because it depends only on $\sup_{u\le t}|\widetilde W_u|$. The integral $I$ changes its sign. Hence
\[
\E_{\mathbf Q}[\,I\,\one_G\,]=\E_{\mathbf Q}[\,(-I)\,\one_G\,],\qquad\text{so that}\qquad\E_{\mathbf Q}[\,I\,\one_G\,]=0 .
\]
Finally we insert \eqref{eq:xs5bh} and the value of $\mathbf Q(G)$:
\[
\P(G)\ge\P\Big(\sup_{u\le t}|W_u|<\frac r{\sqrt{2T}}\Big)\exp\Big(-\frac{J^{C'}_t(\phi)}{2T}\Big).
\]
This completes the proof.
\end{proof}

In the next two lemmas $\{i,j\}$ is a direct pair and $(a,b)$ is a nearest pair satisfying \eqref{eq:direct}. We put
\[
\Delta=\Delta_{ij},\qquad e=\frac{b-a}{\Delta},\qquad w_\lambda=a+\lambda e,\qquad z=w_{\Delta/2},\qquad g=d_{C_i}-d_{C_j} .
\]

\begin{lemma}\label{lem:ridge}
There are $\mathsf r_0\in(0,\Delta/8]$ and $T_1>0$ with the following properties. Put $\mathsf r_1=\mathsf r_0/5$, $t_3=4\mathsf r_1/\Delta$ and
\[
K_{ij}=\{x\in\R^d:\ |x-z|\le\mathsf r_0,\ g(x)\ge\mathsf r_1\}.
\]
\begin{enumerate}[label=\textup{(\alph*)}]
\item In the ball $\bar B(z,\mathsf r_0)$ the basin $B_i$ is the set $\{g<0\}$, the basin $B_j$ is the set $\{g>0\}$ and the ridge $\cR$ is the set $\{g=0\}$. The set $K_{ij}$ is a compact subset of $B_j$, and $d_{C_j}\le5\Delta/8$ on it.
\item For all $T<T_1$ and all $y$ with $|y-z|\le\mathsf r_0/6$ and $g(y)=0$,
\[
\P\big(X^y_\varsigma\in K_{ij}\text{ for some }\varsigma\le t_3\big)\ge\frac2{35} .
\]
\end{enumerate}
\end{lemma}

\begin{proof}
\emph{Step 1: the choice of $\mathsf r_0$.} For $x\notin C_i$ let $n_i(x)=(x-P_ix)/d_{C_i}(x)$. This is a unit vector. The map $n_i$ is continuous on the complement of $C_i$, because $P_i$ and $d_{C_i}$ are continuous and $d_{C_i}>0$ there. By Lemma \ref{lem:direct}(iv) we have $P_iz=a$, $P_jz=b$ and
\[
d_{C_i}(z)=d_{C_j}(z)=\frac\Delta2<d_{C_l}(z)\qquad(l\notin\{i,j\}).
\]
Since $z-a=\frac\Delta2\,e$ and $z-b=-\frac\Delta2\,e$, it follows that
\[
n_i(z)=\frac{z-a}{\Delta/2}=e,\qquad n_j(z)=\frac{z-b}{\Delta/2}=-e,\qquad\ip{n_i(z)}{n_j(z)}=-1 .
\]
The functions $d_{C_i}$ and $d_{C_j}$ are nonexpansive by \ref{P3}. Hence, for every $x$ with $|x-z|\le\Delta/8$,
\[
\frac\Delta2-\frac\Delta8\le d_{C_i}(x)\le\frac\Delta2+\frac\Delta8,\qquad\frac\Delta2-\frac\Delta8\le d_{C_j}(x)\le\frac\Delta2+\frac\Delta8 .
\]
So the ball $\bar B(z,\Delta/8)$ does not meet $C_i\cup C_j$. Hence $n_i$ and $n_j$ are continuous on it. The functions $d_{C_l}-\max\{d_{C_i},d_{C_j}\}$, $l\notin\{i,j\}$, and the function $\ip{n_i}{n_j}$ are continuous on this ball. At the point $z$ the first ones are positive and the last one is equal to $-1$. Therefore there is $\mathsf r_0\in(0,\Delta/8]$ such that on the ball $\mathfrak B=\bar B(z,\mathsf r_0)$
\begin{equation}\label{eq:ridge1}
d_{C_l}>\max\{d_{C_i},d_{C_j}\}\quad(l\notin\{i,j\}),\qquad\ip{n_i}{n_j}\le-\frac12 .
\end{equation}
Since $\mathsf r_0\le\Delta/8$, the bounds obtained above hold on $\mathfrak B$, that is,
\begin{equation}\label{eq:ridge2}
\frac{3\Delta}8\le d_{C_i}\le\frac{5\Delta}8,\qquad\frac{3\Delta}8\le d_{C_j}\le\frac{5\Delta}8 .
\end{equation}
We prove (a). Let $x\in\mathfrak B$. By the first part of \eqref{eq:ridge1}, the inequalities $d_{C_i}(x)<d_{C_l}(x)$ and $d_{C_j}(x)<d_{C_l}(x)$ hold for all $l\notin\{i,j\}$. So $x\in B_i$ if and only if $d_{C_i}(x)<d_{C_j}(x)$, that is, $g(x)<0$. In the same way $x\in B_j$ if and only if $g(x)>0$. The point $x$ belongs to no other basin $B_l$, because $d_{C_l}(x)>d_{C_i}(x)$. Hence $x\in\cR$ if and only if $g(x)=0$. Since $g$ is continuous, the set $K_{ij}$ is closed. It is also bounded, hence compact. It is contained in $\mathfrak B\cap\{g>0\}\subset B_j$. The bound $d_{C_j}\le5\Delta/8$ on $K_{ij}$ is part of \eqref{eq:ridge2}.

\emph{Step 2: the process $g(X)$.} By \ref{P4}, on the open set $\{d_{C_i}>0\}$ the function $d_{C_i}=(2V_{C_i})^{1/2}$ is of class $C^1$ with gradient $\nabla V_{C_i}/d_{C_i}=n_i$. The map $n_i$ is locally Lipschitz there, because $x\mapsto x-P_ix$ is Lipschitz and $d_{C_i}$ is locally bounded from below. So the second derivative of $d_{C_i}$ exists almost everywhere. We bound it. The map $\nabla V_{C_i}(x)=x-P_ix$ is nonexpansive by \ref{P2}. Hence $\nabla^2V_{C_i}\le I$ almost everywhere, where $I$ is the identity matrix. Differentiating the identity $\nabla V_{C_i}=d_{C_i}\,n_i$ we get, almost everywhere on $\{d_{C_i}>0\}$,
\[
\nabla^2V_{C_i}=n_i\otimes n_i+d_{C_i}\,\nabla^2d_{C_i},\qquad\text{hence}\qquad d_{C_i}\,\nabla^2d_{C_i}=\nabla^2V_{C_i}-n_i\otimes n_i\le I-n_i\otimes n_i .
\]
On the other hand $\nabla^2d_{C_i}\ge0$, because $d_{C_i}$ is convex. So
\[
0\le\nabla^2d_{C_i}\le\frac{I-n_i\otimes n_i}{d_{C_i}}\qquad\text{almost everywhere on }\{d_{C_i}>0\}.
\]
The trace of $I-n_i\otimes n_i$ is $d-|n_i|^2=d-1$. So $0\le\tr\nabla^2d_{C_i}\le(d-1)/d_{C_i}$ almost everywhere. The same holds for $j$. The number $\tr\nabla^2g=\tr\nabla^2d_{C_i}-\tr\nabla^2d_{C_j}$ is the difference of two numbers which, by \eqref{eq:ridge2}, lie between $0$ and $(d-1)\cdot\frac8{3\Delta}$. Hence
\[
|\tr\nabla^2g|\le\kappa_g:=\frac{8(d-1)}{3\Delta}\qquad\text{almost everywhere on }\mathfrak B .
\]
Let $y$ be as in (b), $X=X^y$, and let $\varsigma_0$ be the exit time of $X$ from the interior of $\mathfrak B$. The It\^o formula holds for the function $g$, which has a locally Lipschitz gradient near $\mathfrak B$ \cite[Section 2.10]{Krylov1980}. It also follows by mollifying $g$, because $X_u$ has a density for $u>0$ by property (S2) in \cite[Section~2.3]{DarapaneniConvex}. The drift of $X$ is $b(x)=-(x-P_Cx)$ and its diffusion coefficient is $\sqrt{2T}$. Since also $g(y)=0$, for $t\le\varsigma_0$ we have
\[
\mathsf g_t:=g(X_t)=\int_0^t\mathsf b_u\,du+\mathsf N_t,
\]
with
\[
\mathsf b_u=\ip{\nabla g(X_u)}{b(X_u)}+T\,\tr\nabla^2g(X_u),\qquad \mathsf N_t=\sqrt{2T}\int_0^t\ip{\nabla g(X_u)}{dW_u} .
\]
We have $\nabla g=n_i-n_j$. By the second part of \eqref{eq:ridge1} and since $n_i$, $n_j$ are unit vectors, $-1\le\ip{n_i}{n_j}\le-\frac12$ on $\mathfrak B$. Hence
\[
|\nabla g|^2=|n_i|^2+|n_j|^2-2\ip{n_i}{n_j}=2-2\ip{n_i}{n_j}\in[3,4]\qquad\text{on }\mathfrak B .
\]
Since $d\langle \mathsf N\rangle_u=2T\,|\nabla g(X_u)|^2\,du$, this gives
\begin{equation}\label{eq:ridgeqv}
6T\,du\le d\langle \mathsf N\rangle_u\le8T\,du .
\end{equation}
We now bound the drift $\mathsf b_u$. On $\mathfrak B\cap\{g>0\}$ we are in $B_j$ by (a), so $P_C=P_j$ there by Lemma \ref{lem:geometry}(i) and $b=-d_{C_j}n_j$. We use $1-\ip{n_i}{n_j}\in[\frac32,2]$ and \eqref{eq:ridge2}:
\[
\ip{\nabla g}{b}=-d_{C_j}\ip{n_i-n_j}{n_j}=d_{C_j}\big(1-\ip{n_i}{n_j}\big)\in\Big[\frac{3\Delta}8\cdot\frac32,\ \frac{5\Delta}8\cdot2\Big]=\Big[\frac{9\Delta}{16},\ \frac{5\Delta}4\Big].
\]
On $\mathfrak B\cap\{g<0\}$ we are in $B_i$, so $b=-d_{C_i}n_i$ and
\[
\ip{\nabla g}{b}=-d_{C_i}\ip{n_i-n_j}{n_i}=-d_{C_i}\big(1-\ip{n_i}{n_j}\big)\in\Big[-\frac{5\Delta}4,\ -\frac{9\Delta}{16}\Big].
\]
Let $T$ be so small that $T\kappa_g\le\Delta/16$, that is, $128(d-1)T\le3\Delta^2$. Then $T|\tr\nabla^2g|\le\Delta/16$ on $\mathfrak B$. When $\mathsf g_u>0$ we get
\[
\mathsf b_u\ge\frac{9\Delta}{16}-\frac\Delta{16}=\frac\Delta2,\qquad\mathsf b_u\le\frac{5\Delta}4+\frac\Delta{16}=\frac{21\Delta}{16}\le2\Delta ,
\]
and when $\mathsf g_u<0$ we get
\[
\mathsf b_u\le-\frac{9\Delta}{16}+\frac\Delta{16}=-\frac\Delta2,\qquad\mathsf b_u\ge-\frac{5\Delta}4-\frac\Delta{16}=-\frac{21\Delta}{16}\ge-2\Delta .
\]
In short,
\begin{equation}\label{eq:ridgedrift}
\frac\Delta2\le\mathsf b_u\le2\Delta\ \text{ when }\mathsf g_u>0,\qquad-2\Delta\le\mathsf b_u\le-\frac\Delta2\ \text{ when }\mathsf g_u<0 .
\end{equation}
It remains to see that the times with $\mathsf g_u=0$ do not count. Let $\Lambda^s_t$ be the local time of the semimartingale $\mathsf g$ at the level $s$. By the occupation times formula, for $t\le\varsigma_0$,
\[
\int_0^{t}\one_{\{\mathsf g_u=0\}}\,d\langle \mathsf N\rangle_u=\int_{\R}\one_{\{s=0\}}\,\Lambda^s_t\,ds=0 .
\]
By \eqref{eq:ridgeqv} the left side is at least $6T$ times the Lebesgue measure of the set $\{u\le t:\ \mathsf g_u=0\}$. So the set of times $u\le\varsigma_0$ with $\mathsf g_u=0$ has Lebesgue measure zero almost surely.

\emph{Step 3: the state leaves the ridge quickly and stays in $\mathfrak B$.} We write $\mathrm{sgn}(s)=1$ for $s>0$ and $\mathrm{sgn}(s)=-1$ for $s\le0$. Let
\[
H_1=\Big\{\sup_{t\le t_3\wedge\varsigma_0}|\widetilde{\mathsf N}_t|<\mathsf r_1\Big\},\qquad H_2=\Big\{\sqrt{2T}\sup_{t\le t_3}|W_t|<\frac{\mathsf r_0}3\Big\},\qquad\widetilde{\mathsf N}_t=\int_0^t\mathrm{sgn}(\mathsf g_u)\,d\mathsf N_u .
\]
The process $\widetilde{\mathsf N}$ is a martingale with $d\langle\widetilde{\mathsf N}\rangle_u=d\langle \mathsf N\rangle_u\le8T\,du$ by \eqref{eq:ridgeqv}. By the Doob inequality,
\[
\P(H_1^c)\le\frac1{\mathsf r_1^2}\,\E\,\widetilde{\mathsf N}_{t_3\wedge\varsigma_0}^2=\frac1{\mathsf r_1^2}\,\E\,\langle\widetilde{\mathsf N}\rangle_{t_3\wedge\varsigma_0}\le\frac{8Tt_3}{\mathsf r_1^2} .
\]
In the same way, since $|W_t|^2$ is a submartingale and $\E|W_{t_3}|^2=d\,t_3$,
\[
\P(H_2^c)=\P\Big(\sup_{t\le t_3}|W_t|\ge\frac{\mathsf r_0}{3\sqrt{2T}}\Big)\le\frac{18T}{\mathsf r_0^2}\,\E|W_{t_3}|^2=\frac{18\,d\,T\,t_3}{\mathsf r_0^2} .
\]
Both probabilities tend to zero as $T\to0$.

We show that $\varsigma_0>t_3$ on $H_2$. On $\mathfrak B$ we have $|b|=d_C\le d_{C_i}\le5\Delta/8$ by \eqref{eq:ridge2}. Also $t_3=4\mathsf r_1/\Delta$ and $\mathsf r_1=\mathsf r_0/5$, so that $\frac{5\Delta}8\,t_3=\frac52\mathsf r_1=\frac{\mathsf r_0}2$. From $X_t=y+\int_0^tb(X_u)\,du+\sqrt{2T}\,W_t$ we get on $H_2$, for $t\le t_3\wedge\varsigma_0$,
\[
|X_t-z|\le|y-z|+\frac{5\Delta}8\,t_3+\sqrt{2T}\,|W_t|<\frac{\mathsf r_0}6+\frac52\mathsf r_1+\frac{\mathsf r_0}3=\frac{\mathsf r_0}6+\frac{\mathsf r_0}2+\frac{\mathsf r_0}3=\mathsf r_0 .
\]
If we had $\varsigma_0\le t_3$, this would give $|X_{\varsigma_0}-z|<\mathsf r_0$, which is false, because $X_{\varsigma_0}$ is on the boundary of $\mathfrak B$. So $\varsigma_0>t_3$ on $H_2$.

Next we use the Tanaka formula. Since $\mathsf g_0=0$, it reads, for $t\le\varsigma_0$,
\[
|\mathsf g_t|=\int_0^t\mathrm{sgn}(\mathsf g_u)\,d\mathsf g_u+\Lambda^0_t=\int_0^t\mathrm{sgn}(\mathsf g_u)\,\mathsf b_u\,du+\widetilde{\mathsf N}_t+\Lambda^0_t .
\]
The local time $\Lambda^0_t$ is nonnegative. By \eqref{eq:ridgedrift} we have $\mathrm{sgn}(\mathsf g_u)\,\mathsf b_u\ge\Delta/2$ whenever $\mathsf g_u\ne0$. By Step 2 this holds for almost every $u$. Hence, for $t\le\varsigma_0$,
\[
|\mathsf g_t|\ge\int_0^t\mathrm{sgn}(\mathsf g_u)\,\mathsf b_u\,du+\widetilde{\mathsf N}_t\ge\frac\Delta2\,t+\widetilde{\mathsf N}_t .
\]
On $H_1\cap H_2$ we may take $t=t_3$ and get
\[
|\mathsf g_{t_3}|>\frac\Delta2\,t_3-\mathsf r_1=2\mathsf r_1-\mathsf r_1=\mathsf r_1 .
\]
Let $\varsigma$ be the smallest of $\varsigma_0$, $t_3$ and the first time at which $|\mathsf g|=\mathsf r_1$. The process $|\mathsf g|$ is continuous and starts at $0$. On $H_1\cap H_2$ it is larger than $\mathsf r_1$ at the time $t_3<\varsigma_0$. So we have shown that on $H_1\cap H_2$ the time $\varsigma$ is the first time at which $|\mathsf g|=\mathsf r_1$.

\emph{Step 4: the side.} Let $\lambda_+=\Delta/(8T)$, $\lambda_-=2\Delta/(3T)$ and
\[
\kappa=\frac{\lambda_+}{\lambda_-}=\frac{\Delta}{8T}\cdot\frac{3T}{2\Delta}=\frac3{16} ,
\]
and define
\[
\Phi(s)=1-e^{-\lambda_+s}\quad(s\ge0),\qquad\Phi(s)=\kappa\,\big(e^{\lambda_-s}-1\big)\quad(s<0).
\]
The derivatives are
\begin{align*}
\Phi'(s)&=\lambda_+e^{-\lambda_+s},&\Phi''(s)&=-\lambda_+^2e^{-\lambda_+s}&&(s>0),\\
\Phi'(s)&=\kappa\lambda_-e^{\lambda_-s},&\Phi''(s)&=\kappa\lambda_-^2e^{\lambda_-s}&&(s<0).
\end{align*}
The two one sided limits of $\Phi'$ at $0$ are $\lambda_+$ and $\kappa\lambda_-=\lambda_+$. Since $\Phi''$ is also bounded, $\Phi$ is of class $C^1$ with a Lipschitz derivative. We also have $\Phi(0)=0$, $0\le\Phi<1$ on $[0,\infty)$ and $-\kappa<\Phi<0$ on $(-\infty,0)$, so that $-\kappa\le\Phi\le1$. Let $\sigma^2\in[6T,8T]$. For $s>0$ and $\mathsf b\ge\Delta/2$ we use $\sigma^2/2\le4T$:
\[
\frac{\sigma^2}2\,\Phi''(s)+\mathsf b\,\Phi'(s)=\lambda_+e^{-\lambda_+s}\Big(\mathsf b-\frac{\sigma^2}2\lambda_+\Big)\ge\lambda_+e^{-\lambda_+s}\Big(\frac\Delta2-4T\cdot\frac{\Delta}{8T}\Big)=0 .
\]
For $s<0$ and $\mathsf b\ge-2\Delta$ we use $\sigma^2/2\ge3T$:
\[
\frac{\sigma^2}2\,\Phi''(s)+\mathsf b\,\Phi'(s)=\kappa\lambda_-e^{\lambda_-s}\Big(\frac{\sigma^2}2\lambda_-+\mathsf b\Big)\ge\kappa\lambda_-e^{\lambda_-s}\Big(3T\cdot\frac{2\Delta}{3T}-2\Delta\Big)=0 .
\]
The It\^o formula holds for $\Phi$, because $\Phi'$ is Lipschitz and $d\langle\mathsf g\rangle=d\langle \mathsf N\rangle$ is absolutely continuous. With $\sigma_u^2=2T|\nabla g(X_u)|^2$ it gives
\[
\Phi(\mathsf g_{t\wedge\varsigma})=\int_0^{t\wedge\varsigma}\Big(\frac{\sigma_u^2}2\,\Phi''(\mathsf g_u)+\mathsf b_u\,\Phi'(\mathsf g_u)\Big)\,du+\int_0^{t\wedge\varsigma}\Phi'(\mathsf g_u)\,d\mathsf N_u .
\]
By \eqref{eq:ridgeqv} we have $\sigma_u^2\in[6T,8T]$. By \eqref{eq:ridgedrift} we have $\mathsf b_u\ge\Delta/2$ when $\mathsf g_u>0$ and $\mathsf b_u\ge-2\Delta$ when $\mathsf g_u<0$. So the two computations above show that the first integrand is nonnegative for almost every $u$. The second integral is a martingale, because $0\le\Phi'\le\lambda_+$ and $d\langle \mathsf N\rangle_u\le8T\,du$. Hence the process $\Phi(\mathsf g_{t\wedge\varsigma})$ is a bounded submartingale. Since $\varsigma\le t_3$, we may take $t=t_3$ and take expectations. We get
\[
\E\,\Phi(\mathsf g_\varsigma)\ge\Phi(\mathsf g_0)=\Phi(0)=0 .
\]
Let $H_+=\{\mathsf g_\varsigma=\mathsf r_1\}$ and $H_-=\{\mathsf g_\varsigma=-\mathsf r_1\}$. By Step 3 the complement of $H_+\cup H_-$ is contained in $(H_1\cap H_2)^c$. We write $q_0=\P((H_1\cap H_2)^c)$, so that
\[
\P(H_+)+\P(H_-)\ge1-q_0,\qquad q_0\le\P(H_1^c)+\P(H_2^c)\le\frac{8Tt_3}{\mathsf r_1^2}+\frac{18\,d\,T\,t_3}{\mathsf r_0^2} .
\]
Now $\Phi\le1$ everywhere, and
\[
\Phi(-\mathsf r_1)=-\kappa\,\big(1-e^{-\lambda_-\mathsf r_1}\big)\le-\frac\kappa2
\]
when $T$ is so small that $e^{-\lambda_-\mathsf r_1}\le\frac12$. We split the expectation of $\Phi(\mathsf g_\varsigma)$ over the events $H_+$, $H_-$ and the complement of their union:
\[
0\le\E\,\Phi(\mathsf g_\varsigma)\le\P(H_+)-\frac\kappa2\,\P(H_-)+q_0\le\P(H_+)-\frac\kappa2\,\big(1-\P(H_+)-q_0\big)+q_0 .
\]
The right side is equal to $(1+\frac\kappa2)\,\P(H_+)-\frac\kappa2+(1+\frac\kappa2)\,q_0$. Dividing by $1+\frac\kappa2$ and using $\frac\kappa2=\frac3{32}$ and $1+\frac\kappa2=\frac{35}{32}$, we obtain
\[
\P(H_+)\ge\frac{\kappa/2}{1+\kappa/2}-q_0=\frac3{35}-q_0 .
\]
Choose $T_1$ so small that for $T<T_1$ all the conditions on $T$ above hold and $q_0\le1/35$. Then $\P(H_+)\ge2/35$. On $H_+$ we have $\varsigma\le\varsigma_0$, so that $X_\varsigma\in\mathfrak B$. Also $g(X_\varsigma)=\mathsf r_1$. This means that $X_\varsigma\in K_{ij}$. Since $\varsigma\le t_3$, this proves (b).
\end{proof}

\begin{lemma}\label{lem:attempt}
Let $\{i,j\}$ be a direct pair and $\eta>0$. There are numbers $S\ge0$, $t_0>0$, $p_0>0$ and $T_0>0$ such that for all $T<T_0$ and all $y\in\widetilde U_i$
\[
\P\Big(X^y\text{ reaches }\widetilde U_j\text{ before the time }\frac ST+t_0\Big)\ge p_0\exp\Big(-\frac{\Delta_{ij}^2/8+\eta}{T}\Big).
\]
\end{lemma}

\begin{proof}
We use the notation introduced before Lemma \ref{lem:ridge}, and $\mathsf r_0$, $\mathsf r_1$, $t_3$, $K_{ij}$, $T_1$ are as in that lemma. We write $D_i=\diam C_i$.

\emph{Step 0: parameters and the path.} By \eqref{eq:direct}, for $0\le\lambda\le\Delta/2$ and $l\notin\{i,j\}$ we have $d_{C_l}(w_\lambda)>\min\{\lambda,\Delta-\lambda\}=\lambda$. The function $\lambda\mapsto d_{C_l}(w_\lambda)-\lambda$ is continuous and positive on the compact interval $[0,\Delta/2]$. There are finitely many indices $l$. So, by Lemma \ref{lem:direct}(iv) and by compactness,
\[
m_1:=\min\big\{d_{C_l}(w_\lambda)-\lambda:\ 0\le\lambda\le\Delta/2,\ l\notin\{i,j\}\big\}>0
\]
(we put $m_1=\infty$ if $m=2$). Choose $\eps_1\in(0,\min\{\mathsf r_0/12,\delta_*\}]$, with $\eps_1\le D_i$ in addition if $k_i\ge1$, and $\eps_2\in(0,\min\{1,\Delta/2\}]$ such that
\begin{equation}\label{eq:epschoice}
\Big(\frac\Delta2+\eps_1\Big)^2+\eps_2^2\le\frac{\Delta^2}4+\eta .
\end{equation}
This is possible, because the left side tends to $\Delta^2/4$ when $\eps_1$ and $\eps_2$ tend to zero. We put
\[
t_B=1+\log\frac{\Delta/2+\eps_1}{\eps_2},\qquad\varpi=\min\Big\{\frac{\eps_1}2,\ \frac{m_1}4\Big\},\qquad r=\frac\varpi2\,e^{-t_B} .
\]
Since $\eps_2\le\Delta/2$, we have $t_B>1$. Define $\psi:[0,t_B]\to\R^d$ by
\[
\psi(u)=a+u\,\eps_2\,e\quad(0\le u\le1),\qquad\psi(u)=a+\eps_2\exp(u-1)\,e\quad(1\le u\le t_B).
\]
Then $\psi(0)=a$ and
\[
\psi(t_B)=a+\eps_2\,\frac{\Delta/2+\eps_1}{\eps_2}\,e=a+\Big(\frac\Delta2+\eps_1\Big)e=w_{\Delta/2+\eps_1} .
\]
Every value of $\psi$ is a point $w_\lambda$ with $0\le\lambda\le\Delta/2+\eps_1$. Since $\eps_1\le\mathsf r_0/12<\Delta/2$, these $\lambda$ are smaller than $\Delta$. Hence $P_i\psi(u)=a$ by Lemma \ref{lem:direct}(iv).

We compute the energy \eqref{eq:energy} of $\psi$ for the set $C'=C_i$. On $[0,1]$ we have $\psi'=\eps_2e$ and $\psi-P_i\psi=u\,\eps_2e$. So $|\psi'|=\eps_2$ and $|\psi-P_i\psi|\le\eps_2$. Hence $|\psi'+\psi-P_i\psi|\le2\eps_2$. The energy of this piece is therefore at most
\[
\frac12\int_0^1(2\eps_2)^2\,du=2\eps_2^2 .
\]
On $[1,t_B]$ we have $\psi'=\eps_2\exp(u-1)\,e=\psi-a$. So $\psi'+\psi-P_i\psi=2(\psi-a)$ and $|\psi'+\psi-P_i\psi|^2=4\eps_2^2\exp(2(u-1))$. The energy of this piece is
\begin{align*}
\frac12\int_1^{t_B}4\,\eps_2^2\exp(2(u-1))\,du&=\eps_2^2\big(\exp(2(t_B-1))-1\big)\\
&=\eps_2^2\Big(\frac{(\Delta/2+\eps_1)^2}{\eps_2^2}-1\Big)=\Big(\frac\Delta2+\eps_1\Big)^2-\eps_2^2 .
\end{align*}
Adding the two pieces and using \eqref{eq:epschoice},
\begin{equation}\label{eq:psienergy}
J^{C_i}_{t_B}(\psi)\le2\eps_2^2+\Big(\frac\Delta2+\eps_1\Big)^2-\eps_2^2=\Big(\frac\Delta2+\eps_1\Big)^2+\eps_2^2\le\frac{\Delta^2}4+\eta .
\end{equation}

\emph{Step 1: reaching a neighbourhood of the point $a$.} We claim that there are $S\ge0$, $s_1\ge0$, $p_1>0$ and $T_2>0$ such that, with $t_A=S/T+s_1$, for all $T<T_2$ and $y\in\widetilde U_i$ the event
\[
G_A=\big\{X^y_t=X^{i,y}_t\in U_i\text{ for all }t\le t_A,\ |X^y_{t_A}-a|<r\big\}
\]
has probability at least $p_1$. First we show that $X^y$ and $X^{i,y}$ coincide up to the time $t_A$ with large probability. We have $d_{C_i}(y)\le\delta_*/4\le\delta_i/4$, so that $\delta_i/2-d_{C_i}(y)\ge\delta_i/4$. By Lemma \ref{lem:confine}(i) and its proof,
\[
\P\big(X^y_t=X^{i,y}_t\in U_i\text{ for all }t\le t_A\big)\ge1-\P\Big(\sup_{t\le t_A}|\zeta_t|\ge\frac{\delta_i}4\Big).
\]
The last probability tends to zero as $T\to0$ by Proposition \ref{prop:oumax}, because $t_A$ grows only like $1/T$. Indeed, we apply \eqref{eq:oumax} with $N=t_A$, with the radius $\delta_i/4$, with the number $\frac12$ in the place of $\eta$ and with $\tr Q=d$, $\|Q\|=1$. This gives
\[
\P\Big(\sup_{t\le t_A}|\zeta_t|\ge\frac{\delta_i}4\Big)\le\big(64\,t_A+1\big)\bigg[\exp\Big(-\frac{\big(\frac{\delta_i}8-\sqrt{dT}\big)_+^2}{2T}\Big)+\exp\Big(\frac d2-\frac{\delta_i^2}{32T}\Big)\bigg].
\]
For $T\le\delta_i^2/(256d)$ we have $\sqrt{dT}\le\delta_i/16$, and then the first exponential is at most
\[
\exp\Big(-\frac1{2T}\Big(\frac{\delta_i}{16}\Big)^2\Big)=\exp\Big(-\frac{\delta_i^2}{512\,T}\Big).
\]
Since $t_A=S/T+s_1$, the right side of the bound tends to zero. It remains to bound $\P(|X^{i,y}_{t_A}-a|<r)$ from below. There are two cases.

\emph{The case $k_i=0$.} Here $C_i=\{a\}$. Then $X^{i,y}_t-a=e^{-t}(y-a)+\zeta_t$, because both sides solve the linear equation $dY=-Y\,dt+\sqrt{2T}\,dW$ with $Y_0=y-a$. Take $S=0$ and $s_1=\log_+(\delta_*/(2r))$, so that $t_A=s_1$. We have $|y-a|=d_{C_i}(y)\le\delta_*/4$. If $\delta_*\ge2r$, then $e^{-s_1}|y-a|\le\frac{2r}{\delta_*}\cdot\frac{\delta_*}4=\frac r2$. If $\delta_*<2r$, then $s_1=0$ and $|y-a|\le\delta_*/4<r/2$. So $e^{-s_1}|y-a|\le r/2$ in both cases. By \eqref{eq:oucov} we have $\E|\zeta_{s_1}|^2=dT(1-e^{-2s_1})\le dT$. The Chebyshev inequality then gives
\[
\P\big(|X^{i,y}_{s_1}-a|<r\big)\ge\P\Big(|\zeta_{s_1}|<\frac r2\Big)\ge1-\frac{4dT}{r^2} .
\]
This is at least $\frac12$ for small $T$.

\emph{The case $k_i\ge1$.} Take $s_1=0$, so that $t_A=S/T$, and let $S>0$. We apply Theorem \ref{thm:diagram}(b) to the convex set $C_i$ and the initial law $\delta_y$, for which $m_0=d_{C_i}(y)\le\delta_*/4$, with $s=S$. This is allowed when $T$ is so small that $T\le1/e$ and $T\log(1/T)\le S$. The limit law in that theorem is $\mu_S[\delta_{P_iy}]$. Lemma \ref{lem:ergodic} bounds its distance from $u_{C_i}$. By the triangle inequality,
\[
\Wt\big(\Law(X^{i,y}_{S/T}),u_{C_i}\big)\le\Psi_T(S)^{1/2}+2D_i\exp\Big(-\frac{S}{(k_i+4)c_*D_i^2}\Big),
\]
where $\Psi_T(S)$ is the right side of Theorem \ref{thm:diagram}(b) with $s=S$ and $m_0=\delta_*/4$. This right side increases with $m_0$. For fixed $S$ it tends to zero as $T\to0$. Squaring and using $(\alpha+\beta)^2\le2\alpha^2+2\beta^2$,
\[
\Wt\big(\Law(X^{i,y}_{S/T}),u_{C_i}\big)^2\le2\,\Psi_T(S)+8D_i^2\exp\Big(-\frac{2S}{(k_i+4)c_*D_i^2}\Big).
\]
Let $(\xi,\xi')$ be an optimal coupling of $\Law(X^{i,y}_{S/T})$ and $u_{C_i}$, so that $\E|\xi-\xi'|^2$ is the square of the distance on the left. We have $r/2\le\eps_1/8\le D_i$, so \eqref{eq:volA} applies with the radius $r/2$ and gives $u_{C_i}(\bar B(a,r/2))\ge(r/(2D_i))^{k_i}$. If $|\xi'-a|\le r/2$ and $|\xi-\xi'|<r/2$, then $|\xi-a|<r$. Therefore, by the Chebyshev inequality,
\begin{align*}
\P(|\xi-a|<r)&\ge\P\big(|\xi'-a|\le r/2\big)-\P\big(|\xi-\xi'|\ge r/2\big)\\
&\ge\Big(\frac{r}{2D_i}\Big)^{k_i}-\frac4{r^2}\,\Wt\big(\Law(X^{i,y}_{S/T}),u_{C_i}\big)^2 .
\end{align*}
Choose first $S$ and then $T_2$ so that the last term is at most half of the first one for $T<T_2$. For this it is enough that
\[
\frac4{r^2}\cdot8D_i^2\exp\Big(-\frac{2S}{(k_i+4)c_*D_i^2}\Big)\le\frac14\Big(\frac{r}{2D_i}\Big)^{k_i}
\]
and that, for this $S$,
\[
\frac4{r^2}\cdot2\,\Psi_T(S)\le\frac14\Big(\frac{r}{2D_i}\Big)^{k_i}\qquad(T<T_2).
\]
Then $\P(|X^{i,y}_{S/T}-a|<r)\ge\frac12(r/(2D_i))^{k_i}$ for all $T<T_2$ and all $y\in\widetilde U_i$.

\emph{End of Step 1.} The event $G_A$ contains the intersection of the events $\{X^y_t=X^{i,y}_t\in U_i\text{ for all }t\le t_A\}$ and $\{|X^{i,y}_{t_A}-a|<r\}$. Hence
\[
\P(G_A)\ge\P\big(|X^{i,y}_{t_A}-a|<r\big)-\P\Big(\sup_{t\le t_A}|\zeta_t|\ge\frac{\delta_i}4\Big).
\]
In both cases we decrease $T_2$ so that the last probability is at most half of the lower bound found for the first one. The claim follows, with $p_1=\frac14(r/(2D_i))^{k_i}$ if $k_i\ge1$ and $p_1=\frac14$ if $k_i=0$.

\emph{Step 2: climbing to the ridge.} Let $y'$ be a point with $|y'-a|<r$. We have $r\le\varpi/2\le\eps_1/4<\delta_i/2$, so that $d_{C_i}(y')<\delta_i/2$ and $y'\in U_i$. We write $X=X^{y'}$, $X^i=X^{i,y'}$ and $\tau$ for the exit time of $X$ from $B_i$. By \eqref{eq:coincide}, $\tau$ is also the exit time of $X^i$ from $B_i$ and $X=X^i$ on $[0,\tau]$. We apply Lemma \ref{lem:tube} with $C'=C_i$, $\phi=\psi$, $t=t_B$ and the starting point $y'$. Since $\psi(0)=a$ and $r=\frac\varpi2e^{-t_B}$, the radius of the tube in that lemma is
\[
\big(|y'-a|+r\big)e^{t_B}\le2r\,e^{t_B}=\varpi .
\]
By the Doob inequality, as in Step 3 of the proof of Lemma \ref{lem:ridge},
\[
\P\Big(\sup_{u\le t_B}|W_u|\ge\frac r{\sqrt{2T}}\Big)\le\frac{2T}{r^2}\,\E|W_{t_B}|^2=\frac{2\,d\,T\,t_B}{r^2}\le\frac12\qquad\text{for }T<T_3:=\frac{r^2}{4\,d\,t_B} .
\]
By \eqref{eq:psienergy} we have $J^{C_i}_{t_B}(\psi)/(2T)\le(\Delta^2/8+\eta/2)/T$. So Lemma \ref{lem:tube} shows that the event
\[
G_B=\Big\{\sup_{u\le t_B}|X^{i}_u-\psi(u)|\le\varpi\Big\}\qquad\text{satisfies}\qquad\P(G_B)\ge\frac12\exp\Big(-\frac{\Delta^2/8+\eta/2}{T}\Big)
\]
for $T<T_3$.

We show that the tube forces the exit from $B_i$. On $G_B$ we have $|X^i_{t_B}-w_{\Delta/2+\eps_1}|\le\varpi<\eps_1$. Let $x'$ be a point with $|x'-w_{\Delta/2+\eps_1}|\le\varpi$. By \eqref{eq:wlambda} we have $d_{C_i}(w_{\Delta/2+\eps_1})=\Delta/2+\eps_1$ and $d_{C_j}(w_{\Delta/2+\eps_1})=\Delta/2-\eps_1$. Since the distance functions are nonexpansive and $2\varpi\le\eps_1<2\eps_1$,
\[
d_{C_i}(x')\ge\frac\Delta2+\eps_1-\varpi>\frac\Delta2-\eps_1+\varpi\ge d_{C_j}(x') .
\]
On the closure of $B_i$ we have $d_{C_i}\le d_{C_j}$. So $x'$ does not belong to the closure of $B_i$. On $G_B$ this applies to $x'=X^i_{t_B}$. The path $X^i$ starts at $y'\in B_i$ and is outside the closure of $B_i$ at the time $t_B$. Hence it leaves $B_i$ before the time $t_B$, that is, $\tau<t_B$ on $G_B$. Since $X=X^i$ on $[0,\tau]$, the event $G_B$ is contained in the event
\[
G_B'=\{\tau<t_B\}\cap\Big\{\sup_{u\le\tau}|X_u-\psi(u)|\le\varpi\Big\},
\]
which belongs to the family of events prior to $\tau$. Hence $\P(G_B')\ge\P(G_B)$.

Next we locate $X_\tau$ on $G_B'$. Let $\lambda\in[0,\Delta/2+\eps_1]$ be such that $\psi(\tau)=w_\lambda$, so that $|X_\tau-w_\lambda|\le\varpi$. By \eqref{eq:wlambda} this gives
\[
d_{C_i}(X_\tau)\le\lambda+\varpi,\qquad d_{C_j}(X_\tau)\ge\Delta-\lambda-\varpi .
\]
Since $X_\tau$ is on the boundary of $B_i$, there is $l\ne i$ with $d_{C_l}(X_\tau)=d_{C_i}(X_\tau)$. Suppose first that $\lambda\le\Delta/2$. If $l\notin\{i,j\}$, then by the definition of $m_1$ and since $2\varpi\le m_1/2<m_1$,
\[
d_{C_l}(X_\tau)\ge d_{C_l}(w_\lambda)-\varpi\ge\lambda+m_1-\varpi>\lambda+\varpi\ge d_{C_i}(X_\tau),
\]
which is impossible. So $l=j$, and
\[
\Delta-\lambda-\varpi\le d_{C_j}(X_\tau)=d_{C_i}(X_\tau)\le\lambda+\varpi ,
\]
that is, $\lambda\ge\Delta/2-\varpi$. Then $|w_\lambda-z|=\Delta/2-\lambda\le\varpi\le\eps_1$. Suppose now that $\lambda>\Delta/2$. Then $|w_\lambda-z|=\lambda-\Delta/2\le\eps_1$. In both cases $|w_\lambda-z|\le\eps_1$. Since $\varpi\le\eps_1$ and $\eps_1\le\mathsf r_0/12$, it follows that
\[
|X_\tau-z|\le\varpi+\eps_1\le2\eps_1\le\frac{\mathsf r_0}6 .
\]
The point $X_\tau$ is on the boundary of the open set $B_i$, so it is not in $B_i$. It is not in another basin $B_l$ either, because $B_l$ is open and does not meet $B_i$. So it belongs to no basin, that is, $X_\tau\in\cR$. Lemma \ref{lem:ridge}(a) gives $g(X_\tau)=0$.

\emph{Step 3: leaving the ridge.} On $G_B'$ the point $X_\tau$ satisfies the two conditions of Lemma \ref{lem:ridge}(b). We use the strong Markov property at the time $\tau$ and that lemma. For $T<T_1$ and conditionally on $G_B'$, the state then reaches $K_{ij}$ before the time $\tau+t_3<t_B+t_3$ with probability at least $2/35$.

\emph{Step 4: descending to $C_j$.} Let $t_C=\log(5\Delta/\delta_*)$, which is positive since $\Delta\ge\delta_*$. Let $x'\in K_{ij}$. By \ref{P5} and Lemma \ref{lem:ridge}(a),
\[
d_{C_j}\big(\varphi^j_{t_C}(x')\big)=e^{-t_C}d_{C_j}(x')\le\frac{\delta_*}{5\Delta}\cdot\frac{5\Delta}8=\frac{\delta_*}8 .
\]
The set $\{\varphi^j_u(x'):x'\in K_{ij},\ 0\le u\le t_C\}$ is compact, because it is the image of the compact set $K_{ij}\times[0,t_C]$ under a continuous map. By Lemma \ref{lem:geometry}(ii) it is contained in the open set $B_j$. Let $\mathsf d_0>0$ be smaller than half of its distance from the complement of $B_j$ and smaller than $\delta_*/8$. By Lemma \ref{lem:ballistic} and the Chebyshev inequality,
\[
\P\Big(\sup_{u\le t_C}|X^{j,x'}_u-\varphi^j_u(x')|\ge\mathsf d_0\Big)\le\frac1{\mathsf d_0^2}\,\E\sup_{u\le t_C}|X^{j,x'}_u-\varphi^j_u(x')|^2\le\frac{8dTt_C}{\mathsf d_0^2}\le\frac12\qquad(T<T_4),
\]
where $T_4=\mathsf d_0^2/(16\,d\,t_C)$. On the complementary event $X^{j,x'}$ stays in $B_j$ up to the time $t_C$. Hence $X^{x'}=X^{j,x'}$ on $[0,t_C]$ by \eqref{eq:coincide} for the index $j$, and
\[
d_{C_j}(X^{x'}_{t_C})\le d_{C_j}\big(\varphi^j_{t_C}(x')\big)+\mathsf d_0\le\frac{\delta_*}8+\frac{\delta_*}8=\frac{\delta_*}4 ,
\]
that is, $X^{x'}_{t_C}\in\widetilde U_j$. So from every point of $K_{ij}$ the state is in $\widetilde U_j$ after the time $t_C$ with probability at least $\frac12$, for $T<T_4$.

We now put the four steps together. Let $T_0=\min\{T_1,T_2,T_3,T_4\}$ and $t_0=s_1+t_B+t_3+t_C$. Let $T<T_0$ and $y\in\widetilde U_i$. We follow the state through four successive stages. By Step 1, with probability at least $p_1$ the state at the time $t_A=S/T+s_1$ is a point $y'$ with $|y'-a|<r$. We restart at the time $t_A$ and use the Markov property. Step 2 shows that, with conditional probability at least $\frac12\exp(-(\Delta^2/8+\eta/2)/T)$, the state then leaves $B_i$ within a further time $t_B$, at a point of the ridge at distance at most $\mathsf r_0/6$ from $z$. At this exit time we use the strong Markov property and Step 3. With conditional probability at least $2/35$ the state then reaches $K_{ij}$ within a further time $t_3$. Finally we use the strong Markov property at the first hitting time of $K_{ij}$ after the exit. Step 4 shows that, with conditional probability at least $\frac12$, the state is in $\widetilde U_j$ after a further time $t_C$. The total time is smaller than $t_A+t_B+t_3+t_C=S/T+t_0$. Multiplying the four bounds,
\begin{align*}
\P\Big(X^y\text{ reaches }\widetilde U_j\text{ before }\frac ST+t_0\Big)&\ge p_1\cdot\frac12\exp\Big(-\frac{\Delta^2/8+\eta/2}{T}\Big)\cdot\frac2{35}\cdot\frac12\\
&=\frac{p_1}{70}\exp\Big(-\frac{\Delta^2/8+\eta/2}{T}\Big).
\end{align*}
Since $\eta/2<\eta$, the right side is at least $\frac{p_1}{70}\exp(-(\Delta^2/8+\eta)/T)$, which is the claim with $p_0=p_1/70$.
\end{proof}

\begin{proof}[Proof of Theorem \ref{thm:exit}]
\emph{Lower bound.} Let $\eta>0$. The event in the theorem becomes larger when $\eta$ increases, so we may assume that $\eta<\delta_i^2/8$. Let $c=\delta_i^2/8-\eta$. The proof of Lemma \ref{lem:confine}(ii) shows that, with probability tending to one, the process $X^{i,x}$ stays in $B_i$ up to the time $e^{c/T}$. By \eqref{eq:coincide} the same is then true for $X^x$. So $\P(\tau_i\ge e^{c/T})\to1$. Recall also that $\theta_i\ge\tau_i$.

\emph{Upper bound for a start in $\widetilde U_i$.} Let $j$ be an index with $\Delta_{ij}=\delta_i$. The pair $\{i,j\}$ is direct by Lemma \ref{lem:direct}(i). Let $S,t_0,p_0,T_0$ be given by Lemma \ref{lem:attempt} for this pair and for $\eta/2$, and let
\begin{gather*}
p_T=p_0\exp\Big(-\frac{\delta_i^2/8+\eta/2}{T}\Big),\qquad L=2\delta_i,\qquad M=\{y:\ d_C(y)\le L\},\\
t_D=\log\frac{16L}{\delta_*},\qquad\ell=t_D+\frac ST+t_0 .
\end{gather*}
We first show that one round of length $\ell$ succeeds with probability at least $\frac78p_T$ from every point of $M$:
\begin{equation}\label{eq:cycle}
\P\big(X^y\text{ reaches }\bigcup_{l\ne i}\widetilde U_l\text{ before the time }\ell\big)\ge\frac78\,p_T\qquad(y\in M,\ T<T_5),
\end{equation}
where $T_5=\min\{T_0,\delta_*^2/(256d)\}$. To see this, let $y\in M$. By the Chebyshev inequality and the first inequality in \eqref{eq:msq},
\[
\P\Big(d_C(X^y_{t_D})>\frac{\delta_*}4\Big)\le\frac{16}{\delta_*^2}\,\E\,d_C(X^y_{t_D})^2\le\frac{16}{\delta_*^2}\big(e^{-2t_D}L^2+dT\big).
\]
By the choice of $t_D$ we have $e^{-2t_D}L^2=(\delta_*/(16L))^2L^2=\delta_*^2/256$. For $T\le\delta_*^2/(256d)$ we also have $16dT/\delta_*^2\le1/16$. Hence
\[
\P\Big(d_C(X^y_{t_D})>\frac{\delta_*}4\Big)\le\frac{16}{\delta_*^2}\big(e^{-2t_D}L^2+dT\big)=\frac1{16}+\frac{16dT}{\delta_*^2}\le\frac18
\]
for $T<T_5$. Suppose that $d_C(X^y_{t_D})\le\delta_*/4$. Since $d_C=\min_ld_{C_l}$, the point $X^y_{t_D}$ belongs to one of the sets $\widetilde U_l$. If $l\ne i$, the state has reached the union at the time $t_D<\ell$. If $l=i$, then by the Markov property at the time $t_D$ and Lemma \ref{lem:attempt} the state reaches $\widetilde U_j$ within a further time $S/T+t_0$ with conditional probability at least $p_T$. Since $p_T\le1$, in both cases the conditional probability of reaching the union before the time $t_D+S/T+t_0=\ell$ is at least $p_T$. So the probability in \eqref{eq:cycle} is at least $(1-\frac18)\,p_T$. This proves \eqref{eq:cycle}.

Now let $x\in\widetilde U_i$, $T<T_5$ and $u_n=\P(\theta_i>n\ell)$. Put $q_T:=\sup_{t\ge0}\P(d_C(X^x_t)>L)$. On the event $\{\theta_i>n\ell\}$ we have $\theta_i>(n+1)\ell$ only if the state, restarted at the time $n\ell$ from the point $X^x_{n\ell}$, does not reach the union before the time $\ell$. By \eqref{eq:cycle} the conditional probability of this is at most $1-\frac78p_T$ when $X^x_{n\ell}\in M$ and at most $1$ otherwise. So, by the Markov property at the time $n\ell$,
\[
u_{n+1}\le\E\Big[\one_{\{\theta_i>n\ell\}}\Big(1-\frac78\,p_T\,\one_{\{X^x_{n\ell}\in M\}}\Big)\Big]=u_n-\frac78\,p_T\,\P\big(\theta_i>n\ell,\ X^x_{n\ell}\in M\big).
\]
The last probability is at least $u_n-\P(X^x_{n\ell}\notin M)\ge u_n-q_T$. Therefore
\[
u_{n+1}\le\Big(1-\frac78\,p_T\Big)u_n+\frac78\,p_T\,q_T .
\]
We prove by induction that $u_n\le(1-\frac78p_T)^n+q_T$. For $n=0$ this holds because $u_0\le1$. If it holds for $n$, then
\[
u_{n+1}\le\Big(1-\frac78\,p_T\Big)^{n+1}+\Big(1-\frac78\,p_T\Big)q_T+\frac78\,p_T\,q_T=\Big(1-\frac78\,p_T\Big)^{n+1}+q_T .
\]
We bound $q_T$. We have $d_C(x)\le\delta_*/4\le\delta_i/4$, so that $L-d_C(x)\ge2\delta_i-\delta_i/4\ge\delta_i$. By the second inequality in \eqref{eq:msq}, $q_T\le\sup_t\P(|\zeta_t|\ge\delta_i)$. For $T\le\delta_i^2/(4d)$ we have $\sqrt{dT}\le\delta_i/2$. Then \eqref{eq:gaussconc} with $\tr Q=d$, $\|Q\|=1$ gives, for every $t$,
\[
\P\big(|\zeta_t|\ge\delta_i\big)\le\P\Big(|\zeta_t|\ge\sqrt{dT}+\frac{\delta_i}2\Big)\le\exp\Big(-\frac{\delta_i^2}{8T}\Big).
\]
So $q_T$ tends to zero as $T\to0$.

Take $n$ to be the integer part of $e^{(\delta_i^2/8+\eta)/T}/\ell$. Then $n\ell\le e^{(\delta_i^2/8+\eta)/T}$ and $n\ge e^{(\delta_i^2/8+\eta)/T}/\ell-1$. Since $\ell$ grows like $1/T$,
\[
n\,p_T\ge\Big(\frac{e^{(\delta_i^2/8+\eta)/T}}{\ell}-1\Big)p_0\,e^{-(\delta_i^2/8+\eta/2)/T}\longrightarrow\infty .
\]
Indeed, the right side is equal to $p_0\,e^{\eta/(2T)}/\ell-p_T$, where $p_T\to0$. For $T\le1$ we have $\ell\le(t_D+S+t_0)/T$. Since also $e^{s}\ge s^2/2$ for $s\ge0$, we get
\[
\frac{p_0\,e^{\eta/(2T)}}{\ell}\ge\frac{p_0\,T}{t_D+S+t_0}\cdot\frac{\eta^2}{8T^2}=\frac{p_0\,\eta^2}{8\,(t_D+S+t_0)\,T}\longrightarrow\infty .
\]
We now use $1-s\le e^{-s}$ and the bound for $u_n$:
\[
\P\big(\theta_i>e^{(\delta_i^2/8+\eta)/T}\big)\le\P(\theta_i>n\ell)=u_n\le\Big(1-\frac78\,p_T\Big)^n+q_T\le e^{-\frac78np_T}+q_T\longrightarrow0 .
\]
The convergence is uniform in $x\in\widetilde U_i$, because the bounds for $n\,p_T$ and for $q_T$ do not depend on $x$.

\emph{Upper bound for a start in $B_i$.} Let $x\in B_i$ and $t_1=\log_+(8d_{C_i}(x)/\delta_*)$, so that $d_{C_i}(\varphi^i_{t_1}(x))=e^{-t_1}d_{C_i}(x)\le\delta_*/8$. We argue as in the proof of Lemma \ref{lem:confine}(iii). Let $\mathsf d'\in(0,\delta_*/8]$ be smaller than the distance of the curve $\{\varphi^i_u(x):u\le t_1\}$ from the complement of $B_i$. By Lemma \ref{lem:ballistic} and the Chebyshev inequality,
\[
\P\Big(\sup_{u\le t_1}|X^{i,x}_u-\varphi^i_u(x)|\ge\mathsf d'\Big)\le\frac{8dTt_1}{\mathsf d'^2}\longrightarrow0\qquad(T\to0).
\]
On the complementary event $X^{i,x}$ stays in $B_i$ up to the time $t_1$. So $X^x=X^{i,x}$ on $[0,t_1]$ by \eqref{eq:coincide}, and
\[
d_{C_i}(X^x_{t_1})\le\frac{\delta_*}8+\mathsf d'\le\frac{\delta_*}4 ,
\]
that is, $X^x_{t_1}\in\widetilde U_i$. By the Markov property at the time $t_1$ and the previous step, which is uniform in the starting point in $\widetilde U_i$,
\[
\P\big(\theta_i>t_1+e^{(\delta_i^2/8+\eta)/T}\big)\le\P\big(X^x_{t_1}\notin\widetilde U_i\big)+\sup_{x'\in\widetilde U_i}\P\big(\theta_i^{x'}>e^{(\delta_i^2/8+\eta)/T}\big)\longrightarrow0 ,
\]
where $\theta_i^{x'}$ denotes the time $\theta_i$ for the starting point $x'$. This holds for every $\eta>0$. The number $t_1$ does not depend on $T$, so $t_1+e^{(\delta_i^2/8+\eta/2)/T}\le e^{(\delta_i^2/8+\eta)/T}$ for small $T$. Applying the last display with $\eta/2$ in the place of $\eta$ we get $\P(\theta_i>e^{(\delta_i^2/8+\eta)/T})\to0$. This proves the upper bound. The theorem follows from it, from the lower bound and from $\tau_i\le\theta_i$.
\end{proof}

\begin{remark}
The proof shows that in each attempt the state leaves the basin through a neighbourhood of the midpoint between $C_i$ and a nearest component with a probability which is at least $p_0\exp(-(\delta_i^2/8+\eta)/T)$. We have not determined the law of the exit point $X_{\tau_i}$. The proof also gives no information on the prefactor of the exit time. For flat wells we do not expect the prefactor to be a constant. By Theorem \ref{thm:gibbs} the integral of $e^{-d_C^2/2T}$ over the basin $B_i$ is of order $T^{(d-k_i)/2}$, so one expects a power of $T$ which depends on the dimension of the component.
\end{remark}

\subsection{Convergence to equilibrium}
Theorem \ref{thm:meta}(c) gives the limit when $t\to\infty$ is taken first. We now show that the law of the state is close to the equilibrium $\gamma_T^C$, in total variation, after a time $\exp(c/T)$ with an explicit number $c$. This number cannot be lowered for two components. For three components it cannot be lowered when a component $C_l$ with $\delta_l=\Delta_*$ has $w_l<1$ (Remark \ref{rem:equil}). With $\Delta_*$ as in Definition \ref{def:direct}, we put
\[
n_m=\lfloor m/2\rfloor .
\]

\begin{theorem}\label{thm:equil}
Let $C$ be as in \eqref{eq:union}, let $c>n_m\Delta_*^2/8$ and $L>0$. Then
\[
\lim_{T\to0}\ \sup_{x:\,d_C(x)\le L}\ \sup_{t\ge e^{c/T}}\TV\big(P_t(x,\cdot),\gamma_T^C\big)=0 .
\]
Hence, for every probability measure $\rho_0$ on $\R^d$,
\[
\sup_{t\ge e^{c/T}}\TV\big(\rho_{T,t},\gamma_T^C\big)\longrightarrow0\qquad(T\to0),
\]
and $\rho_{T,t_T}\to\sum_iw_iu_{C_i}$ weakly whenever $t_T\ge e^{c/T}$, with $w_i$ as in Proposition \ref{prop:gibbsunion}.
\end{theorem}

\begin{proof}
In this proof $P_t(x,\cdot)$ always denotes the transition probabilities. We write $P_{C_1}$ for the projection onto $C_1$. The supremum over the set $\{d_C\le L\}$ increases with $L$, so we may assume $L\ge\max_i\delta_i$. Since $c>n_m\Delta_*^2/8$, we can fix $\eta>0$ with $c>n_m(\Delta_*^2/8+2\eta)$. Put $M=\{y:\ d_C(y)\le2L\}$. The proof follows the plan of the proof of \cite[Theorem~5.10]{DarapaneniConvex}. We divide the time into rounds of the same length $r$. After showing that the laws at the end of one round, started from two points of $M$, overlap, we iterate with Lemma \ref{lem:doeblin}. There are two differences. First, a round now contains the transitions which bring the state to one fixed well. Second, the overlap is exponentially small in $1/T$.

\emph{Step 0: the target well and the constants.} By Lemma \ref{lem:direct}(iii) the graph whose edges are the direct pairs $\{i,j\}$ with $\Delta_{ij}\le\Delta_*$ is connected. Let $\mathcal T$ be a spanning tree of it, that is, a connected subgraph without cycles which contains all the vertices $1,\dots,m$. In a tree two vertices are joined by exactly one path with distinct vertices. Their distance is the number of edges of this path. We show that a tree with $m$ vertices has a vertex whose distance in the tree to every other vertex is at most $n_m$.

Let $v_0,v_1,\dots,v_N$ be a path in $\mathcal T$ with distinct vertices and with the largest possible number $N$ of edges. We have $N\ge1$ because $m\ge2$ and $\mathcal T$ is connected. Also $N\le m-1$, because the $N+1$ vertices of the path are distinct. Put
\[
h=\lceil N/2\rceil,\qquad i_0=v_h .
\]
The vertex $i_0$ is in the middle of the path. It divides the path into the two halves
\[
\Gamma_1=(v_0,\dots,v_h),\qquad\Gamma_2=(v_h,\dots,v_N),
\]
which have $\lceil N/2\rceil$ and $\lfloor N/2\rfloor$ edges, respectively. Only $i_0$ lies on both halves. We claim that every vertex of $\mathcal T$ is at distance at most $h$ from $i_0$. Suppose that this is false. Then there is a vertex $l$ whose distance $\mathsf n(l)$ from $i_0$ is larger than $h$. Let
\[
l=z_0,\ z_1,\ \dots,\ z_{\mathsf n(l)}=i_0
\]
be the path in $\mathcal T$ from $l$ to $i_0$. It arrives at $i_0$ through the neighbour $z_{\mathsf n(l)-1}$ of $i_0$. This vertex is different from $i_0$, so it lies on at most one of the two halves. Let $\Gamma$ be a half on which $z_{\mathsf n(l)-1}$ does not lie. The number of edges of $\Gamma$ is at least $\lfloor N/2\rfloor$. We show that the path from $l$ to $i_0$ meets $\Gamma$ only at $i_0$. Indeed, let $s<\mathsf n(l)$ and suppose that $z_s$ is a vertex of $\Gamma$. Then $z_s,z_{s+1},\dots,z_{\mathsf n(l)}$ is a path in $\mathcal T$ from $z_s$ to $i_0$. The part of $\Gamma$ between $z_s$ and $i_0$ is another one. Since there is only one path between two vertices of a tree, the two coincide. So $z_{\mathsf n(l)-1}$ lies on $\Gamma$, which contradicts the choice of $\Gamma$. Hence the path $z_0,\dots,z_{\mathsf n(l)}$ followed by the half $\Gamma$, which we run through starting from $i_0$, is a path in $\mathcal T$ with distinct vertices. Its number of edges is at least
\[
\mathsf n(l)+\lfloor N/2\rfloor>\lceil N/2\rceil+\lfloor N/2\rfloor=N .
\]
This is impossible, because $N$ is the largest number of edges of a path in $\mathcal T$. So the claim holds. Finally,
\[
h=\Big\lceil\frac N2\Big\rceil\le\Big\lceil\frac{m-1}2\Big\rceil=\Big\lfloor\frac m2\Big\rfloor=n_m .
\]
For the second equality we distinguish two cases. If $m=2b$, then $\lceil(2b-1)/2\rceil=b$. If $m=2b+1$, then $\lceil2b/2\rceil=b$.

After relabelling the components we may assume that $i_0=1$. So for every index $l$ there is a path $l=l_0,l_1,\dots,l_{\mathsf n(l)}=1$ with $\mathsf n(l)\le n_m$ such that each $\{l_{s-1},l_s\}$ is a direct pair with $\Delta_{l_{s-1}l_s}\le\Delta_*$. We fix one such path for each $l$ (with $\mathsf n(l)=0$ if $l=1$). We call the well $C_1$ the target well.

Lemma \ref{lem:attempt}, applied with the given $\eta$ to each of the finitely many pairs $(i,j)=(l_{s-1},l_s)$ on these paths, gives four numbers for each pair. We take the largest of the numbers $S$, the largest of the numbers $t_0$, the smallest of the numbers $p_0$ and the smallest of the numbers $T_0$. In this way we get numbers $S\ge0$, $t_0>0$, $p_0>0$, $T_0>0$ which work for all these pairs. Put
\[
t_*=\frac ST+t_0,\qquad p_T=p_0\exp\Big(-\frac{\Delta_*^2/8+\eta}{T}\Big).
\]
Since $\Delta_{ij}\le\Delta_*$ for these pairs, Lemma \ref{lem:attempt} gives, for each of them,
\begin{equation}\label{eq:xs5ca}
\P\big(X^{y}\text{ reaches }\widetilde U_j\text{ before the time }t_*\big)\ge p_0\exp\Big(-\frac{\Delta_{ij}^2/8+\eta}{T}\Big)\ge p_T\qquad(y\in\widetilde U_i,\ T<T_0).
\end{equation}
Let $t_D=\log(32L/\delta_*)$, $t_R=t_D+n_mt_*$ and $t_a=\log(\delta_1/\sqrt T)$. We take $\lambda_d$ as in Theorem \ref{thm:tvmix}. The time $t_c$ is defined as in the proof of \cite[Theorem~5.10]{DarapaneniConvex} for the convex set $C_1$, but with the accuracy $\frac18$ in place of $\frac14$:
\[
t_c=\frac{c_*(2D_1^2+2d)}{T}\,\log\Big(4\Big(\frac{\sqrt e\,\cZ_1(C_1)}{\omega_d}\Big)^{1/2}T^{-d/4}\Big),\qquad D_1=\diam C_1 .
\]
The time $t_D$ is positive because $L\ge\delta_*$. The time $t_a$ is positive for $T<\delta_1^2$. As in the proof of \cite[Theorem~5.10]{DarapaneniConvex}, the time $t_c$ is positive for $T\le1$. The length of a round is $r=t_R+t_a+1+t_c$, which is at most a constant multiple of $\log(1/T)/T$. Indeed,
\begin{align*}
r&=t_D+\frac{n_mS}{T}+n_mt_0+\log\delta_1+\frac12\log\frac1T+1\\
&\quad+\frac{c_*(2D_1^2+2d)}{T}\Big[\log\Big(4\Big(\frac{\sqrt e\,\cZ_1(C_1)}{\omega_d}\Big)^{1/2}\Big)+\frac d4\log\frac1T\Big].
\end{align*}
For $T\le1/e$ we have $1\le\log(1/T)$ and $1\le1/T$, so each of the terms on the right is at most a constant multiple of $\log(1/T)/T$. Hence there is a number $A_0$, which does not depend on $T$, such that
\begin{equation}\label{eq:xs5cb}
r\le\frac{A_0}T\,\log\frac1T\qquad(T\le1/e).
\end{equation}
All statements below hold for $T$ small enough, uniformly in the starting points considered.

\emph{Step 1: recovery.} Let $y\in M$, so that $d_C(y)\le2L$. As in the proof of \eqref{eq:cycle}, the first inequality in \eqref{eq:msq} gives
\[
\E\,d_C\big(X^y_{t_D}\big)^2\le e^{-2t_D}d_C(y)^2+dT\,\big(1-e^{-2t_D}\big)\le e^{-2t_D}4L^2+dT .
\]
By the definition of $t_D$ we have $e^{-t_D}=\delta_*/(32L)$. Thus
\[
e^{-2t_D}\,4L^2=\frac{\delta_*^2}{1024\,L^2}\cdot4L^2=\frac{\delta_*^2}{256} .
\]
By the Chebyshev inequality, for $T\le\delta_*^2/(256d)$,
\[
\P\Big(d_C\big(X^y_{t_D}\big)>\frac{\delta_*}4\Big)\le\frac{16}{\delta_*^2}\big(e^{-2t_D}4L^2+dT\big)=\frac{16}{\delta_*^2}\cdot\frac{\delta_*^2}{256}+\frac{16\,dT}{\delta_*^2}\le\frac1{16}+\frac1{16}=\frac18 .
\]
Now $d_C=\min_ld_{C_l}$. So, if $d_C(X^y_{t_D})\le\delta_*/4$, then $d_{C_l}(X^y_{t_D})\le\delta_*/4$ for some $l$. Thus with probability at least $\frac78$ the state at the time $t_D$ belongs to $\widetilde U_l$ for some $l$.

\emph{Step 2: routing.} Take a point $y'\in\widetilde U_l$. Let $l=l_0,\dots,l_{\mathsf n(l)}=1$ be the path fixed in Step 0. We show that there is a stopping time $\hat\tau\le n_mt_*$ such that
\begin{equation}\label{eq:xs5cc}
\P\big(X^{y'}_{\hat\tau}\in\widetilde U_1\big)\ge p_T^{\,\mathsf n(l)}\ge p_T^{\,n_m}.
\end{equation}
If $\mathsf n(l)=0$, then $l=1$ and we take $\hat\tau=0$. Let $\mathsf n(l)\ge1$. Put $\hat\tau_0=0$ and, for $1\le s\le \mathsf n(l)$, let $\hat\tau_s$ be the first time after $\hat\tau_{s-1}$ at which the state $X^{y'}$ is in $\widetilde U_{l_s}$. These are the successive hitting times of $\widetilde U_{l_1},\dots,\widetilde U_{l_{\mathsf n(l)}}$. They are stopping times. When $\hat\tau_s$ is finite, the state at the time $\hat\tau_s$ belongs to $\widetilde U_{l_s}$, because this set is closed. Let
\[
H_s=\{\hat\tau_s\le s\,t_*\}\qquad(0\le s\le \mathsf n(l)).
\]
On $H_{s-1}$ the state at the time $\hat\tau_{s-1}$ is a point of $\widetilde U_{l_{s-1}}$. We use the strong Markov property at the time $\hat\tau_{s-1}$ and \eqref{eq:xs5ca} for the pair $(l_{s-1},l_s)$. They show that the conditional probability of the event $\{\hat\tau_s-\hat\tau_{s-1}\le t_*\}$, given the events prior to $\hat\tau_{s-1}$, is at least $p_T$ on $H_{s-1}$. Since $H_{s-1}\cap\{\hat\tau_s-\hat\tau_{s-1}\le t_*\}\subset H_s$, we get
\[
\P(H_s)\ge\P\big(H_{s-1}\cap\{\hat\tau_s-\hat\tau_{s-1}\le t_*\}\big)\ge p_T\,\P(H_{s-1})\qquad(1\le s\le \mathsf n(l)).
\]
Since $\P(H_0)=1$, induction on $s$ gives $\P(H_{\mathsf n(l)})\ge p_T^{\,\mathsf n(l)}$. The number $p_T$ is a lower bound for a probability, so $p_T\le1$. As $\mathsf n(l)\le n_m$, this gives $p_T^{\,\mathsf n(l)}\ge p_T^{\,n_m}$. Now put $\hat\tau=\min\{\hat\tau_{\mathsf n(l)},\,n_mt_*\}$. This is a stopping time with $\hat\tau\le n_mt_*$. On $H_{\mathsf n(l)}$ we have $\hat\tau=\hat\tau_{\mathsf n(l)}\le \mathsf n(l)\,t_*$ and $X^{y'}_{\hat\tau}\in\widetilde U_1$, which proves \eqref{eq:xs5cc}.

\emph{Step 3: waiting and approach to $C_1$.} Let $y''\in\widetilde U_1$, so that $d_{C_1}(y'')\le\delta_*/4\le\delta_1/4$. Let $N_T$ be a time which grows at most polynomially in $1/T$. By Lemma \ref{lem:confine}(i) and Proposition \ref{prop:oumax}, with probability tending to one the process $X^{y''}$ coincides with $X^{1,y''}$ on $[0,N_T]$. We write this out. For $N>0$ put
\[
\eps_T(N)=(64N+1)\Big[\exp\Big(-\frac{\delta_1^2}{512\,T}\Big)+\exp\Big(\frac d2-\frac{\delta_1^2}{32\,T}\Big)\Big].
\]
Let $x$ be a point with $d_{C_1}(x)\le\delta_1/4$. Then $x\in U_1$ and $\delta_1/2-d_{C_1}(x)\ge\delta_1/4$, so Lemma \ref{lem:confine}(i) gives the first inequality in
\begin{equation}\label{eq:xs5cd}
\P\big(X^x_t\ne X^{1,x}_t\text{ for some }t\le N\big)\le\P\Big(\sup_{0\le t\le N}|\zeta_t|\ge\frac{\delta_1}4\Big)\le\eps_T(N)\qquad\Big(T\le\frac{\delta_1^2}{256\,d}\Big).
\end{equation}
The second inequality is Proposition \ref{prop:oumax} in $\R^d$, where $\tr Q=d$ and $\|Q\|=1$, with the radius $\delta_1/4$ and with the parameter $\frac12$ in place of the parameter $\eta$ of that proposition. The factor in front is then $16N/(1/2)^2+1=64N+1$. For $T\le\delta_1^2/(256d)$ we have $\sqrt{dT}\le\delta_1/16$. Therefore
\[
\frac12\cdot\frac{\delta_1}4-\sqrt{dT}\ge\frac{\delta_1}{16},\qquad\frac1{2T}\Big(\frac{\delta_1}{16}\Big)^2=\frac{\delta_1^2}{512\,T},\qquad\frac1{2T}\Big(\frac{\delta_1}4\Big)^2=\frac{\delta_1^2}{32\,T} .
\]
If $N_T$ is at most a constant multiple of a power of $1/T$, then $\eps_T(N_T)\to0$ as $T\to0$, because a power of $1/T$ times $e^{-\kappa/T}$ tends to zero for every $\kappa>0$. This proves the statement on $X^{y''}$ and $X^{1,y''}$.

Next, by Theorem \ref{thm:envelope} for the convex set $C_1$, for every $u\ge t_a$,
\[
d_{C_1}\big(X^{1,y''}_{u}\big)\le e^{-u}\,\frac{\delta_1}4+|\zeta_{u}|\le\frac{\sqrt T}4+|\zeta_{u}| .
\]
The first inequality is \eqref{eq:env3} with $v=0$ and $d_{C_1}(y'')\le\delta_1/4$. For the second one, $u\ge t_a$ gives
\[
e^{-u}\le e^{-t_a}=\frac{\sqrt T}{\delta_1},\qquad\text{hence}\qquad e^{-u}\,\frac{\delta_1}4\le\frac{\sqrt T}{\delta_1}\cdot\frac{\delta_1}4=\frac{\sqrt T}4 .
\]
By \eqref{eq:oucov}, $\E|\zeta_u|^2=dT(1-e^{-2u})\le dT$. So the Chebyshev inequality gives
\[
\P\big(|\zeta_u|>\sqrt{2dT}\big)\le\frac{\E|\zeta_u|^2}{2dT}\le\frac{dT}{2dT}=\frac12 .
\]
Hence $\P(|\zeta_{u}|\le\sqrt{2dT})\ge\frac12$. On the event $\{|\zeta_u|\le\sqrt{2dT}\}$ we get
\[
d_{C_1}\big(X^{1,y''}_{u}\big)\le\frac{\sqrt T}4+\sqrt{2dT}=\Big(\frac14+\sqrt{2d}\Big)\sqrt T\le\Big(\frac12+\sqrt{2d}\Big)\sqrt T=\lambda_d\sqrt T .
\]

We apply these estimates after the routing of Step 2. Let the routing start at the time $t_D$, so that it ends at the stopping time $t_D+\hat\tau\le t_D+n_mt_*=t_R$. Let $y''$ be the state at this stopping time. The time which remains until $t_R+t_a$ is
\[
u=t_R+t_a-(t_D+\hat\tau),\qquad\text{so that}\qquad t_a\le u\le t_R+t_a .
\]
The number $u$ and the point $y''$ are determined by the events prior to $t_D+\hat\tau$. After this time the state is driven by the increments of $W$, which form a Brownian motion independent of these events. We restart the envelope process from zero. So, conditionally on the events prior to $t_D+\hat\tau$, the estimates above hold with $y''$ and $u$ fixed. The time $N_T=t_R+t_a$ is at most a constant multiple of $1/T$. Clearly $u\le N_T$. By \eqref{eq:xs5cd} the state and the state process of $C_1$, both restarted from $y''$ at the time $t_D+\hat\tau$, coincide up to the time $u$ with conditional probability at least $1-\eps_T(N_T)$. The restarted envelope process is at most $\sqrt{2dT}$ at the time $u$ with conditional probability at least $\frac12$. Both events occur together with conditional probability at least
\[
1-\eps_T(N_T)-\frac12\ge\frac13 ,
\]
as soon as $T$ is so small that $\eps_T(N_T)\le\frac16$. Hence, with conditional probability at least $\frac13$, the state at the time $t_R+t_a$ is a point $y^\circ$ with
\begin{equation}\label{eq:cool}
d_{C_1}(y^\circ)\le\lambda_d\sqrt T .
\end{equation}

\emph{Step 4: a lower bound for the transition probability, and mixing in the well.} Let $y^\circ$ satisfy \eqref{eq:cool}. We denote by $\nu_{y^\circ}$ the uniform probability measure on the ball $\bar B(P_{C_1}y^\circ,\sqrt T)$, whose volume is $\omega_dT^{d/2}$. By Lemma \ref{lem:geometry}(i) the hypothesis of Lemma \ref{lem:density} holds with $C'=C_1$ and $\vartheta=\delta_1/2$. Indeed, a point $x$ with $d_{C_1}(x)<\delta_1/2$ belongs to $U_1\subset B_1$. For such a point $\Pi_C(x)=\{P_{C_1}x\}$. Let $T$ be so small that $(\lambda_d+1+8\sqrt d)\sqrt T<\delta_1/2$. Then Lemma \ref{lem:density} with $\lambda=\lambda_d$ says that, for every Borel set $A'\subset \bar B(P_{C_1}y^\circ,\sqrt T)$,
\[
P_1(y^\circ,A')\ge c_{d,\lambda_d}\,T^{-d/2}\,\Leb(A')=c_{d,\lambda_d}\,\omega_d\,\frac{\Leb(A')}{\omega_dT^{d/2}}=c_{d,\lambda_d}\,\omega_d\,\nu_{y^\circ}(A') .
\]
For an arbitrary Borel set $A$ we apply this to $A'=A\cap \bar B(P_{C_1}y^\circ,\sqrt T)$ and use $P_1(y^\circ,A)\ge P_1(y^\circ,A')$ and $\nu_{y^\circ}(A)=\nu_{y^\circ}(A')$. Therefore, for small $T$, the transition probability at time one satisfies
\[
P_1(y^\circ,\cdot)\ge c_{d,\lambda_d}\,\omega_d\,\nu_{y^\circ} .
\]
Let $P^{(1)}_t$ denote the transition probabilities of the state process of the convex set $C_1$. As in Step 3 of the proof of \cite[Theorem~5.10]{DarapaneniConvex}, with the accuracy $\frac18$, Lemma \ref{lem:L2} gives $\TV(\nu_{y^\circ}P^{(1)}_{t_c},\gamma_T^{C_1})\le\frac18$. We repeat the calculation. Let $T\le1$. The centre of the ball $\bar B(P_{C_1}y^\circ,\sqrt T)$ is a point of $C_1$, so on this ball $d_{C_1}\le\sqrt T$ and $V_{C_1}\le T/2$. Also $\cZ_T(C_1)\le\cZ_1(C_1)$. The density of $\nu_{y^\circ}$ with respect to $\gamma_T^{C_1}$ is
\[
h_0=\frac{\cZ_T(C_1)\,e^{V_{C_1}/T}}{\omega_dT^{d/2}}\,\one_{\bar B(P_{C_1}y^\circ,\sqrt T)},
\]
and it satisfies
\[
\|h_0-1\|_{L^2(\gamma_T^{C_1})}^2\le\int h_0^2\,d\gamma_T^{C_1}=\frac{\cZ_T(C_1)}{(\omega_dT^{d/2})^2}\int_{\bar B(P_{C_1}y^\circ,\sqrt T)}e^{V_{C_1}/T}\,dx\le\frac{\sqrt e\,\cZ_1(C_1)}{\omega_d\,T^{d/2}} .
\]
By the definition of $t_c$,
\[
\exp\Big(-\frac{T\,t_c}{c_*(2D_1^2+2d)}\Big)=\frac14\Big(\frac{\sqrt e\,\cZ_1(C_1)}{\omega_d}\Big)^{-1/2}T^{d/4} .
\]
We apply Lemma \ref{lem:L2} to the convex set $C_1$. The number $M_T$ of that lemma is at most $2D_1^2+2dT\le2D_1^2+2d$. Thus
\[
\TV\big(\nu_{y^\circ}P^{(1)}_{t_c},\gamma_T^{C_1}\big)\le\frac12\,\exp\Big(-\frac{T\,t_c}{c_*(2D_1^2+2d)}\Big)\Big(\frac{\sqrt e\,\cZ_1(C_1)}{\omega_d\,T^{d/2}}\Big)^{1/2}=\frac12\cdot\frac14=\frac18 .
\]

Each point of $\bar B(P_{C_1}y^\circ,\sqrt T)$ is at distance at most $\sqrt T$ from $C_1$. So, by Lemma \ref{lem:confine}(i) and Proposition \ref{prop:oumax}, the process started from $\nu_{y^\circ}$ coincides up to the time $t_c$ with the state process of $C_1$ with probability tending to one. More precisely, let $T$ be so small that $\sqrt T\le\delta_1/4$ and $T\le\delta_1^2/(256d)$. Then every point $x$ of the ball satisfies $d_{C_1}(x)\le\delta_1/4$. Thus \eqref{eq:xs5cd} holds for it with $N=t_c$. For every Borel set $A$ this gives
\[
\big|\nu_{y^\circ}P_{t_c}(A)-\nu_{y^\circ}P^{(1)}_{t_c}(A)\big|\le\int\P\big(X^x_{t_c}\ne X^{1,x}_{t_c}\big)\,\nu_{y^\circ}(dx)\le\eps_T(t_c) .
\]
The time $t_c$ is at most a constant multiple of $\log(1/T)/T$, hence of $T^{-2}$. So $\eps_T(t_c)\to0$. In particular $\eps_T(t_c)\le\frac18$ for small $T$. By the triangle inequality,
\[
\TV\big(\nu_{y^\circ}P_{t_c},\gamma_T^{C_1}\big)\le\TV\big(\nu_{y^\circ}P_{t_c},\nu_{y^\circ}P^{(1)}_{t_c}\big)+\TV\big(\nu_{y^\circ}P^{(1)}_{t_c},\gamma_T^{C_1}\big)\le\frac18+\frac18 .
\]
Hence, for small $T$,
\begin{equation}\label{eq:wellmix}
\TV\big(\nu_{y^\circ}P_{t_c},\gamma_T^{C_1}\big)\le\frac14 .
\end{equation}

\emph{Step 5: one round.} Let $y\in M$. We define an event $G$ on which Steps 1, 2 and 3 succeed. Let $G_1$ be the event that $d_C(X^y_{t_D})\le\delta_*/4$. By Step 1, $\P(G_1)\ge\frac78$. On $G_1$ the state at the time $t_D$ belongs to exactly one of the sets $\widetilde U_l$, because $\widetilde U_l\subset B_l$ by Lemma \ref{lem:geometry}(i) and the basins are pairwise disjoint. Let $G_2\subset G_1$ be the event that, in addition, the routing of Step 2, started at the time $t_D$ from the state at this time, brings the state into $\widetilde U_1$ at the time $t_D+\hat\tau$. By the Markov property at the time $t_D$ and by \eqref{eq:xs5cc}, $\P(G_2)\ge p_T^{\,n_m}\,\P(G_1)$. Let $G\subset G_2$ be the event that, in addition, the state at the time $t_R+t_a$ satisfies \eqref{eq:cool}. By the strong Markov property at the time $t_D+\hat\tau$ and by Step 3, $\P(G)\ge\frac13\,\P(G_2)$. Combining the three bounds, we get
\[
\P(G)\ge\frac13\,\P(G_2)\ge\frac13\,p_T^{\,n_m}\,\P(G_1)\ge\frac78\cdot p_T^{\,n_m}\cdot\frac13=\frac7{24}\,p_T^{\,n_m} .
\]

Whether $G$ occurs depends only on the path up to the time $t_R+t_a$. Recall that $r=(t_R+t_a)+1+t_c$. Let $A$ be a Borel set. If $y^\circ$ satisfies \eqref{eq:cool}, then by Step 4
\[
P_{1+t_c}(y^\circ,A)=\int P_1(y^\circ,dz)\,P_{t_c}(z,A)\ge c_{d,\lambda_d}\omega_d\int\nu_{y^\circ}(dz)\,P_{t_c}(z,A)=c_{d,\lambda_d}\omega_d\,\nu_{y^\circ}P_{t_c}(A) .
\]
By the Markov property at the time $t_R+t_a$,
\[
P_r(y,A)=\E\big[P_{1+t_c}\big(X^y_{t_R+t_a},A\big)\big]\ge\E\big[\one_G\,P_{1+t_c}\big(X^y_{t_R+t_a},A\big)\big].
\]
So, by the Markov property at the times $t_R+t_a$ and $t_R+t_a+1$ and by Step 4,
\[
P_r(y,\cdot)\ge\E\big[\one_G\,c_{d,\lambda_d}\omega_d\,\nu_{X^y_{t_R+t_a}}P_{t_c}\big]=c_{d,\lambda_d}\omega_d\,\P(G)\,m_y,\qquad m_y:=\E\big[\nu_{X^y_{t_R+t_a}}P_{t_c}\,\big|\,G\big],
\]
and $\TV(m_y,\gamma_T^{C_1})\le\frac14$ by \eqref{eq:wellmix}. Indeed, $m_y$ is an average of probability measures $\nu_{y^\circ}P_{t_c}$ with $y^\circ$ satisfying \eqref{eq:cool}. So, for every Borel set $A$,
\[
\big|m_y(A)-\gamma_T^{C_1}(A)\big|\le\E\Big[\big|\nu_{X^y_{t_R+t_a}}P_{t_c}(A)-\gamma_T^{C_1}(A)\big|\,\Big|\,G\Big]\le\frac14 .
\]
The same holds for another point $\tilde y\in M$. Hence $\TV(m_y,m_{\tilde y})\le\frac14+\frac14=\frac12$. Put
\[
a=\frac7{24}\,c_{d,\lambda_d}\omega_d\,p_T^{\,n_m}.
\]
By the lower bound for $\P(G)$, which holds for $y$ and for $\tilde y$, we have $P_r(y,\cdot)\ge a\,m_y$ and $P_r(\tilde y,\cdot)\ge a\,m_{\tilde y}$. Here $a\le1$, because $P_r(y,\cdot)$ and $m_y$ are probability measures. Since $1-a\,(1-\TV(m_y,m_{\tilde y}))\le1-\frac a2$, we get by Lemma \ref{lem:doeblin}(i), for all $y,\tilde y\in M$,
\[
\TV\big(P_r(y,\cdot),P_r(\tilde y,\cdot)\big)\le1-\alpha_T,\qquad\alpha_T:=\frac12\cdot\frac7{24}\,c_{d,\lambda_d}\omega_d\,p_T^{\,n_m}.
\]

\emph{Step 6: iteration.} Let $d_C(x)\le L$ and $e_n=\TV(P_{nr}(x,\cdot),\gamma_T^C)$. Since $\gamma_T^C$ is invariant (Lemma \ref{lem:invariant}), Lemma \ref{lem:doeblin}(ii) gives $e_{n+1}\le(1-\alpha_T)e_n+\alpha_T(P_{nr}(x,M^c)+\gamma_T^C(M^c))$. Here we took $\Pi=P_r$, $\alpha=\alpha_T$, $\mu=P_{nr}(x,\cdot)$ and $\mu'=\gamma_T^C$, so that $\mu\Pi=P_{(n+1)r}(x,\cdot)$ and $\mu'\Pi=\gamma_T^C$. By the second inequality in \eqref{eq:msq}, $P_t(x,M^c)\le\P(|\zeta_t|\ge L)$ for all $t$. Indeed,
\[
P_t(x,M^c)=\P\big(d_C(X^x_t)>2L\big)\le\P\big(|\zeta_t|\ge2L-d_C(x)\big)\le\P\big(|\zeta_t|\ge L\big).
\]
By Lemma \ref{lem:invariant}, $\gamma_T^C(M^c)=\lim_{t\to\infty}P_t(x,M^c)$, so the same bound holds for $\gamma_T^C(M^c)$. Put $q_T=\sup_t\P(|\zeta_t|\ge L)$. By \eqref{eq:gaussconc} with $\tr Q=d$ and $\|Q\|=1$,
\[
q_T\le\exp\Big(-\frac{(L-\sqrt{dT})_+^2}{2T}\Big),
\]
which tends to zero as $T\to0$. The recursion becomes
\[
e_{n+1}\le(1-\alpha_T)\,e_n+2\alpha_T\,q_T .
\]
We solve it by induction, starting from $e_0\le1$. If $e_n\le(1-\alpha_T)^n+2q_T$, then
\[
e_{n+1}\le(1-\alpha_T)\big[(1-\alpha_T)^n+2q_T\big]+2\alpha_T\,q_T=(1-\alpha_T)^{n+1}+2q_T .
\]
So, using $1-\alpha_T\le e^{-\alpha_T}$, for every $n\ge0$,
\[
e_n\le(1-\alpha_T)^n+2q_T\le\exp(-n\alpha_T)+2q_T .
\]
Let $n$ be the integer part of $e^{c/T}/r$. Since $r$ is at most a constant multiple of $\log(1/T)/T$ and
\[
\alpha_T=\frac7{48}\,c_{d,\lambda_d}\omega_d\,p_0^{\,n_m}\exp\Big(-\frac{n_m(\Delta_*^2/8+\eta)}{T}\Big),
\]
the choice of $\eta$ ensures that $n\alpha_T\to\infty$. To see this, write $\kappa_0=\frac7{48}c_{d,\lambda_d}\omega_d\,p_0^{\,n_m}$. By \eqref{eq:xs5cb},
\[
n\ge\frac{e^{c/T}}{r}-1\ge\frac{T\,e^{c/T}}{A_0\log(1/T)}-1 .
\]
By the choice of $\eta$, we have $c-n_m(\Delta_*^2/8+\eta)\ge n_m\eta$. Hence
\[
n\,\alpha_T\ge\kappa_0\Big[\frac{T}{A_0\log(1/T)}\exp\Big(\frac{c-n_m(\Delta_*^2/8+\eta)}{T}\Big)-1\Big]\ge\kappa_0\Big[\frac{T}{A_0\log(1/T)}\,e^{n_m\eta/T}-1\Big].
\]
Since $n_m\ge1$ and $\eta>0$, the right side tends to infinity as $T\to0$.

Because $\TV(\mu P_s,\gamma_T^C)=\TV(\mu P_s,\gamma_T^CP_s)\le\TV(\mu,\gamma_T^C)$ for every probability measure $\mu$ and every $s\ge0$, the total variation distance to the invariant law is nonincreasing in time. If $t\ge e^{c/T}$, then $t\ge nr$. Therefore
\[
\TV\big(P_t(x,\cdot),\gamma_T^C\big)\le e_n\le\exp(-n\alpha_T)+2q_T .
\]
The right side does not depend on $t$ and on the point $x$ with $d_C(x)\le L$. It tends to zero as $T\to0$, which proves the first statement.

\emph{Step 7: general initial laws.} Let $\rho_0$ be a probability measure and $L>0$. We have $\rho_{T,t}=\int P_t(x,\cdot)\,\rho_0(dx)$. So, for every Borel set $A$,
\[
\big|\rho_{T,t}(A)-\gamma_T^C(A)\big|\le\int\big|P_t(x,A)-\gamma_T^C(A)\big|\,\rho_0(dx)\le\int\TV\big(P_t(x,\cdot),\gamma_T^C\big)\,\rho_0(dx).
\]
We divide the last integral into the integral over $\{d_C\le L\}$ and the integral over $\{d_C>L\}$. On the second set we use that the total variation distance is at most one. This gives
\[
\TV\big(\rho_{T,t},\gamma_T^C\big)\le\int\TV\big(P_t(x,\cdot),\gamma_T^C\big)\,\rho_0(dx)\le\sup_{x:\,d_C(x)\le L}\TV\big(P_t(x,\cdot),\gamma_T^C\big)+\rho_0\big(d_C>L\big).
\]
We take the supremum over $t\ge e^{c/T}$ and let $T\to0$. By the first statement,
\[
\limsup_{T\to0}\ \sup_{t\ge e^{c/T}}\TV\big(\rho_{T,t},\gamma_T^C\big)\le\rho_0\big(d_C>L\big).
\]
The left side does not depend on $L$. Also $\rho_0(d_C>L)\to0$ as $L\to\infty$, because $d_C$ is finite at every point. Thus letting first $T\to0$ and then $L\to\infty$ gives the second statement.

For the third statement let $f$ be a continuous function with $|f|\le1$ and let $t_T\ge e^{c/T}$. By the formula for the total variation distance recalled before \eqref{eq:TVlsc},
\[
\Big|\int f\,d\rho_{T,t_T}-\int f\,d\gamma_T^C\Big|\le2\,\TV\big(\rho_{T,t_T},\gamma_T^C\big)\longrightarrow0\qquad(T\to0),
\]
and $\int f\,d\gamma_T^C\to\sum_iw_i\int f\,du_{C_i}$ by Proposition \ref{prop:gibbsunion}. Thus the third statement follows from the second one and from Proposition \ref{prop:gibbsunion}. This completes the proof.
\end{proof}

\begin{corollary}\label{cor:two}
Let $C=C_1\cup C_2$ with $C_1,C_2$ nonempty, compact, convex and disjoint, and let $\delta=\dist(C_1,C_2)$. Let $\rho_0$ have a finite second moment and $\rho_0(\cR)=0$, and let $t_T$ be such that $T\,t_T\to\infty$ as $T\to0$.
\begin{enumerate}[label=\textup{(\roman*)}]
\item If $t_T\le e^{c/T}$ for some $c<\delta^2/8$, then $\rho_{T,t_T}\to a_1u_{C_1}+a_2u_{C_2}$ weakly, where $a_i=\rho_0(B_i)$.
\item If $t_T\ge e^{c/T}$ for some $c>\delta^2/8$, then $\rho_{T,t_T}\to w_1u_{C_1}+w_2u_{C_2}$ weakly, where $w_1,w_2$ are the weights of Proposition \ref{prop:gibbsunion}.
\end{enumerate}
\end{corollary}

\begin{proof}
For $m=2$ we have $\delta_*=\delta_1=\delta_2=\Delta_*=\delta$ and $n_2=1$. Indeed, each of the two components has only one other component. So $\delta_1=\delta_2=\dist(C_1,C_2)$, and $\Delta_*=\delta_1$ by Lemma \ref{lem:direct}(iii). In (i) we may assume that $c>0$. This is because $e^{c/T}$ increases with $c$, so the hypothesis $t_T\le e^{c/T}$ remains true when $c$ is replaced by a larger number which is still smaller than $\delta^2/8$. Thus $0<c<\delta_*^2/8$. Also $t_T\to\infty$, because $T\,t_T\to\infty$. In (ii) we have $c>n_2\Delta_*^2/8$. So (i) is Theorem \ref{thm:meta}(b), applied along every sequence $T_n\to0$ with $s=\infty$, and (ii) is Theorem \ref{thm:equil}.
\end{proof}

\begin{remark}\label{rem:equil}
\emph{Optimality.} Let $i$ be an index with $w_i<1$ and let $c<\delta_i^2/8$. Then the conclusion of Theorem \ref{thm:equil} fails for this $c$. For $m=2$ and $m=3$ we have $n_m=1$ and $\Delta_*=\max_i\delta_i$. So for two and three components the constant $\Delta_*^2/8$ of Theorem \ref{thm:equil} cannot be lowered as soon as a component $C_i$ with $\delta_i=\Delta_*$ has $w_i<1$. Then the conclusion holds for $c>\Delta_*^2/8$ and fails for $c<\Delta_*^2/8$. For two components the condition is always satisfied, because $\delta_1=\delta_2=\Delta_*$ and $w_1+w_2=1$, so that $w_1$ and $w_2$ cannot both be equal to one. Thus for two components the weights of the two wells change from the initial weights $a_1,a_2$ to the weights $w_1,w_2$ of dimension and volume when $T\log t$ passes the value $\Delta_{12}^2/8$.

\emph{Four and more components.} For $m\ge4$ we do not expect the constant $n_m\Delta_*^2/8$ to be optimal. It comes from asking that the transitions along a path of wells succeed one after the other in a single round. By Remark \ref{rem:connected} and by analogy with the smooth case \cite{FreidlinWentzell2012,HKS1989}, we expect that the conclusion of Theorem \ref{thm:equil} holds for every $c>\Delta_*^2/8$ and every $m$.
\end{remark}

We verify the two claims of Remark \ref{rem:equil}. Let $i$ be an index with $w_i<1$ and let $c<\delta_i^2/8$. By Theorem \ref{thm:exit}, a state started at $x\in C_i$ is in $B_i$ at the time $e^{c/T}$ with probability tending to one. Indeed, choose $\eta>0$ with $c<\delta_i^2/8-\eta$. With probability tending to one the exit time $\tau_i$ from $B_i$ is at least $\exp((\delta_i^2/8-\eta)/T)$, which is larger than $e^{c/T}$. On the other hand, the closure $F$ of $B_i$ does not meet the components $C_l$, $l\ne i$, so $\limsup_{T\to0}\gamma_T^C(B_i)\le w_i$ by Proposition \ref{prop:gibbsunion}. To see this, let $l\ne i$. The open set $B_l$ contains $C_l$ and does not meet $B_i$, so it does not meet $F$. Thus $u_{C_l}(F)=0$ for $l\ne i$. Also $u_{C_i}(F)=1$, because $C_i\subset B_i$. For measures which converge weakly, the upper limit of the masses of a closed set is at most its mass for the limit. Since $\gamma_T^C\to\sum_lw_lu_{C_l}$ weakly,
\[
\limsup_{T\to0}\gamma_T^C(B_i)\le\limsup_{T\to0}\gamma_T^C(F)\le\sum_{l=1}^mw_l\,u_{C_l}(F)=w_i .
\]
Hence
\[
\TV\big(P_{e^{c/T}}(x,\cdot),\gamma_T^C\big)\ge P_{e^{c/T}}(x,B_i)-\gamma_T^C(B_i)\ge1-w_i-o(1),
\]
and the conclusion of Theorem \ref{thm:equil} fails for this $c$, because $d_C(x)=0$ and $1-w_i>0$.

It remains to prove the second claim, namely that $n_m=1$ and $\Delta_*=\max_i\delta_i$ for $m=2$ and $m=3$. For $m=2$ this was shown in the proof of Corollary \ref{cor:two}. Also $n_3=\lfloor3/2\rfloor=1$. Now let $m=3$ and let $\beta_1\le\beta_2\le\beta_3$ be the three distances $\Delta_{12},\Delta_{13},\Delta_{23}$ in increasing order. Two different edges on three vertices have a common vertex and contain all three vertices. So they form a spanning tree. One edge alone does not connect the three vertices. Hence the graph with the edges $\{i,j\}$, $\Delta_{ij}\le\Delta$, is connected exactly when at least two of the three distances are at most $\Delta$, that is, when $\Delta\ge\beta_2$. So $\Delta_*=\beta_2$. The two smallest of the three distances form a spanning tree. Let $\{i,j\}$ be a pair with $\Delta_{ij}=\beta_1$ and let $l$ be the third index. Then $\delta_i=\delta_j=\beta_1$ and $\delta_l=\min\{\Delta_{li},\Delta_{lj}\}=\beta_2$. Thus the second smallest distance is the largest of the $\delta_i$. This gives $\max_i\delta_i=\beta_2=\Delta_*$.

\section{The control problem for disjoint unions}\label{sec:hjbdisc}

We return to the value function $v_T$ of Section \ref{subsec:Bcontrol}, defined by \eqref{eq:log}. For a compact convex set the two iterated limits of $v_T$ are equal by Theorem \ref{thm:hjb}. Here $C$ is a disjoint union as in \eqref{eq:union}, with the notation of Section \ref{sec:meta}. In particular $U_i=\{d_{C_i}<\delta_i/2\}$. The process $X^{i,x}$ is the state process of the convex set $C_i$, driven by the same Brownian motion as $X^x$.

\begin{theorem}\label{thm:hjbdisc}
Let $g:\R^d\to\R$ be bounded and Lipschitz and let $v_T$ be defined by \eqref{eq:log}. Let $C$ be a disjoint union as in \eqref{eq:union}, and let $x\in C_i$. Then for every $t>0$
\begin{equation}\label{eq:hjblower}
\liminf_{T\to0}v_T(t,x)\ge\min\Big\{\inf_{U_i}g,\ \inf_{\R^d}g+\frac{\delta_i^2}4\Big\}.
\end{equation}
For the function $g(y)=\frac{\delta_i^2}4\min\{1,\,2\dist(y,C\setminus C_i)/\delta_i\}$ we have, for every $x\in C_i$,
\[
\lim_{t\to\infty}\ \lim_{T\to0}v_T(t,x)=\frac{\delta_i^2}4\ \ne\ 0=\lim_{T\to0}\ \lim_{t\to\infty}v_T(t,x).
\]
\end{theorem}

\begin{proof}
Let $x\in C_i$, $t>0$ and $\eta\in(0,1)$. By Lemma \ref{lem:confine}(i) and its proof, with $d_{C_i}(x)=0$, on the event $\{\sup_{u\le t}|\zeta_u|<\delta_i/2\}$ the process $X^x$ coincides on $[0,t]$ with $X^{i,x}$, which stays in $U_i$ on this event. Indeed, by Theorem \ref{thm:envelope} for the convex set $C_i$ we have $d_{C_i}(X^{i,x}_u)\le|\zeta_u|<\delta_i/2$ for $u\le t$ on this event. So the event $\{X^x_t\notin U_i\}$ is contained in the event $\{\sup_{u\le t}|\zeta_u|\ge\delta_i/2\}$. We estimate the probability of the second event by \eqref{eq:oumax}, used with $N=t$, the radius $\delta_i/2$, $\tr Q=d$ and $\|Q\|=1$. Here $(\delta_i/2)^2/2T=\delta_i^2/8T$. Hence, by Proposition \ref{prop:oumax},
\[
\pi_T:=\P\big(X^x_t\notin U_i\big)\le\Big(\frac{16t}{\eta^2}+1\Big)\bigg[\exp\Big(-\frac{\big((1-\eta)\frac{\delta_i}2-\sqrt{dT}\big)_+^2}{2T}\Big)+\exp\Big(\frac d2-\frac{\delta_i^2}{8T}\Big)\bigg].
\]
We bound the sum of the two exponentials by twice the larger one and take $-2T\log$. This gives
\[
-2T\log\pi_T\ge-2T\log\Big(\frac{16t}{\eta^2}+1\Big)-2T\log2+\min\Big\{\Big((1-\eta)\frac{\delta_i}2-\sqrt{dT}\Big)_+^2,\ \frac{\delta_i^2}4-dT\Big\} .
\]
As $T\to0$ the first two terms tend to zero. The minimum tends to $\min\{(1-\eta)^2\delta_i^2/4,\ \delta_i^2/4\}=(1-\eta)^2\delta_i^2/4$. Hence $\liminf_{T\to0}(-2T\log\pi_T)\ge(1-\eta)^2\delta_i^2/4$.

Now we split the expectation in \eqref{eq:log} according to whether $X^x_t\in U_i$. On the event $\{X^x_t\in U_i\}$ we have $g(X^x_t)\ge\inf_{U_i}g$. On the complement, which has probability $\pi_T$, we have $g(X^x_t)\ge\inf g$. Thus
\begin{align*}
\E\exp\Big(-\frac{g(X^x_t)}{2T}\Big)&\le\exp\Big(-\frac{\inf_{U_i}g}{2T}\Big)+\exp\Big(-\frac{\inf g}{2T}\Big)\,\pi_T\\
&\le2\max\Big\{\exp\Big(-\frac{\inf_{U_i}g}{2T}\Big),\ \exp\Big(-\frac{\inf g}{2T}\Big)\pi_T\Big\}.
\end{align*}
In the second line we used that a sum of two nonnegative numbers is at most twice the larger one. Taking $-2T\log$, which reverses the inequality and turns the maximum into a minimum, we get
\[
v_T(t,x)\ge-2T\log2+\min\Big\{\inf_{U_i}g,\ \inf g-2T\log\pi_T\Big\} .
\]
The lower limit of a minimum of two functions is the minimum of their lower limits. Letting $T\to0$,
\[
\liminf_{T\to0}v_T(t,x)\ge\min\Big\{\inf_{U_i}g,\ \inf g+(1-\eta)^2\frac{\delta_i^2}4\Big\},
\]
and $\eta\to0$ gives \eqref{eq:hjblower}.

For the example, the function $g$ is Lipschitz, $0\le g\le\delta_i^2/4$ and $g=0$ on $C\setminus C_i$, which is not empty since $m\ge2$. The Lipschitz property holds because the distance to the set $C\setminus C_i$ is nonexpansive. If $y\in U_i$, then for every $c\in C_j$, $j\ne i$, we have $|y-c|\ge\delta_i-d_{C_i}(y)>\delta_i/2$, so $\dist(y,C\setminus C_i)\ge\delta_i/2$ and $g(y)=\delta_i^2/4$. Here the first inequality is the triangle inequality,
\[
|y-c|\ge|P_iy-c|-|y-P_iy|\ge\delta_i-d_{C_i}(y) ,
\]
and the second holds by the definition of $U_i$. Hence $\inf_{U_i}g=\delta_i^2/4$ and $\inf g=0$. So \eqref{eq:hjblower} gives $\liminf_{T\to0}v_T(t,x)\ge\delta_i^2/4$. Since also $v_T\le\sup g=\delta_i^2/4$ by \eqref{eq:log}, we get $\lim_{T\to0}v_T(t,x)=\delta_i^2/4$ for every $t>0$. The bound $v_T\le\sup g$ holds because $\exp(-g/2T)\ge\exp(-\sup g/2T)$. The limit $\delta_i^2/4$ does not depend on $t$, so the first iterated limit is $\delta_i^2/4$. On the other hand, by Theorem \ref{thm:hjb}(a) the other iterated limit is $g_C$. We have $g_C=0$, because $g+d_C^2$ vanishes on $C\setminus C_i$ and is nonnegative everywhere. This completes the proof.
\end{proof}

\begin{corollary}\label{cor:hjbconvex}
Let $C\subset\R^d$ be a finite disjoint union of nonempty compact convex sets. The following are equivalent.
\begin{enumerate}[label=\textup{(\roman*)}]
\item $C$ is convex.
\item For every bounded Lipschitz function $g$ and every $x\in C$, both iterated limits of $v_T(t,x)$ exist and
\[
\lim_{t\to\infty}\ \lim_{T\to0}v_T(t,x)=\lim_{T\to0}\ \lim_{t\to\infty}v_T(t,x).
\]
\end{enumerate}
\end{corollary}

\begin{proof}
A finite disjoint union of $m$ nonempty compact convex sets is convex if and only if $m=1$, since a convex set is connected. Indeed, for $m\ge2$ the set $C$ is the union of the two disjoint nonempty closed sets $C_1$ and $C_2\cup\dots\cup C_m$, so it is not connected. If $m=1$, (ii) holds by Theorem \ref{thm:hjb}(b). If $m\ge2$, (ii) fails by Theorem \ref{thm:hjbdisc}. For the function $g$ of that example and $x\in C_i$ the two iterated limits exist and are different.
\end{proof}

\begin{remark}
We compare Corollary \ref{cor:commute} and Corollary \ref{cor:hjbconvex}. For the law of the state the two limits do not commute even for a convex set, unless the initial law is special. For the value function they commute for every convex set and every terminal cost. The proof of part (b) of \cite[Theorem~6.3]{DarapaneniConvex} shows the reason. Inside a convex set the drift vanishes. Over a long horizon a control can cross the set with an energy which tends to zero, so in the limit the value does not depend on the initial state. Without control the state needs a time of order $1/T$ to cross the set. When the feasible set is disconnected, crossing from $C_i$ to another component costs at least the energy $\delta_i^2/4$, which by \cite[Proposition~6.1]{DarapaneniConvex} is the energy needed to reach the boundary of the basin. This energy is exactly the difference of the two limits in Theorem \ref{thm:hjbdisc}. It is also twice the constant $\delta_i^2/8$ of Theorem \ref{thm:exit}. This agrees with Lemma \ref{lem:tube}, where a path of energy $J$ is followed with a probability of order $\exp(-J/2T)$.
\end{remark}

\section{Open questions}\label{sec:open}

Some questions remain open.

\begin{enumerate}[label=\arabic*.,leftmargin=2em]
\item \emph{The threshold for convergence to equilibrium with four and more components.} Prove that the conclusion of Theorem \ref{thm:equil} holds for every $c>\Delta_*^2/8$ and every $m$. See Remark \ref{rem:equil}. For this one has to describe the order in which the wells are visited, as in the hierarchy of cycles of \cite{FreidlinWentzell2012}, for flat wells of different dimensions.
\item \emph{Prefactors.} Find the asymptotic law of the transition time between two components, with the prefactor. By Theorem \ref{thm:gibbs}, up to an exponentially small error, the integral of $e^{-d_C^2/2T}$ over a basin is a polynomial in $\sqrt T$ whose coefficients are multiples of intrinsic volumes. So the prefactor should be a power of $T$ which depends on the dimensions of the components and on the local shape of the ridge at its lowest points. Compare \cite{BerglundGentz2010,AvelinJulinViitasaari2023}.
\item \emph{Sets which are not convex.} In the Euclidean case, the local Steiner formula holds for small radii for sets of positive reach \cite{Federer1959}. This suggests versions of Theorem \ref{thm:gibbs} and Theorem \ref{thm:meta} for unions of such sets.
\end{enumerate}

\end{document}